\documentclass[11pt]{amsart}
\usepackage{amscd} 
\usepackage[dvips]{graphicx, color}
\usepackage{amssymb}
\usepackage{psfrag}
\usepackage{enumerate}

\numberwithin{equation}{section} 

\newcommand{\bp}{p}
\newcommand{\bu}{u}
\newcommand{\bv}{v}
\newcommand{\bH}{{\mathbf H}}
\newcommand{\bV}{{\mathbf V}}
\newcommand{\bE}{{\mathbf E}}

\newcommand{\bC}{\mathbb{C}}
\newcommand{\bN}{\mathbb{N}}
\newcommand{\bR}{\mathbb{R}}
\newcommand{\bZ}{\mathbb{Z}}

\newcommand{\Ht}{{\mathrm{ht}}}
\newcommand{\HV}{{\mathrm{HV}}}
\newcommand{\id}{{\mathrm{id}}}
\newcommand{\odd}{{\mathrm{odd}}}

\newcommand{\sgn}{{\mathrm{sgn}}}
\newcommand{\Tab}{{\mathrm{Tab}}}
\newcommand{\HVTab}{{\mathrm{Tab}_{\scriptscriptstyle \mathrm HV}}}
\newcommand{\HVP}{P_{\scriptscriptstyle \mathrm HV}}
\newcommand{\Rep}{{\mathrm{Rep}}}

\newcommand{\cA}{\mathcal{A}}

\newcommand{\cH}{\mathcal{H}}
\newcommand{\cLU}{\mathcal{LU}}
\newcommand{\cP}{\mathcal{P}}
\newcommand{\cTv}{\mathcal{T}_{\mathrm v}}
\newcommand{\cU}{\mathcal{U}}
\newcommand{\cW}{\mathcal{W}}
\newcommand{\cZ}{\mathcal{Z}}

\newcommand{\fg}{\mathfrak{g}}
\newcommand{\hg}{\hat{\fg}}
\newcommand{\fP}{\mathfrak{P}}
\newcommand{\fS}{\mathfrak{S}}

\newcommand{\Uqhg}{U_q(\hg )}

\newcommand{\lprod}{\overset{\leftarrow}{\prod }}
\newcommand{\rprod}{\overset{\rightarrow}{\prod }}

\newcommand{\Position}[1]{\mathrm{hp}( #1 )}

\newcommand{\path}[3]{#1 \overset{#2}{\to} #3}
\newcommand{\I}{\mathrm{I}}
\newcommand{\II}{\mathrm{II}}
\newcommand{\one}{(1)}
\newcommand{\two}{(2)}
\newcommand{\three}{(3)}
\newcommand{\four}{(4)}
\newcommand{\five}{(5)}
\newcommand{\six}{(6)}
\newcommand{\seven}{(7)}

\usepackage{amsthm}
\newtheorem{thm}{Theorem}[section]
\newtheorem{lem}[thm]{Lemma}
\newtheorem{prop}[thm]{Proposition}

\theoremstyle{definition}
\newtheorem{defn}[thm]{Definition}
\newtheorem{exmp}[thm]{Example}
\theoremstyle{remark}
\newtheorem{rem}[thm]{Remark}

\begin{document}
\title[Jacobi-Trudi determinant of type $D_n$]{
Paths and tableaux descriptions of \\
Jacobi-Trudi determinant associated with\\ 
quantum affine algebra of type $D_n$}

\author[W.~Nakai]{}
%
\author[T.~Nakanishi]{}

\maketitle

\begin{center}
\vspace{-8pt}
{\sc Wakako Nakai\footnote{e-mail: {\tt m99013c@math.nagoya-u.ac.jp}}
and Tomoki Nakanishi\footnote{e-mail: {\tt nakanisi@math.nagoya-u.ac.jp}}}\\
{\scriptsize {\it Graduate School of Mathematics, 
Nagoya University,
Nagoya 464-8602, Japan}}
\end{center}
\vspace{13pt}

{\small 
{\sc Abstract}.
We study the Jacobi-Trudi-type determinant which is conjectured to be
the $q$-character of a certain, in many cases irreducible, finite-dimensional
representation of the quantum affine algebra of type $D_n$.
Unlike the $A_n$ and $B_n$ cases, a simple application of the Gessel-Viennot
path method does not yield an expression of the
determinant by a positive sum over a set of tuples of paths. 
However, applying an additional involution
and a deformation of paths, we obtain an expression by a positive sum over
a set of tuples of paths,
which is naturally translated into the one over a set of tableaux
on a skew diagram.
}
\vspace{11pt}


\section{Introduction}
Let $\fg$ be the simple Lie algebra over $\bC$, and let 
$\hg$ be the corresponding untwisted affine Lie algebra. 
Let $\Uqhg$ be the quantum affine algebra, namely, 
the quantized universal enveloping algebra 
of $\hg$ \cite{D1,J}. 
In order to investigate the 
finite-dimensional representations of $\Uqhg$
\cite{D2,CP}, 
an injective ring homomorphism
\begin{equation}
\chi_q: \Rep \, \Uqhg \to \bZ[Y_{i,a}^{\pm 1}]_{i=1, \dots, n; a \in \bC^{\times}},
\end{equation}
called the {\it  $q$-character} of $\Uqhg$, was introduced
and studied in \cite{FR1,FM}, 
where $\Rep \, \Uqhg$ is the Grothendieck ring of the category of the 
finite-dimensional representations of $\Uqhg$. 
The $q$-character contains the essential data of 
each representation $V$. 
So far, however, the explicit description of 
$\chi_q(V)$ is available only for
a limited type of representations (e.g., the fundamental representations)
\cite{FM,CM},
and the description  for 
general $V$ is an open problem.
See \cite{N,H} for related results.

In our previous work \cite{NN}, 
for $\fg = A_n$, $B_n$, $C_n$, and $D_n$,
we conjecture
that the $q$-characters
of a certain family of, in many cases irreducible, finite-dimensional
representations  are given by
the determinant form $\chi_{\lambda/\mu,a}$,
where $\lambda/\mu$ is a skew diagram and $a$ is a complex
parameter.
We call $\chi_{\lambda/\mu,a}$ the Jacobi-Trudi determinant.
For $A_n$ and $B_n$, this is a reinterpretation
of the conjecture for the spectra
of the transfer matrices of the vertex models
associated with the corresponding representations \cite{BR,KOS}.
See also \cite{KS} for related results for $C_n$ and $D_n$.
Let us briefly summarize the result of \cite{NN}.
Following the standard Gessel-Viennot method \cite{GV},
we represent the Jacobi-Trudi determinant by {\it paths},
and apply an involution for intersecting paths.
For $A_n$ and $B_n$, this immediately reproduces
the known tableaux descriptions of $\chi_{\lambda/\mu,a}$
by \cite{BR,KOS}.
Here,  by {\it tableaux description\/} we mean
an expression  of 
$\chi_{\lambda/\mu,a}$ by a positive sum over 
a certain set of tableaux on $\lambda/\mu$.
For $A_n$, the relevant tableaux are nothing but the
semistandard tableaux as  the usual character for $\fg=A_n$.
For $B_n$, the tableaux are given by 
the `horizontal' and `vertical' rules
similar to the ones for the semistandard tableaux.
In contrast, we find that the tableaux description of 
$\chi_{\lambda/\mu,a}$ for $C_n$
is not as simple as the former cases. 
The main difference is that
a simple application of the Gessel-Viennot method does not
yield an expression  of $\chi_{\lambda/\mu,a}$ by 
a positive sum. 
Nevertheless, in some special cases (i.e.,
a skew diagram $\lambda/\mu$ of at most two columns or
of at most three rows),
one can further work out the cancellation of the remaining
negative contribution, and obtain a tableaux description of
$\chi_{\lambda/\mu,a}$.
Besides the horizontal and vertical rules,
we have an additional rule, which we call the {\it extra rule},
due to the above process.

In this paper, we consider the same problem
for $\chi_{\lambda/\mu,a}$ for $D_n$,
where the situation is quite parallel to $C_n$.
By extending the idea of \cite{NN},
we now successfully obtain a tableaux
description of $\chi_{\lambda/\mu,a}$ 
for a general skew diagram $\lambda/\mu$.
The resulting tableaux description 
shows nice compatibility
with the proposed algorithm to generate the $q$-character
by \cite{FM}, and it is expected to be useful to study the
$q$-characters further.
We also hope that our tableaux will be useful to
parameterize the much-awaited crystal basis for
the Kirillov-Reshetikhin representations \cite{KN,OSS},
where $\lambda/\mu$ is a rectangular shape.
To support it, for a two-row rectangular diagram $\lambda/\mu$,
our tableaux agree with the ones
for the proposed crystal graph by \cite{SS}.
Meanwhile, our tableau rule is rather different from the 
one for the non-quantum case \cite{FK} due to the 
different nature of the determinant and the generating function. 
The method herein is also applicable to a general skew diagram
$\lambda/\mu$ for  $C_n$,
and it will be reported in a separate publication \cite{NN2}.

Now let us explain the organization and the main idea
of the paper.

In Section \ref{sec:JT-det}, following \cite{NN}, 
we define the Jacobi-Trudi determinant $\chi_{\lambda/\mu,a}$
for $D_n$. The procedure to derive the tableaux
description of $\chi_{\lambda/\mu,a}$
consists of three steps.

In Section \ref{sec:GV} we do the first step.
Here we apply the standard method by  \cite{GV}
for the determinant $\chi_{\lambda/\mu,a}$.
Namely, first, we introduce lattice paths,
and express  $\chi_{\lambda/\mu, a}$
as a sum
over a set of $l$-tuples of paths $p=(p_i)$ with fixed end points.
Secondly, we define the weight-preserving, sign-reversing
involution $\iota_1$ (the {\it first involution})
so that for an intersecting tuple of paths $p$
the contributions  from $p$ and $\iota_1(\bp)$
cancel each other in the sum.
Unlike $A_n$ and $B_n$, however, this involution cannot be
defined on the entire set of the intersecting tuples of paths,
and the resulting expression for $\chi_{\lambda/\mu, a}$
(the first sum, Proposition \ref{prop:first-involution}) 
still includes negative terms.

In Section \ref{sec:second-involution} we do the second step.
Extending the idea of \cite{NN} for $C_n$,
we define {\it another} weight-preserving, sign-reversing involution
$\iota_2$ (the {\it second involution}).
Then, the resulting expression (the second sum,
Theorem \ref{thm:positive-sum})
no longer includes 
negative terms.
However, the contribution  from the tuples of paths with 
`transposed' pairs still remains,
and one cannot naturally translate such a tuple of paths
into a tableau on the skew diagram $\lambda/\mu$.

In Sections \ref{sec:tableau-description} and \ref{sec:folding-map}, 
we do the last step.
In Section \ref{sec:tableau-description}, we claim the existence
of a weight-preserving deformation $\phi$ of
the paths (the {\it folding map}),
where $\phi$ `resolves' transposed pairs by folding.
The resulting expression (the third sum, Theorem \ref{thm:tilde-positive-sum})
is now naturally  translated into the tableaux description
whose tableaux are determined by the
horizontal, vertical, and extra rules (Theorems \ref{thm:tableau-description} 
and \ref{thm:configuration}).
The explicit list of the extra rule has wide
variety,
and examples are given for $\lambda/\mu$ with at most two columns
or at most three rows.
The construction of the folding map $\phi$ is 
the most technical part of the work. 
We provide the details in 
Section \ref{sec:folding-map}.

We remark that
while the explicit list of the
extra rule for tableaux looks
rather complicated and disordered,
it is a simple and easily recognizable 
graphical rule in the path language.
Therefore, the paths description (especially, the third sum)
may be  as important as the tableaux description for applications.

\section{The Jacobi-Trudi determinant of type $D_n$}\label{sec:JT-det}

In this section, we define the Jacobi-Trudi determinant 
$\chi_{\lambda/\mu,a}$, following \cite{NN}. 
See \cite{NN} for more information. 

A {\it partition} is a sequence of weakly decreasing 
non-negative integers $\lambda=(\lambda_1, \lambda_2, \dots)$
with finitely many non-zero terms $\lambda_1 \ge \lambda_2 \ge \cdots \ge \lambda_l>0$. 
The {\it length} $l(\lambda)$ of $\lambda$ is 
the number of the non-zero integers. 
The {\it conjugate} of $\lambda$ 
is denoted by $\lambda'=(\lambda'_1,\lambda'_2, \dots)$. 
As usual, we identify a partition $\lambda$ with a {\it Young diagram}
$\lambda = \{ (i,j)\in \bN^2 \mid 1 \le j \le \lambda_i \}$, 
and also identify a pair of partitions such that $\lambda_i \ge \mu_i$ 
for any $i$, with a {\it skew diagram}
$\lambda/\mu = \{ (i,j) \in \bN^2 \mid \mu_i+1 \le j \le \lambda_i\}$. 

Let 
\begin{equation}\label{eq:entries}
I=\{1, 2, \dots, n, \overline{n}, \dots,\overline{2},  \overline{1} \}. 
\end{equation}
Let $\cZ$ be the commutative ring over $\bZ$ generated by 
$z_{i,a}$'s, $i \in I$, $a \in \bC$, with the following generating relations: 
\begin{gather}\label{eq:relations}
z_{i,a}z_{\overline{i},a-2n+2i} = z_{i-1,a}z_{\overline{i-1},a-2n+2i}
\quad (i=1, \dots, n), \quad 
z_{0,a}=z_{\overline{0},a}=1.
\end{gather}

Let $\bZ[[X]]$ be the formal power series ring 
over $\bZ$ with variable $X$. 
Let $\cA$ be the {\it non-commutative} ring
generated by $\cZ$ and $\bZ[[X]]$ with relations 
$$Xz_{i,a}=z_{i,a-2}X, \qquad i \in I, a \in \bC. $$
For any $a \in \bC$, we define $E_a(z,X)$, $H_a(z,X)\in \cA$ 
as 
\begin{align}
E_a(z,X) & := \label{eq:E}
\big\{ \prod_{\scriptscriptstyle 1 \le k \le n}^{\rightarrow} (1 + z_{k,a}X) \big\}
(1-z_{\overline{n},a}Xz_{n,a}X)^{-1}
\big\{ \prod_{\scriptscriptstyle 1 \le k \le n}^{\leftarrow} 
(1 + z_{\overline{k},a}X) \big\}, \\
H_a(z,X) & :=
\big\{ \prod_{\scriptscriptstyle 1 \le k \le n}^{\rightarrow} 
(1 - z_{\overline{k},a}X)^{-1} \big\}
(1 - z_{\overline{n},a}Xz_{n,a}X)
\big\{\prod_{\scriptscriptstyle 1 \le k \le n}^{\leftarrow} (1 - z_{k,a}X)^{-1} \big\}, 
\end{align} 
where $\rprod _{\scriptscriptstyle 1 \le k \le n}A_k=A_1\dots A_n$
and $\lprod _{\scriptscriptstyle 1 \le k \le n}A_k=A_n\dots A_1$.
Then we have
\begin{equation}\label{eq:HE}
H_a(z, X)E_a(z,-X) = E_a(z,-X)H_a(z,X)=1.
\end{equation}

For any $i \in \bZ$ and $a \in \bC$, 
we define $e_{i,a}$, $h_{i,a}\in \cZ$ as 
$$
E_a(z,X) = \sum_{i=0}^{\infty}e_{i,a}X^i, \qquad 
H_a(z,X) = \sum_{i=0}^{\infty}h_{i,a}X^i, $$
with $e_{i,a}=h_{i,a}=0$ for $i < 0$. 

Due to relation \eqref{eq:HE}, we have \cite[Eq.\ (2.9)]{M}
\begin{equation}\label{eq:determinant}
\det (h_{\lambda_i-\mu_j-i+j, a+2(\lambda_i - i)})_{1 \le i,j \le l}
= \det(e_{\lambda'_i-\mu'_j-i+j, a-2(\mu'_j-j+1)})_{1 \le i,j \le l'}
\end{equation}
for any pair of partitions $(\lambda, \mu)$, where $l$ and $l'$ are 
any non-negative integers such that $l \ge l(\lambda), l(\mu)$
and $l' \ge l(\lambda'), l(\mu')$. 
For any skew diagram $\lambda/\mu$, 
let $\chi_{\lambda/\mu,a}$ denote the determinant on the left- or 
right-hand side of \eqref{eq:determinant}. 
We call it the {\it Jacobi-Trudi determinant} associated with  
the quantum affine algebra $\Uqhg$ of type $D_n$. 

Let $d(\lambda/\mu):=\max \{ \lambda'_i-\mu'_i\}$ be the {\it depth} 
of $\lambda/\mu$. 
We conjecture that, if $d(\lambda/\mu) \le n$, 
the determinant $\chi_{\lambda/\mu, a}$ 
is the {\it $q$-character} for a certain finite-dimensional representations 
$V$ of quantum affine algebras. 
We further expect that $\chi_{\lambda/\mu, a}$ is 
the $q$-character for an irreducible $V$, if $d(\lambda/\mu) \le n-1$ and 
$\lambda/\mu$ is connected \cite{NN}. 

\begin{rem}
The above conjecture and the ones for types $B_n$ and $C_n$ in \cite{NN}
tell that the irreducible character of $\Uqhg$, corresponding 
to a connected skew diagram, is always expressed by the 
same determinant \eqref{eq:determinant} regardless of the 
type of the algebra. This is a remarkable contrast to 
the non-quantum case \cite{KT}. 
For example, the tensor product of two first fundamental 
modules of $\fg$ has two irreducible submodules for 
type $A_n$ and three ones for type $B_n$, $C_n$, or $D_n$. 
On the other hand, 
under the appropriate choice of the values for the 
spectral parameters, the tensor product of two first fundamental representations 
of $\Uqhg$ has exactly two irreducible subquotients, 
one of which corresponds to two-by-one rectangular diagram 
and the other of which corresponds to one-by-two rectangular 
diagram, regardless of the type of the algebra. 
In fact, this is the simplest example of the conjecture. 
\end{rem}

\section{Gessel-Viennot paths and the first involution}\label{sec:GV}

Following \cite{NN}, 
let us apply the method by \cite{GV} to the determinant 
$\chi_{\lambda/\mu,a}$ in \eqref{eq:determinant} 
and the generating function $E_a(z, X)$ in \eqref{eq:E}. 

\begin{figure}
\begin{psfrags}
\psfrag{x}{$\scriptstyle x$}
\psfrag{y}{$\scriptstyle y$}
\psfrag{s_1}{$s_1$}
\psfrag{s_2}{$s_2$}
\psfrag{s_3}{$s_3$}
\psfrag{L_a^1}{$L_a^1(s_i)$}
\psfrag{L_a^2}{$L_a^2(s_i)$}
\psfrag{height}{height}
\psfrag{-n}{$\scriptstyle -n$}
\psfrag{-n+1}{$\scriptstyle -n+1$}
\psfrag{-1}{$\scriptstyle -1$}
\psfrag{0}{$\scriptstyle 0$}
\psfrag{1}{$\scriptstyle 1$}
\psfrag{2}{$\scriptstyle 2$}
\psfrag{n-1}{$\scriptstyle n-1$}
\psfrag{n}{$\scriptstyle n$}
\psfrag{bn, n}{$\scriptstyle \overline{n}, n$}
\psfrag{b1}{$\scriptstyle \overline{1}$}
\psfrag{b2}{$\scriptstyle \overline{2}$}
\psfrag{bn1}{$\scriptstyle \overline{n-1}$}
\psfrag{bn}{$\scriptstyle \overline{n}$}
\psfrag{d}{$\cdots$}
\psfrag{v}{$\vdots$}
\psfrag{a}{$\scriptstyle a$}
\psfrag{a-2}{$\scriptstyle a-2$}
\psfrag{a-4}{$\scriptstyle a-4$}
\psfrag{N}{$\scriptstyle N$}
\psfrag{W}{$\scriptstyle W$}
\psfrag{S}{$\scriptstyle S$}
\psfrag{E}{$\scriptstyle E$}
\psfrag{O}{$\scriptscriptstyle (0,0)$}
\includegraphics[width=12cm, clip]{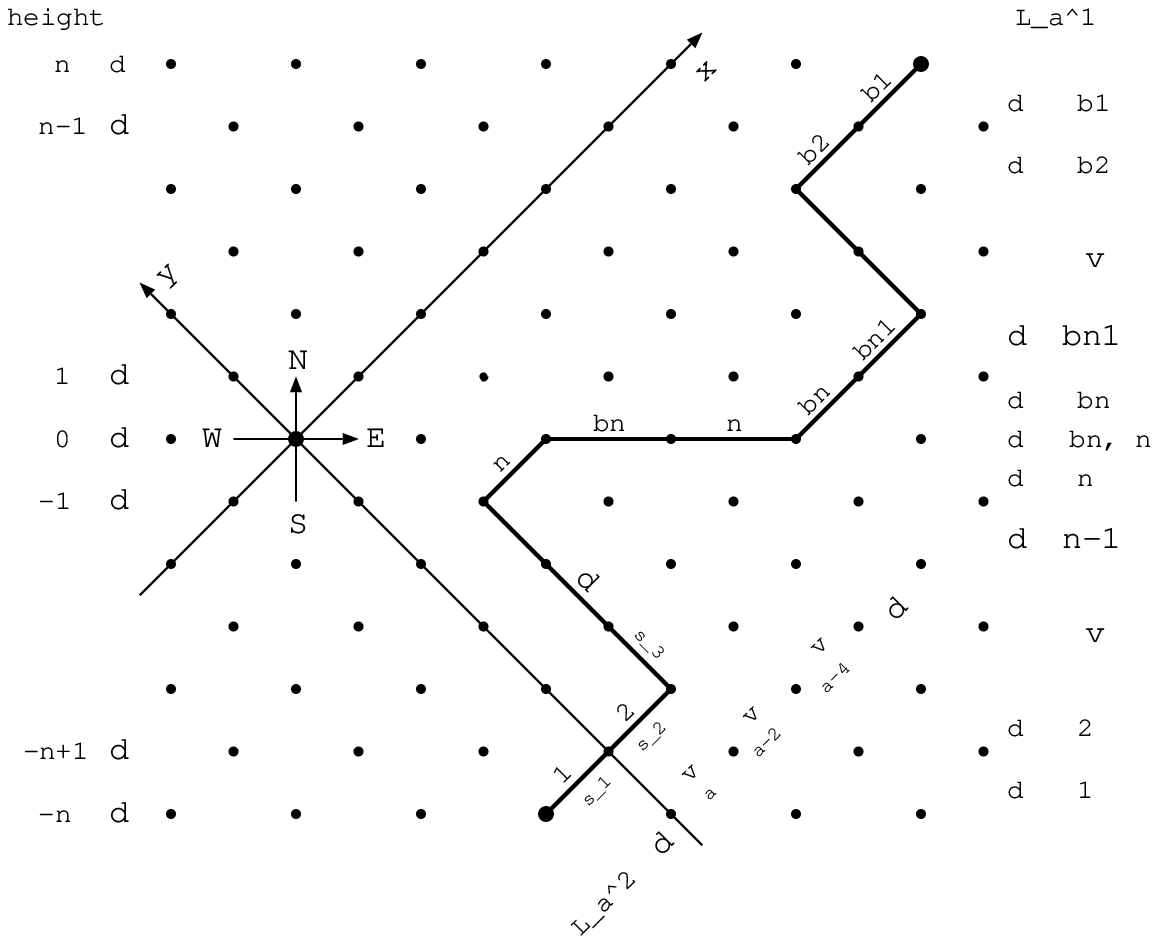}
\end{psfrags}
\caption{An example of a path of type $D_n$ and its $e$-labeling.}\label{fig:path}
\end{figure}

Consider the lattice $\bZ^2$ 
and rotate it by $45^{\circ}$ as in Figure \ref{fig:path}. 
An {\it E-step} $s$ is a step between two points in 
the lattice of length $\sqrt{2}$ in east direction. 
Similarly, an {\it NE-step} (resp.~an {\it NW-step}) is 
a step between two points in the lattice 
of unit length in northeast direction 
(resp.~northwest direction).  
For any point $(x,y)\in \bR^2$,  
we define the {\it height} as $\Ht(x,y):=x+y$, 
and the {\it horizontal position} as $\Position{x,y}:=\frac{1}{2}(x-y)$. 
Due to \eqref{eq:E}, we define a {\it path $p$} (of type $D_n$)  
as a sequence of consecutive steps 
$(s_1, s_2, \dots )$ which satisfies the following conditions: 
\begin{enumerate} 
\item 
It starts from a point $u$ at height $-n$ 
and ends at a point $v$ 
at height $n$. 
\item 
Each step $s_i$ is an NE-, NW-, 
or E-step. 
\item  \label{item:D-path-three}
The E-steps occur only at   
height $0$, 
and the number of E-steps is even. 
\end{enumerate}
We also write $p$ as $\path{u}{p}{v}$. 
See Figure \ref{fig:path} for an example. 

Let $\cP$ be the set of all paths of type $D_n$. 
For any $p \in \cP$, set 
\begin{equation}\label{eq:label-steps}
\begin{aligned}
E(p)& := \{s \in p \mid 
\text{$s$ is an NE- or E-step} \}, \\
E_0(p)& := \{s \in p \mid 
\text{$s$ is an E-step} \} \subset E(p). 
\end{aligned}
\end{equation}
If $E_0(p)=\{s_j, s_{j+1}, \dots, s_{j+2k-1} \}$, then 
let 
$$E_0^1(p):=\{s_{j+1}, s_{j+3}, \dots, s_{j+2k-1} \}
\subset E_0(p).$$ 
Fix $a \in \bC$. The {\it $e$-labeling} (of type $D_n$) associated 
with $a$ for a path $p \in \cP$ is the pair of maps 
$L_a = (L^1_a, L^2_a)$ on $E(p)$ defined as follows: 
Suppose that a step $s\in E(p)$ starts at a point $w=(x,y)$, and 
let $m:=\Ht(w)$. Then, we set  
\begin{equation}\label{eq:e-label}
\begin{aligned}
L_a^1(s) & =
\begin{cases}
n+1+m, & \text{if $m < 0$}, \\
n, & \text{if $m=0$ and $s \in E_0^1(p)$}, \\
\overline{n-m}, & \text{otherwise}, \\
\end{cases}
\\
L_a^2(s) &= a - 2x. 
\end{aligned}
\end{equation}
See Figure \ref{fig:path}. 

Now we define the {\it weight} of $p \in \cP$ as 
$$
z_a^p := \prod_{s \in E(p)} z_{L_a^1(s), L_a^2(s)} \in \cZ. 
$$
By the definition of $E_a(z,X)$ in \eqref{eq:E}, 
for any $k\in \bZ$, we have  
\begin{equation}\label{eq:e}
e_{r, a-2k}(z) = \sum_pz_a^p, 
\end{equation}
where the sum runs over all $p \in \cP$ 
such that $(k, -n-k) \overset{p}{\to} (k+r, n-k-r)$. 

For any $l$-tuples of distinct points   
$\bu = (u_1, \dots, u_l)$ of height $-n$ and  
$\bv = (v_1, \dots, v_l)$ of height $n$,  
and any permutation $\sigma \in \fS_l$, let 
\begin{gather*}
\fP(\sigma; \bu, \bv):= \{ \bp = (p_1, \dots, p_l) \mid 
p_i\in \cP,\ \path{u_i}{p_i}{v_{\sigma(i)}} \text{ for } i=1, \dots, l \}, 
\end{gather*}
and set 
\begin{gather*}
\fP(\bu, \bv) := \bigsqcup_{\sigma \in \fS_l}\fP(\sigma; \bu, \bv). 
\end{gather*}
We define the {\it weight} $z_a^{\bp}$ and the 
{\it sign} $(-1)^{\bp}$ of 
$\bp \in \fP(\bu, \bv)$ 
as 
\begin{equation}\label{eq:weight}
z_a^{\bp} := \prod_{i=1}^lz_a^{p_i}, \qquad 
(-1)^{\bp}:= \sgn \, \sigma \quad \text{if $\bp \in \fP(\sigma; \bu, \bv)$}. 
\end{equation}

For any skew diagram $\lambda/\mu$, set $l=\lambda_1$, and 
\begin{align*}
\bu_{\mu}:=(u_1, \dots, u_{l}), \qquad 
& u_i:=(\mu'_i + 1 -i, -n - \mu'_i - 1 + i),  \\
\bv_{\lambda}:=(v_1, \dots, v_{l}), \qquad 
& v_i:=(\lambda'_i + 1 -i, n - \lambda'_i - 1 + i). 
\end{align*}
Then, due to \eqref{eq:e}, the determinant \eqref{eq:determinant}
can be written as 
\begin{equation}\label{eq:GV-path-description}
\chi_{\lambda/\mu,a}
= \sum_{\bp \in \fP(\bu_{\mu}, \bv_{\lambda})} (-1)^{\bp}z_a^{\bp}. 
\end{equation}

In the $A_n$ case, one can define a natural weight-preserving, 
sign-reversing involution on the set of all the tuples 
$\bp$ which have some intersecting pair $(p_i, p_j)$. 
However, this does not hold for $D_n$ because of  
Condition \eqref{item:D-path-three} of the definition 
of a path of type $D_n$. Therefore, 
as in the cases of types $B_n$ and $C_n$ \cite{NN}, 
we introduce the following notion: 

\begin{defn}
We say that 
an intersecting pair $(p_i,p_j)$ of paths is  
{\it specially intersecting} if it satisfies the following conditions: 
\begin{enumerate}
\item The intersection of $p_i$ and $p_j$ occurs only at height $0$. 
\item $p_i(0)-p_j(0)$ is odd, where 
$p_i(0)$ is the horizontal position of 
the leftmost point on $p_i$ at height $0$. 
\end{enumerate}
Otherwise, we say that an intersecting pair $(p_i,p_j)$ is 
{\it ordinarily intersecting}. 
\end{defn}
As in the cases of types $B_n$ and $C_n$ \cite{NN}, 
we can define a weight-preserving, sign-reversing involution 
$\iota_1$ on the set of all the tuples $\bp \in \fP(\bu_{\mu}, \bv_{\lambda})$ 
which have some ordinarily intersecting pair $(p_i, p_j)$. 
Therefore, we have 
\begin{prop}\label{prop:first-involution}
For any skew diagram $\lambda/\mu$, 
\begin{equation}\label{eq:first-sum}
\chi_{\lambda/\mu, a}= \sum_{\bp \in P_1(\lambda/\mu)}(-1)^{\bp}z_a^{\bp},
\end{equation}
where $P_1(\lambda/\mu)$ is the set of all 
$\bp \in \fP(\bu_{\mu}, \bv_{\lambda})$ which do not have any ordinarily 
intersecting pair $(p_i,p_j)$ of paths. 
\end{prop}

For $B_n$, the sum \eqref{eq:first-sum} is a positive sum 
because no $\bp \in P_1(\lambda/\mu)$ has a `transposed' pair 
$(p_i, p_j)$. But, this is not so for $C_n$ and $D_n$. 

\section{The second involution}\label{sec:second-involution}
In this section, we define another weight-preserving involution, 
the {\it second involution}. This is defined by using 
the paths deformations called {\it expansion} and 
{\it folding}. 
As a result, the second involution 
cancels all the negative contributions in  
\eqref{eq:first-sum}, and we obtain 
an  expression by a positive sum, see 
\eqref{eq:positive-sum}. 

\subsection{Expansion and folding}\label{sec:expansion-folding}

Let
\begin{align*}
S_+ & :=\{ (x,y) \in \bR^2 \mid 0 \le \Ht(x,y)\le n \}, \\
S_- & :=\{ (x,y) \in \bR^2 \mid -n \le \Ht(x,y)\le 0 \}. 
\end{align*}

For any $w=(x,y) \in S_+$, define $w^*\in S_-$ by 
$$w^*:=(-y+1,-x-1).$$
Then we have $\Ht(w^*)=-\Ht(w)$, $\Position{w^*}=\Position{w}+1$. 
Conversely, we define $(w^*)^*=w$, and we call the correspondence
\begin{equation}\label{eq:dual-map}
S_+\leftrightarrow S_-, \qquad w \leftrightarrow w^*
\end{equation}
the {\it dual map}. 

\begin{defn}
A {\it lower path} $\alpha$ (of type $D_n$)
is a sequence of consecutive steps in $S_-$ which starts at a 
point of height $-n$ and
ends at a point of height $0$, and each step is 
an NE- or NW-step. 
Similarly, an {\it upper path} $\beta$ (of type $D_n$) 
is a sequence of consecutive steps in $S_+$ which 
starts at a point of height $0$ and
ends at a point of height $n$, and each step is
an NE- or NW-step. 
\end{defn}
For any lower path $\alpha$ and an upper path $\beta$, 
let $\alpha(r)$ and $\beta(r)$ be the horizontal positions 
of $\alpha$ and $\beta$ at height $r$, respectively. 
We define an upper path $\alpha^*$ and a lower path $\beta^*$ by 
$$\alpha^*(r)=\alpha(-r)-1, \qquad 
\beta^*(-r)=\beta(r)+1, \quad (0 \le r \le n)$$
and call them the {\it duals} of $\alpha$, $\beta$. 

Let 
$$(\alpha; \beta):=(\alpha_1, \dots, \alpha_l; \beta_1, \dots, \beta_l)$$
be a pair of an $l$-tuple $\alpha$ of lower paths 
and an $l$-tuple $\beta$ of upper paths. 
We say that $(\alpha; \beta)$ is {\it nonintersecting} 
if $(\alpha_i, \alpha_j)$ is not intersecting, and so is 
$(\beta_i, \beta_j)$ for any $i,j$. 

{}From now on, 
let $\lambda/\mu$ be a skew diagram, and we set $l=\lambda_1$. 
Let 
\begin{multline*}
\cH(\lambda/\mu):= \\
\Bigg\{ (\alpha; \beta)=(\alpha_1, \dots, \alpha_l; \beta_1, \dots, \beta_l) 
\Bigg| 
\begin{array}{l}
(\alpha; \beta) \text{ is nonintersecting}, \\
\alpha_i(-n)=\frac{n}{2}+\mu'_i+1-i, \\
\beta_i(n)=-\frac{n}{2}+\lambda'_i+1-i
\end{array}
\Bigg\}. 
\end{multline*}

For any skew diagram $\lambda/\mu$, 
we call the following condition the {\it positivity condition}: 
\begin{equation}\label{eq:positivity}
\lambda'_{i+1}-\mu'_i \le n, \qquad i=1, \dots, l-1. 
\end{equation}
We call 
this the `positivity condition', 
because \eqref{eq:positivity} guarantees that 
$\chi_{\lambda/\mu,a}$ is a positive sum 
(see Theorem \ref{thm:positive-sum}). 
By the definition, we have
\begin{lem}\label{lem:intersecting}
Let $\lambda/\mu$ be a skew diagram satisfying  
the positivity condition \eqref{eq:positivity}, and 
let $(\alpha; \beta) \in \cH(\lambda/\mu)$. 
Then,  
\begin{equation}
\beta_{i+1}(n) \le \alpha^*_i(n), \qquad 
\beta^*_{i+1}(-n) \le \alpha_i(-n).
\end{equation} 
\end{lem}

\begin{figure}
\begin{psfrags}
\psfrag{S+}{$S_+$}
\psfrag{S-}{$S_-$}
\psfrag{-n}{\small $\scriptstyle -n$}
\psfrag{-n+1}{\small $\scriptstyle -n+1$}
\psfrag{-1}{\small $\scriptstyle -1$}
\psfrag{0}{\small $\scriptstyle 0$}
\psfrag{1}{\small $\scriptstyle 1$}
\psfrag{n-1}{\small $\scriptstyle n-1$}
\psfrag{n}{\small $\scriptstyle n$}
\psfrag{d}{\small $\cdots$}
\psfrag{height}{\small height}
\includegraphics[height=6cm, clip]{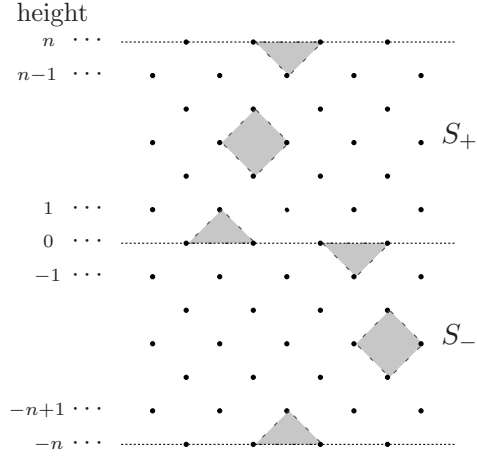}
\end{psfrags}
\caption{Examples of units. }\label{fig:units}
\end{figure}

A {\it unit} $U\subset S_{\pm}$ is either  
a unit square with its vertices on the lattice, 
or half of a unit square with its vertices on the lattice 
and the diagonal line on the boundary of $S_{\pm}$. 
See Figure \ref{fig:units} for examples. 
The {\it height} $\Ht(U)$ of $U$ is given by the height of the 
left vertex of $U$. 

\begin{defn}\label{def:unit}
Let $(\alpha; \beta) \in \cH(\lambda/\mu)$. 
For any unit $U \subset S_{\pm}$, let $\pm r=\Ht(U)$ and 
let $a$ and $a'=a+1$ be the horizontal positions of
the left and the right vertices of $U$.
Then,  
\begin{enumerate}
\item $U$ is called a {\it $\I$-unit} of $(\alpha; \beta)$ if 
there exists some $i$ ($0 \le i \le l$) such that 
\begin{equation}\label{eq:I-unit}
\begin{aligned}
& \alpha^*_i(r) \le a < a' \le \beta_{i+1}(r), \qquad 
\text{if $U \subset S_+$}, \\
& \alpha_i(-r) \le a < a' \le \beta^*_{i+1}(-r), \qquad 
\text{if $U \subset S_-$}. 
\end{aligned}
\end{equation}
\item $U$ is called a {\it $\II$-unit} of $(\alpha; \beta)$ if 
there exists some $i$ ($0 \le i \le l$) such that 
\begin{equation}\label{eq:II-unit}
\begin{aligned}
& \beta_{i+1}(r) \le a < a' \le \alpha_i^*(r), \qquad 
\text{if $U \subset S_+$}, \\
& \beta^*_{i+1}(-r) \le a < a' \le \alpha_i(-r), \qquad 
\text{if $U \subset S_-$}. 
\end{aligned}
\end{equation}
\end{enumerate}
Here, we set $\beta_{l+1}(r)=\beta^*_{l+1}(-r)=-\infty$ 
and $\alpha_0(-r)=\alpha^*_0(r)=+\infty$. 
Furthermore, a $\II$-unit $U$ of $(\alpha; \beta)$ 
is called a {\it boundary} $\II$-unit if \eqref{eq:II-unit}
holds for $i=0, l$, or $r=n$.  
\end{defn}

For a $\I$-unit, actually 
\eqref{eq:I-unit} does not hold for $i=0,l$. 
Also, it does not hold for $r=n$ if $\lambda/\mu$
satisfies the positivity condition \eqref{eq:positivity}, 
by Lemma \ref{lem:intersecting}. 

The {\it dual} $U^*$ of a unit $U$ is 
its image by the dual map \eqref{eq:dual-map}. 
Let $U$ and $U'$ be units. 
If the left or the right vertex of $U$ 
is also a vertex of $U'$, then 
we say that $U$ and $U'$ are {\it adjacent} 
and write $U\diamond U'$.   
It immediately follows from the definition that 
\begin{lem}\label{lem:I-II-unit}
\begin{enumerate}
\item 
A unit $U$ is a $\I$-unit (resp.~a $\II$-unit) if and only if 
the dual $U^*$ is a $\I$-unit (resp.~a $\II$-unit). 
\item No unit is simultaneously a $\I$- and $\II$-unit. 
\item \label{item:not-sim-unit}
If $U$ is a $\I$-unit and $U'$ is a $\II$-unit, 
then $U$ and $U'$ are not adjacent. 
\end{enumerate}
\end{lem}

Fix $(\alpha; \beta)\in \cH(\lambda/\mu)$. 
Let $\cU_{\I}$ be the set of all $\I$-units of $(\alpha; \beta)$, 
and let $\tilde{\cU}_{\I}:=\bigcup_{U\in \cU_{\I}}U$, 
where the union is taken for $U$ as a subset of $S_+ \sqcup S_-$. 
Let $\sim$ be the equivalence relation in $\cU_{\I}$ generated by 
the relation $\diamond$, and $[U]$ be its equivalence class of $U \in \cU_{\I}$. 
We call $\bigcup_{U'\in [U]}U'$  
a {\it connected component} of 
$\tilde{\cU}_{\I}$. 
For $\II$-units, 
$\cU_{\II}$, $\tilde{\cU}_{\II}$ 
and its connected component 
are defined similarly.

\begin{figure}
\begin{psfrags}
\psfrag{S+}{$S_+$}
\psfrag{S-}{$S_-$}
\psfrag{a1}{$\alpha_1$}
\psfrag{a2}{$\alpha_2$}
\psfrag{a3}{$\alpha_3$}
\psfrag{a4}{$\alpha_4$}
\psfrag{a5}{$\alpha_5$}
\psfrag{a6}{$\alpha_6$}
\psfrag{b1}{$\beta_1$}
\psfrag{b2}{$\beta_2$}
\psfrag{b3}{$\beta_3$}
\psfrag{b4}{$\beta_4$}
\psfrag{b5}{$\beta_5$}
\psfrag{b6}{$\beta_6$}
\psfrag{V}{$V$}
\psfrag{iota}{$\iota_1$}
\psfrag{-n}{\small $\scriptstyle -n \cdots$}
\psfrag{-n+1}{\small $\scriptstyle -n+1 \cdots$}
\psfrag{-1}{\small $\scriptstyle -1 \cdots$}
\psfrag{0}{\small $\scriptstyle 0 \cdots$}
\psfrag{1}{\small $\scriptstyle 1 \cdots$}
\psfrag{n-1}{\small $\scriptstyle n-1 \cdots$}
\psfrag{n}{\small $\scriptstyle n \cdots$}
\psfrag{v}{\small $\vdots$}
\psfrag{height}{\small height}
\includegraphics[height=6cm, clip]{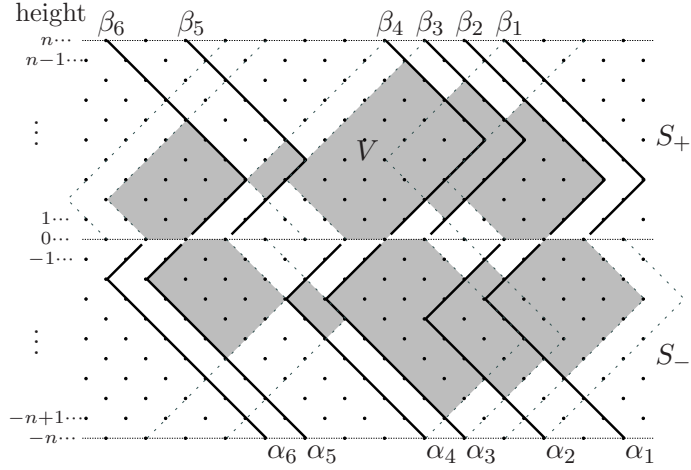}
\end{psfrags}
\caption{The undotted lines represent $\alpha_i$'s and $\beta_i$'s
while the dotted lines represent their duals, 
$\alpha^*_i$'s and $\beta^*_i$'s. The shaded area 
represents a $\I$-region $V$. 
}\label{fig:I-region}
\end{figure}

\begin{defn}
Let $\lambda/\mu$ be a skew diagram satisfying the positivity condition 
\eqref{eq:positivity}, 
and let $(\alpha; \beta) \in \cH(\lambda/\mu)$. 
\begin{enumerate}
\item
A connected component $V$ of $\tilde{\cU}_{\I}$
is called a {\it $\I$-region} of $(\alpha; \beta)$
if it contains at least one $\I$-unit of height $0$. 
\item 
A connected component $V$ of $\tilde{\cU}_{\II}$
is called a {\it $\II$-region} of $(\alpha; \beta)$ 
if it satisfies the following conditions: 
\begin{enumerate}[(i)]
\item $V$ contains at least one $\II$-unit of height $0$. 
\item $V$ does not contain any boundary $\II$-unit. 
\end{enumerate}
\end{enumerate}
\end{defn}
See Figure \ref{fig:I-region} for an example. 

\begin{prop}\label{prop:duality}
If $V$ is a $\I$- or $\II$-region, then 
$V^*=V$, where 
for a union of units $V=\bigcup U_i$, we define $V^*=\bigcup U^*_i$. 
\end{prop}

\begin{proof}
We remark that if two units are adjacent, then 
their duals are also adjacent. 
It follows that, for any $\I$-unit $U\subset V$, 
$U \sim U_0 \diamond U_0^* \sim U^*$ holds, where 
$U_0$ is any $\I$-unit $U \subset V$ of height $0$. 
Therefore, $U^* \subset V$. 
\end{proof}

\begin{figure}
\begin{psfrags}
\psfrag{S+}{$S_+$}
\psfrag{S-}{$S_-$}
\psfrag{a1}{$\alpha'_1$}
\psfrag{a2}{$\alpha'_2$}
\psfrag{a3}{$\alpha'_3$}
\psfrag{a4}{$\alpha'_4$}
\psfrag{a5}{$\alpha'_5$}
\psfrag{a6}{$\alpha'_6$}
\psfrag{b1}{$\beta'_1$}
\psfrag{b2}{$\beta'_2$}
\psfrag{b3}{$\beta'_3$}
\psfrag{b4}{$\beta'_4$}
\psfrag{b5}{$\beta'_5$}
\psfrag{b6}{$\beta'_6$}
\psfrag{iota}{$\iota_1$}
\psfrag{-n}{\small $\scriptstyle -n \cdots$}
\psfrag{-n+1}{\small $\scriptstyle -n+1 \cdots$}
\psfrag{-1}{\small $\scriptstyle -1 \cdots$}
\psfrag{0}{\small $\scriptstyle 0 \cdots$}
\psfrag{1}{\small $\scriptstyle 1 \cdots$}
\psfrag{n-1}{\small $\scriptstyle n-1 \cdots$}
\psfrag{n}{\small $\scriptstyle n \cdots$}
\psfrag{v}{\small $\vdots$}
\psfrag{V}{$V$}
\psfrag{height}{\small height}
\includegraphics[height=6cm, clip]{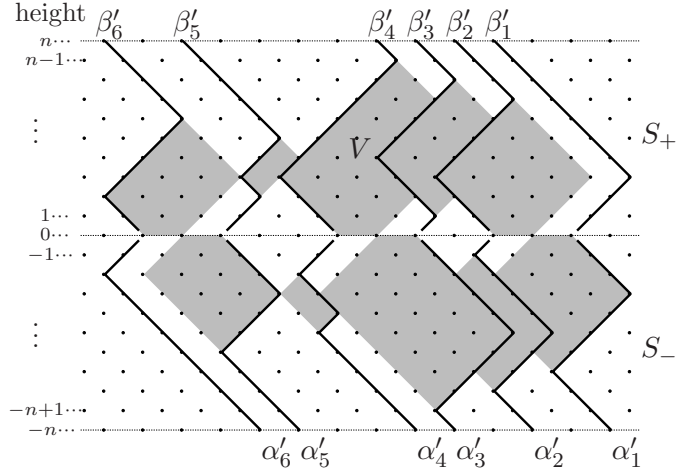}
\end{psfrags}
\caption{The tuple $(\alpha'; \beta'):=\varepsilon_V(\alpha; \beta)$ for
$(\alpha; \beta)$ with respect to $V$ in 
Figure \ref{fig:I-region}.}\label{fig:expansion}
\end{figure}

For any $(\alpha; \beta)\in \cH(\lambda/\mu)$, 
let $V$ be any $\I$- or $\II$-region of $(\alpha; \beta)$. 
Let $\alpha'_i$ be the lower path obtained from 
$\alpha_i$ by replacing the part $\alpha_i\cap V$ with $\beta^*_{i+1}\cap V$, 
and let $\beta'_i$ be the upper path obtained from 
$\beta_i$ by replacing the part $\beta_i\cap V$ with $\alpha^*_{i-1}\cap V$. 
Set 
$\varepsilon_V(\alpha; \beta):=
(\alpha'_1, \dots, \alpha'_l; \beta'_1, \dots, \beta'_l)$. 
See Figure \ref{fig:expansion} for an example.  

\begin{prop}\label{prop:I-II-region}
Let $\lambda/\mu$ be a skew diagram 
satisfying the positivity condition \eqref{eq:positivity}. 
Then, for any $(\alpha; \beta) \in \cH(\lambda/\mu)$, 
we have 
\begin{enumerate}
\item \label{item:prop-one}
For any $\I$- or $\II$-region 
$V$ of $(\alpha; \beta)$,
$\varepsilon_V(\alpha; \beta) \in \cH(\lambda/\mu)$. 
\item \label{item:prop-two}
For any $\I$-region $V$ of $(\alpha; \beta)$, 
$V$ is a $\II$-region of $\varepsilon_V(\alpha; \beta)$. 
\item For any $\II$-region $V$ of $(\alpha; \beta)$, 
$V$ is a $\I$-region of $\varepsilon_V(\alpha; \beta)$. 
\end{enumerate}
\end{prop}

\begin{proof}
We give a proof when $V$ is a $\I$-region. 

\eqref{item:prop-one} 
Set $(\alpha'; \beta'):=\varepsilon_V(\alpha; \beta)$. 
First, since $V$ does not contain any unit of height $\pm n$, we have 
$\alpha'_i(-n)=\alpha_i(-n)=\frac{n}{2}+\mu'_i+1-i$ 
and $\beta'_i(n)=\beta_i(n)=\frac{n}{2}+\lambda'_i+1-i$. 
Secondly, let us prove that $(\alpha'; \beta')$ 
is nonintersecting. 
Suppose, for example, 
if $(\alpha'_i, \alpha'_{i+1})$ is intersecting at a point $w$, 
then it implies that $(\alpha_i, \beta^*_{i+2})$ is intersecting 
at $w$. Set $-r= \Ht(w)$. Since  
$\alpha_i(-r)=\beta^*_{i+2}(-r) < \beta^*_{i+1}(-r)$, 
the unit $U\not\subset V$ whose left vertex is $w$ is a $\I$-unit. 
On the other hand, the unit $U'$ whose right vertex is $w$ is in $V$. 
This contradicts to the fact that $V$ is a connected component 
of $\tilde{\cU}_{\I}$. 

\eqref{item:prop-two}
It is obvious that a unit in $V$ is a $\II$-unit of $(\alpha'; \beta')$, 
and $U \sim U'$ for any two units $U, U' \subset V$. 
Assume that there exist some  
$\II$-unit $U''\not\subset V$ of $(\alpha'; \beta')$ 
which is adjacent to some $U \subset V$.  
Since $U''$ is a $\II$-unit 
of $(\alpha; \beta)$ and $U$ is a $\I$-unit of $(\alpha; \beta)$, 
it contradicts to Lemma \ref{lem:I-II-unit} \eqref{item:not-sim-unit}. 
Therefore, $V$ is a connected component of the $\II$-units of $(\alpha'; \beta')$. 
\end{proof}

We call the correspondence $(\alpha; \beta) \mapsto \varepsilon_V(\alpha; \beta)$ 
the {\it expansion} (resp.~the {\it folding}) with respect to $V$, 
if $V$ is a $\I$-region (resp.~a $\II$-region) of $(\alpha; \beta)$. 
We remark that $\varepsilon_V \circ \varepsilon_V = \id$
for any $\I$- or $\II$-region $V$.

\begin{figure}
\begin{psfrags}
\psfrag{a}{$\alpha_i$}
\psfrag{a*}{$\alpha^*_i$}
\psfrag{a'}{$\alpha'_i$}
\psfrag{b}{$\beta_{i+1}$}
\psfrag{b*}{$\beta^*_{i+1}$}
\psfrag{b'}{$\beta'_{i+1}$}
\psfrag{V}{$V$}
\includegraphics[width=12cm, clip]{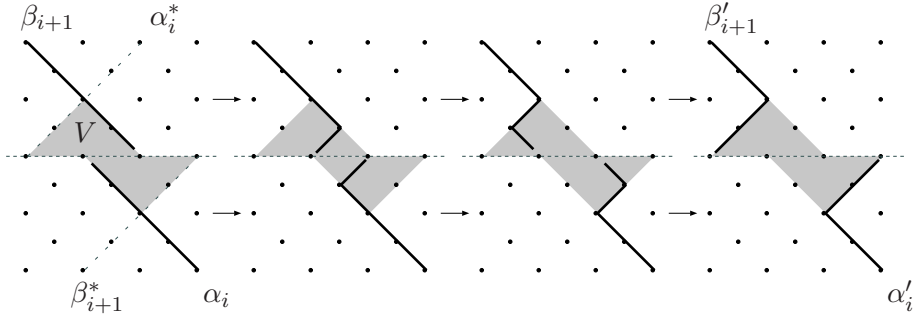}
\end{psfrags}
\caption{An example of a procedure 
of the expansion $(\alpha; \beta) \mapsto \varepsilon_V(\alpha; \beta)$ 
with respect to a $\I$-region $V$ 
at a pair $(\alpha_i, \beta_{i+1})$, 
by each unit.}\label{fig:procedure}
\end{figure}

\begin{rem}\label{rem:transformations}
The expansion and the folding are decomposed into a series of deformations of paths 
along each unit in $V$. See Figure \ref{fig:procedure}.
This is a key fact in the proof of the weight-preserving property of the maps 
$\iota_2$ in Section \ref{sec:positive-sum} and $\phi$ in 
Section \ref{sec:folding-map}. 
\end{rem}

\subsection{The second involution and 
an expression of $\chi_{\lambda/\mu,a}$ 
by a positive sum}\label{sec:positive-sum}

{}From now, we assume that 
$\lambda/\mu$ satisfies the positivity condition \eqref{eq:positivity}. 

Let $\bp \in P_1(\lambda/\mu)$, and 
let $p_i(\pm n)$ be the horizontal position of $p_i$ at height $\pm n$. 
Then $p_i(-n)<p_j(-n)$ for any $i<j$. We call  
a pair $(p_i, p_j)$, $i<j$ {\it transposed} if $p_i(n)>p_j(n)$. 

For each $\bp \in P_1(\lambda/\mu)$, one can uniquely 
associate $(\alpha; \beta) \in \cH(\lambda/\mu)$ 
by removing all the E-steps from $\bp$. 
We write $\pi(\bp)$ for $(\alpha;\beta)$. 
A $\I$- or $\II$-region of $(\alpha; \beta)=\pi(\bp)$ 
is also called a {\it $\I$- or $\II$-region of $\bp$}.

Let $\bp \in P_1(\lambda/\mu)$ and 
$(\alpha; \beta)=\pi(\bp)$. 
If $h:=\alpha_i(0)-\beta_{i+1}(0)$ is 
a non-positive number (resp.~a positive number), 
then we call a pair $(\alpha_i, \beta_{i+1})$ an {\it overlap}  
(resp.~a {\it hole}). 
Furthermore, 
if $h$ is an even number (resp.~an odd number), 
then we say that 
$(\alpha_i, \beta_{i+1})$ is {\it even} (resp.~{\it odd}). 
Using that no triple $(p_i,p_j,p_k)$ 
exists for $\bp \in P_1(\lambda/\mu)$ 
which is intersecting at a point, we have 

\begin{lem}\label{lem:alpha-beta}
Let $(\alpha; \beta)=\pi(\bp)$ for $\bp \in P_1(\lambda/\mu)$. 
Then, for any $i$, 
\begin{enumerate}\label{item:alpha-beta-one}
\item $(\alpha_i, \beta_{i+1})$ is an odd overlap 
if and only if $(p_i, p_j)$ is a specially intersecting, 
non-transposed pair for some $j>i$. 
\item \label{item:alpha-beta-two}
$(\alpha_i, \beta_{i+1})$ is an even overlap
if and only if $(p_i, p_j)$ is a transposed pair for some $j>i$.  
\item \label{item:alpha-beta-three}
$(\alpha_i, \beta_{i+1})$ is a hole  
if and only if $(p_i,p_j)$ is not intersecting for any $j>i$. 
\end{enumerate}
\end{lem}

\begin{figure}
\begin{psfrags}
\psfrag{bn}{$\scriptscriptstyle \overline{n}$}
\psfrag{n}{$\scriptscriptstyle n$}
\psfrag{ai}{$\scriptstyle \alpha_i$}
\psfrag{a*i}{$\scriptstyle \alpha^*_i$}
\psfrag{a'i}{$\scriptstyle \alpha'_i$}
\psfrag{a'*i}{$\scriptstyle {\alpha'_i}^*$}
\psfrag{aj}{$\scriptstyle \alpha_j$}
\psfrag{a*j}{$\scriptstyle \alpha^*_j$}
\psfrag{a'j}{$\scriptstyle \alpha'_j$}
\psfrag{a'*j}{$\scriptstyle {\alpha'_j}^*$}
\psfrag{bi+1}{$\scriptstyle \beta_{i+1}$}
\psfrag{b*i+1}{$\scriptstyle \beta^*_{i+1}$}
\psfrag{b'i+1}{$\scriptstyle \beta'_{i+1}$}
\psfrag{b'*i+1}{$\scriptstyle {\beta'}^*_{i+1}$}
\psfrag{pi}{$p_i$}
\psfrag{p'i}{$p'_i$}
's\psfrag{pj}{$p_j$}
\psfrag{p'j'}{$p'_{k}$}
\psfrag{v}{$\vdots$}
\psfrag{V}{$V$}
\psfrag{ev}{$\scriptstyle \varepsilon_V(\bp)$}
\psfrag{ev'}{$\scriptstyle \varepsilon_V(\bp')$}
\begin{tabular}{ll}
  \begin{tabular}{ll}
    (a) 
    & $(\alpha_i, \beta_{i+1})$: an even overlap 
  \end{tabular}
  & 
  \begin{tabular}{ll}
    (b) 
    & $(\alpha_i, \beta_{i+1})$: an odd overlap 
  \end{tabular}
  \\
  \includegraphics[height=7.5cm, clip]{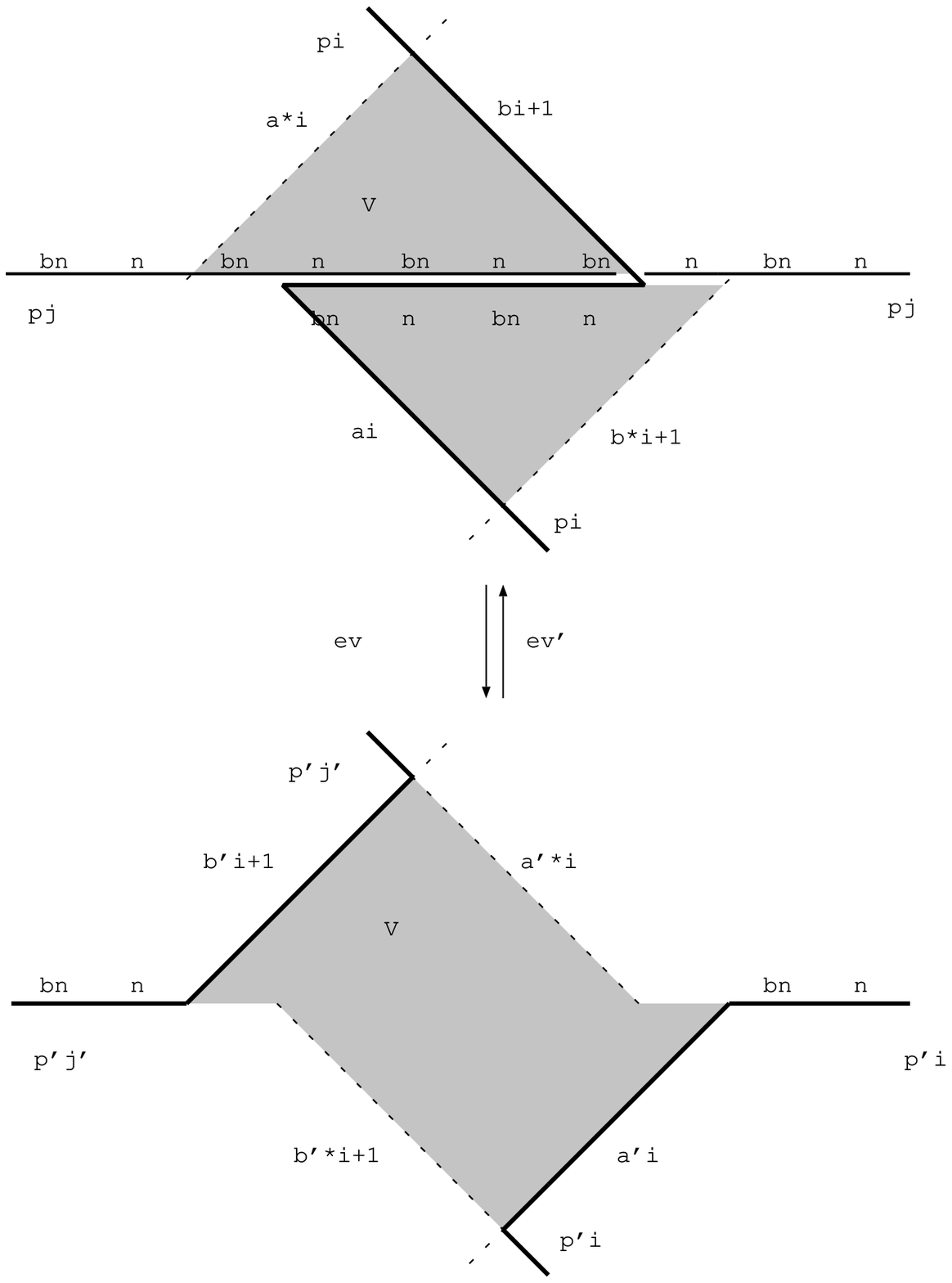}\quad 
  &
  \quad 
  \includegraphics[height=7.5cm, clip]{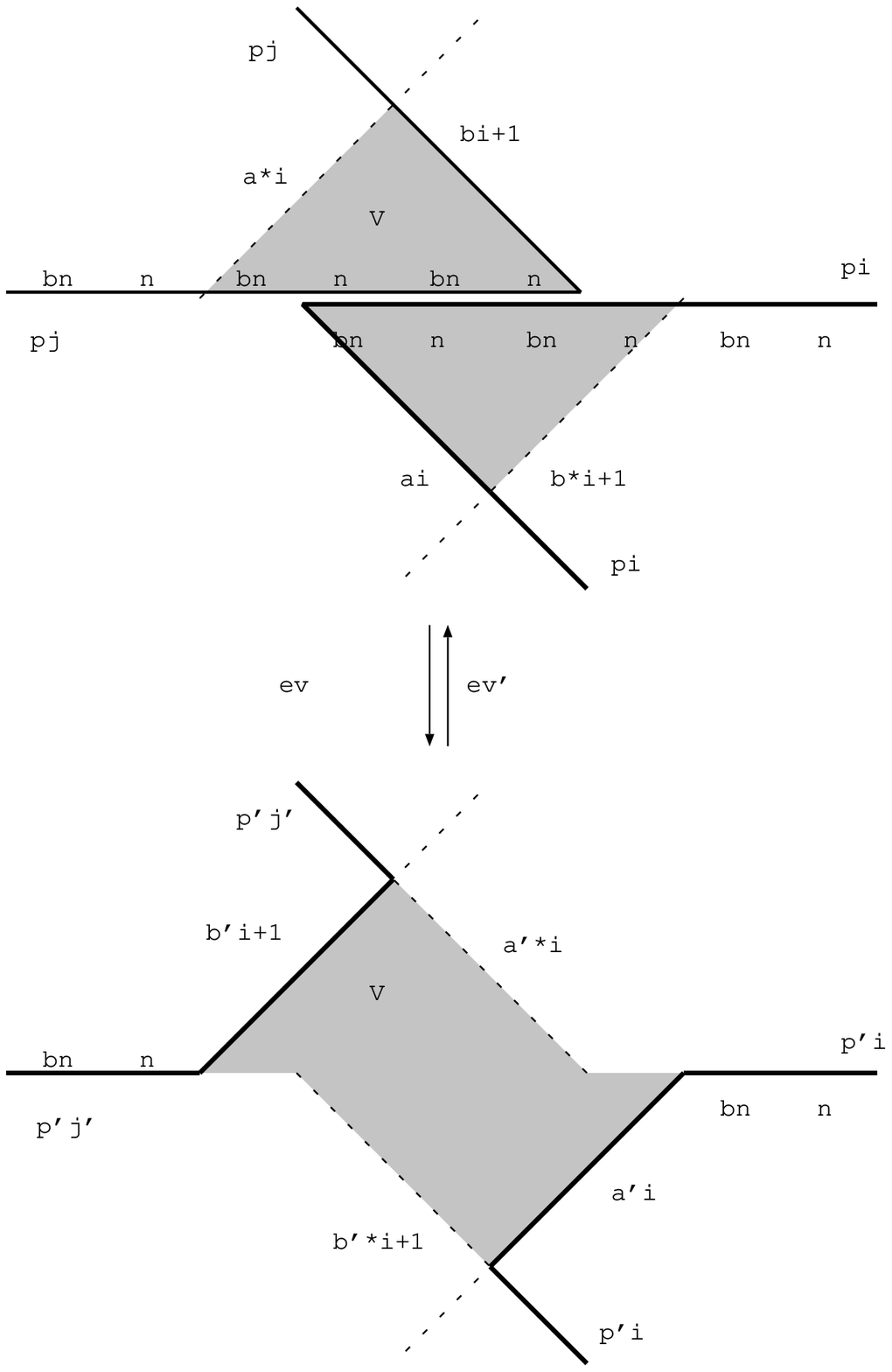}\\
  \begin{tabular}{ll}
    \phantom{(a)} 
    & $(\alpha'_i, \beta'_{i+1})$: an even hole 
  \end{tabular}
  & 
  \begin{tabular}{ll}
    \phantom{(b)} 
    & $(\alpha'_i, \beta'_{i+1})$: an odd hole 
  \end{tabular}
\end{tabular}
\end{psfrags}
\caption{The deformation $\varepsilon_V:\bp \leftrightarrow \bp'$ 
with respect to a $\I$-region $V$ of $\bp$
and a $\II$-region $V$ of $\bp'$.  
}\label{fig:deformation}
\end{figure}

Let $\bp\in P_1(\lambda/\mu)$, 
$(\alpha; \beta) =\pi(\bp)$, and $V$ be a $\I$- or $\II$-region of $\bp$. 
Then, there exists $\bp' \in P_1(\lambda/\mu)$ such that 
$$\varepsilon_V(\alpha;\beta)=\pi(\bp').$$
It is constructed from $\bp$ as follows, 
which is well-defined by Lemma \ref{lem:alpha-beta}:  
\begin{enumerate}[A.]
\item \label{item:expansion}{\it The case of a $\I$-region $V$}.  
For any $i$, replace $(\alpha_i, \beta_{i+1})$ 
in $\bp$ with $(\alpha'_i, \beta'_{i+1})$. 
Furthermore, for any $i$ such that 
$(\alpha_i, \beta_{i+1})$ is 
an overlap and intersects with $V$ 
at height $0$, remove the E-steps between $\beta'_{i+1}(0)$ 
and $\alpha'_i(0)$. See 
Figure \ref{fig:deformation}. 
\item {\it The case of a $\II$-region $V$}.  
This is the reverse operation of Case \ref{item:expansion}. 
Namely, for any $i$, replace $(\alpha_i, \beta_{i+1})$
in $\bp$ with $(\alpha'_i, \beta'_{i+1})$.
Furthermore, for any $i$ such that
$(\alpha_i, \beta_{i+1})$ is
a hole and intersects with $V$
at height $0$, then add the E-steps between $\beta'_{i+1}(0)$
and $\alpha'_i(0)$ 
as in Figure \ref{fig:deformation} (a) (for an even hole) 
and Figure \ref{fig:deformation} (b) (for an odd hole) 
wherein $\{\alpha_i, \beta_{i+1}, p_i, p_j\}$ 
and $\{\alpha'_i, \beta'_{i+1}, p'_i, p'_{k}\}$ are interchanged. 
\end{enumerate}

We call the correspondence $\bp \mapsto \bp'$ the {\it expansion}
(resp. the {\it folding}) of $\bp$ with respect to 
a $\I$-region (resp. a $\II$-region) $V$,
and write $\varepsilon_V(\bp):=\bp'$.

For any $\I$-region $V$ (resp.~$\II$-region $V$) 
of $\bp \in P_1(\lambda/\mu)$ with $(\alpha;\beta)=\pi(\bp)$, we set 
\begin{equation}\label{eq:m}
n(V):=\# \Big\{ i \Big| \text{
\begin{tabular}{l}
$(\alpha_i, \beta_{i+1})$ 
is an even overlap (resp.~an even hole)\\
which intersects with $V$ at height $0$
\end{tabular}
}\Big\}.
\end{equation}
Let $V$ be a $\I$- or $\II$-region of $\bp \in P_1(\lambda/\mu)$. 
By Lemma \ref{lem:alpha-beta}, 
$n(V)$ is equal to the number of the transposed pairs 
$(p_i, p_j)$ in $\bp$ which intersect with 
$V$ at height $0$. Moreover, since 
the expansion (resp.~the folding) $\bp \mapsto \varepsilon_V(\bp)$ 
is a deformation that `resolves' all the transposed pairs 
(resp.~transposes all the even holes) in $\bp$ 
which intersect with $V$ at height $0$, 
we have 
\begin{lem}\label{lem:sign}
Let $\bp \in P_1(\lambda/\mu)$ and $V$ be a 
$\I$- or $\II$-region of $\bp$. Then,  
$$(-1)^{\varepsilon_V(\bp)}=(-1)^{n(V)}\cdot (-1)^{\bp}.$$
\end{lem}

\begin{defn}
We say that a $\I$- or $\II$-region $V$ is 
{\it even} (resp.~{\it odd}) if $n(V)$ is even (resp.~odd). 
\end{defn}

Let 
$P_{\odd}(\lambda/\mu)$ be the set of all 
$\bp \in P_1(\lambda/\mu)$ which have at least one 
odd $\I$- or $\II$-region of $\bp$.
We can define an involution 
$$\iota_2:P_{\odd}(\lambda/\mu) \to P_{\odd}(\lambda/\mu)$$
as follows: 
Let $V$ be the unique odd $\I$- or $\II$-region of 
$\bp \in P_{\odd}(\lambda/\mu)$ such that the value 
$\max \{ \Position{w} \mid w \in V, \  \Ht(w)=0 \}$
is greatest among all the odd $\I$- or $\II$-regions of $\bp$, 
and set $\iota_2(\bp)=\varepsilon_V(\bp)$. 
Then we have 

\begin{prop}\label{prop:second-involution}
The map $\iota_2: P_{\odd}(\lambda/\mu) \to P_{\odd}(\lambda/\mu)$ 
is a weight-preserving, sign-reversing involution. 
\end{prop}

\begin{proof}
The map $\iota_2$ is an involution because 
$\varepsilon_V \circ \varepsilon_V = \id$, and 
sign-reversing by Lemma \ref{lem:sign}. 
We prove that $\iota_2$ is weight-preserving 
in the case where $\bp \mapsto \bp':=\iota_2(\bp)$ is an expansion
with respect to a $\I$-region $V$ of $\bp$. 
Let $(\alpha; \beta)=\pi(\bp)$, 
and we decompose the weights $z_a^{\bp}$ and $z_a^{\bp'}$ 
in \eqref{eq:weight} into two parts as $z_a^{\bp}=HE$ 
and $z_a^{\bp}=H'E'$
where $H$ and $H'$ are the factors from the 
$e$-labeling on $(\alpha; \beta)$ and $(\alpha'; \beta')$, 
while $E$ and $E'$ are the ones from the $e$-labeling 
on the height $0$ part (the E-steps) of $\bp$ and $\bp'$. 
By Remark \ref{rem:transformations}, 
we have $H'=H\delta$, where 
$$\delta:=\prod_{U\subset V:\,\text{unit}}\delta(U), $$
and, for any unit $U\subset V$ in $S_{\pm}$ 
of height $\pm r$ with left vertex $(x,y)$, 
$$\delta(U):=
\begin{cases}
z_{\overline{n-r}, a-2x}/z_{\overline{n-r+1}, a-2x}, & 
\text{if $r\ne 0$ and $U \subset S_+$}, \\
z_{n-r, a-2x}/z_{n-r+1, a-2x}, &
\text{if $r\ne 0$ and $U \subset S_-$}, \\
z_{\overline{n}, a-2x}, & 
\text{if $r=0$ and $U \subset S_+$}, \\
z_{n, a-2x}, & 
\text{if $r=0$ and $U \subset S_-$}. 
\end{cases}
$$ 
Using the relations in \eqref{eq:relations}, 
we have $\delta(U)\cdot \delta(U^*)=1$ for any $U$ 
whose height is not $0$. 
Therefore, combining $\delta(U)$ for all the $\I$-units in $V$, 
we obtain 
$$\delta=\prod_{
\genfrac{}{}{0pt}{1}{U\subset V:\,\text{unit}}{\Ht(U)=0}
}\delta(U)=
\prod_{i=1}^{l-1}
\left(
\prod_{k=\alpha^*_i(0)}^{\beta_{i+1}(0)-1}z_{\overline{n},a-2k}
\prod_{k=\alpha_i(0)}^{\beta^*_{i+1}(0)-1}z_{n,a-2k}
\right). 
$$
See Figure \ref{fig:deformation}.  
On the other hand, we have 
$E'=E\delta^{-1}$, and therefore, we obtain 
$z_a^{\iota_2(\bp)}=H'E'=HE=z_a^{\bp}$. 
\end{proof}

It follows from Proposition \ref{prop:second-involution} 
that the contributions of $P_{\odd}(\lambda/\mu)$
to the sum \eqref{eq:first-sum} cancel each other.

Let 
$P_2(\lambda/\mu):= 
P_1(\lambda/\mu) \backslash P_{\odd}(\lambda/\mu)$, i.e., 
the set of all  
$\bp \in P_1(\lambda/\mu)$ which satisfy the following conditions:  
\begin{enumerate}[(i)]
\item \label{item:Ptwo-no-ord}
$\bp$ does not have any ordinarily intersecting pair $(p_i, p_j)$.
\item \label{item:Ptwo-no-odd}
$\bp$ does not have any odd $\I$- or $\II$-region. 
\end{enumerate}
Every $\bp \in P_2(\lambda/\mu)$ 
has an even number of transposed pairs, which implies that 
$(-1)^{\bp}=1$. 
Thus, the sum \eqref{eq:first-sum} reduces to a positive sum, 
and we have 
\begin{thm}\label{thm:positive-sum}
For any skew diagram $\lambda/\mu$
satisfying the positivity condition \eqref{eq:positivity}, 
we have
\begin{equation}\label{eq:positive-sum}
\chi_{\lambda/\mu, a}= \sum_{\bp \in P_2(\lambda/\mu)}z_a^{\bp}. 
\end{equation}
\end{thm}

\section{The folding map and a tableaux description}
\label{sec:tableau-description}
In this section, we give a tableaux description of 
$\chi_{\lambda/\mu,a}$. Namely, the sum \eqref{eq:positive-sum}
is translated into the one over a set of the tableaux of shape $\lambda/\mu$ 
which satisfy certain conditions called the {\it horizontal}, 
{\it vertical}, and {\it extra} rules. 

\subsection{The folding map}\label{sec:tableau-description-1}
Since a path $\bp \in P_2(\lambda/\mu)$ 
in \eqref{eq:positive-sum} might have (an even number of) 
transposed pairs $(p_i, p_j)$, 
the sum \eqref{eq:positive-sum} cannot be translated 
into a tableaux description yet. 
Therefore, we introduce another set of paths as follows.

Let $P(\lambda/\mu)$ be the set of all 
$\bp =(p_1, \dots, p_l)\in \fP(\id; \bu_{\mu}, \bv_{\lambda})$ 
such that 
\begin{enumerate}[(i)]
\item \label{item:p-no-ord} $p$ does not have any ordinarily intersecting 
{\it adjacent} pair $(p_i, p_{i+1})$. 
\item \label{item:p-no-odd} $p$ does not have any odd $\II$-region. 
\end{enumerate}
Here, an odd $\II$-region of $\bp \in P(\lambda/\mu)$ 
is defined in the same way as that of $\bp \in P_1(\lambda/\mu)$. 
The following fact is not so trivial. 
\begin{prop}\label{prop:folding-map}
There exists a weight-preserving bijection 
$$\phi: P_2(\lambda/\mu) \to P(\lambda/\mu). $$ 
\end{prop}
The map $\phi$ is called the {\it folding map}. 
Roughly speaking, it is an iterated application of 
(some generalization of) the folding in Section \ref{sec:second-involution}. 
The construction of $\phi$ 
is the most technical part of the paper.  
We provide the details in Section \ref{sec:folding-map}. 
Admitting Proposition \ref{prop:folding-map}, 
we immediately have 
\begin{thm}\label{thm:tilde-positive-sum}
For any skew diagram $\lambda/\mu$ satisfying
the positivity condition \eqref{eq:positivity}, 
we have
\begin{equation}\label{eq:path-description}
\chi_{\lambda/\mu, a}= \sum_{\bp \in P(\lambda/\mu)} z_a^{\bp}.
\end{equation}
\end{thm}

\subsection{Tableaux description}\label{sec:extra}
Define a partial order in $I$ in \eqref{eq:entries} by
$$1 \prec 2 \prec \dots \prec n-1 \prec 
\genfrac{}{}{0pt}{0}{n}{\overline{n}}
\prec \overline{n-1}\prec \dots \prec \overline{2} \prec\overline{1}. $$

A {\it tableau} $T$ of shape $\lambda/\mu$ is the skew diagram $\lambda/\mu$ 
with each box filled by one entry of $I$. 
For a tableau $T$ and $a \in \bC$, we define 
the {\it weight} of $T$ as 
$$z_a^T=\prod_{(i,j) \in \lambda/\mu}z_{T(i,j),a+2(j-i)}, $$
where $T(i,j)$ is the entry of $T$ at $(i,j)$. 

\begin{defn}
A tableau $T$ (of shape $\lambda/\mu$) is 
called an {\it $\HV$-tableau} if it satisfies the following conditions: 
\newline
\begin{tabular}{ll}
$(\bH)$ & horizontal rule \ $T(i,j) \preceq T(i,j+1)$ or 
$(T(i,j), T(i,j+1))=(n, \overline{n})$. \\
$(\bV)$ & vertical rule \quad \ $T(i,j) \not\succeq T(i+1,j)$. 
\end{tabular}
\end{defn}
We denote the set of all $\HV$-tableaux of shape $\lambda/\mu$ 
by $\HVTab(\lambda/\mu)$. 

\begin{rem}
The configuration $(T(i,j), T(i,j+1))=(n, \overline{n})$ 
is prohibited later by another rule. 
See Remark \ref{rem:E-one}. 
\end{rem}

Let $\HVP(\lambda/\mu)$ be the set of all 
$\bp\in \fP(\id; \bu_{\mu}, \bv_{\lambda})$
which do not have any ordinarily intersecting 
adjacent pair $(p_i, p_{i+1})$.   
With any $\bp \in \HVP(\lambda/\mu)$, we associate a tableau $T$
of shape $\lambda/\mu$ as follows: 
For any $j=1, \dots, l$, 
let $E(p_j)= \{ s_{i_1}, s_{i_2}, \dots, s_{i_m} \}$ 
$(i_1<i_2< \dots<i_m)$ 
be the set defined as in \eqref{eq:label-steps}, and set 
\begin{equation*}
T(\mu'_j+k, j)=L_a^1(s_{i_k}), \qquad k=1, \dots, m, 
\end{equation*}
where $L_a^1$ is the first component of the $e$-labeling \eqref{eq:e-label}. 
It is easy to see that $T$ satisfies the vertical rule $(\bV)$ 
because of the definition of the $e$-labeling of $p_j$,  
and satisfies the horizontal rule $(\bH)$ because $\bp$ 
does not have any ordinarily intersecting adjacent pair. 
Therefore, if we set $\cTv:\bp \mapsto T$, we have
\begin{prop}\label{prop:path-tableau}
The map 
$$\cTv: \HVP(\lambda/\mu) \to \HVTab(\lambda/\mu)$$
is a weight-preserving bijection. 
\end{prop}

Let $\Tab(\lambda/\mu):=\cTv(P(\lambda/\mu))$. 
In other words, $\Tab(\lambda/\mu)$ is the set of 
all the tableaux $T$ which satisfy 
$(\bH)$, $(\bV)$, and the following {\it extra rule}: 
\vspace*{5pt}

\noindent
\begin{tabular}{ll} 
$(\bE)$ & The corresponding $\bp = \cTv^{-1}(T)$ does not have 
any odd $\II$-region. 
\vspace{5pt}
\end{tabular}

\noindent
By Theorem \ref{thm:tilde-positive-sum} and Proposition \ref{prop:path-tableau}, 
we 
obtain a tableaux description of $\chi_{\lambda/\mu,a}$, 
which is the main result of the paper. 
\begin{thm}\label{thm:tableau-description}
For any skew diagram $\lambda/\mu$ satisfying
the positivity condition \eqref{eq:positivity}, 
we have 
$$\chi_{\lambda/\mu, a} = \sum_{T \in \Tab(\lambda/\mu)} z_a^{T}. $$
\end{thm}

\subsection{Extra rule in terms of tableau}
It is straightforward to translate the extra rule 
$(\bE)$ into tableau language. We only give 
the result here. 

Fix an $\HV$-tableau $T$. 
For any $a_1, \dots, a_m \in I$, 
let $C(a_1, \dots, a_m)$ be a configuration 
in $T$ as follows: 
\begin{equation}
\text{
{\setlength{\unitlength}{0.2mm}
\begin{picture}(25,100)
\multiput(0,0)(25,0){2}{\line(0,1){110}}
\multiput(0,0)(0,25){2}{\line(1,0){25}}
\multiput(0,110)(0,-25){3}{\line(1,0){25}}
\put(4,92){$a_1$}
\put(4,67){$a_2$}
\put(9,34){$\vdots$}
\put(2,7){$a_{\scriptscriptstyle m}$}
\end{picture}
}
}
\end{equation}
If $1 \preceq a_1 \prec \dots \prec a_m \preceq n$, 
then we call it an {\it L-configuration}. 
If $\overline{n} \preceq a_1 \prec \dots \prec a_m \preceq \overline{1}$, 
then we call it a {\it U-configuration}. 
Note that an L-configuration 
corresponds to a part of a lower path, while a U-configuration 
corresponds to a part of an upper path under the map $\cTv$. 

\begin{figure}
{\small 
\begin{psfrags}
\psfrag{r}{$r$}
\psfrag{n}{$n$}
\psfrag{-n}{$-n$}
\psfrag{k}{$n-k+1$}
\psfrag{-k}{$-n+k-1$}
\psfrag{0}{$0$}
\psfrag{v}{$\vdots$}
\psfrag{d}{$\dots$}
\psfrag{a1}{$a_1$}
\psfrag{as}{$a_s$}
\psfrag{a}{$a$}
\psfrag{a'1}{$\overline{a'_1}$}
\psfrag{a's}{$\overline{a'_{s'}}$}
\psfrag{b1}{$\overline{b_1}$}
\psfrag{bt}{$\overline{b_t}$}
\psfrag{b}{$b$}
\psfrag{b'1}{$b'_1$}
\psfrag{b't}{$b'_{t'}$}
\psfrag{pj}{$p_j$}
\psfrag{pj+1}{$p_{j+1}$}
\psfrag{height}{height}
\includegraphics[height=8cm, clip]{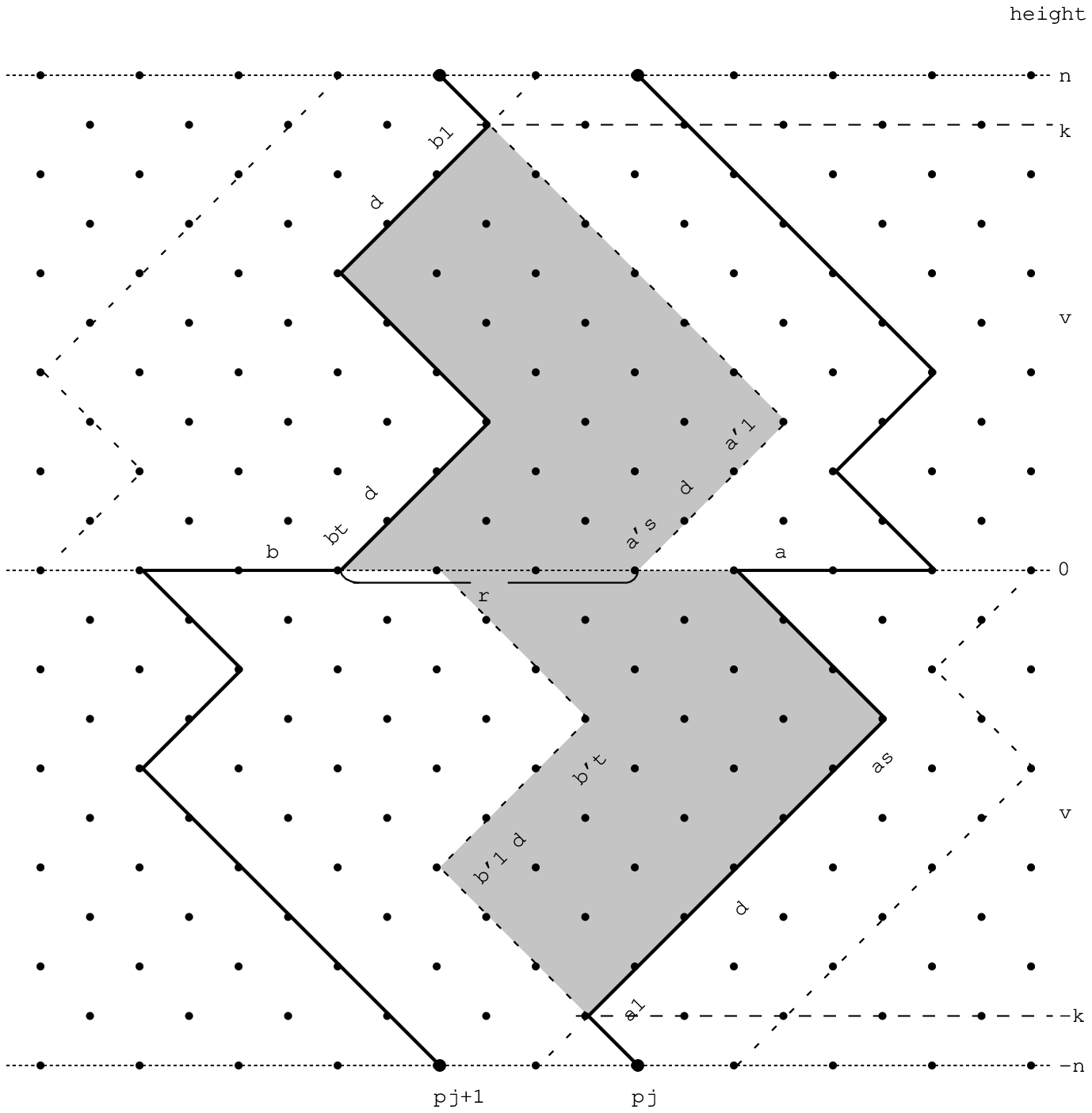}
\end{psfrags}
}
\caption{An example of adjacent paths $(p_j,p_{j+1})$ such that 
a part of it corresponds to an LU-configuration 
of type 1 as in \eqref{eq:type-one}. 
}\label{fig:tableau-path-correspondence}
\end{figure}

Let $(L,U)$ be 
a pair of an L-configuration $L=C(a_1, \dots, a_s)$ in 
the $j$th column and a U-configuration 
$U=C(\overline{b_t}, \dots, \overline{b_1})$ in the $(j+1)$th column. 
We call it an {\it LU-configuration} of $T$ if it satisfies 
one of the following two conditions: 

Condition 1. {\it LU-configuration of type 1}. 
$(L,U)$ has the form  
\begin{equation}\label{eq:type-one}
\text{
{\setlength{\unitlength}{0.2mm}
\begin{picture}(50,170)
{\thicklines
\put(0,55){\line(0,1){115}}
\put(25,0){\line(0,1){170}}
\put(50,0){\line(0,1){115}}
\put(0,55){\line(1,0){25}}
\put(0,80){\line(1,0){25}}
\put(0,145){\line(1,0){25}}
\put(0,170){\line(1,0){25}}
\put(25,0){\line(1,0){25}}
\put(25,25){\line(1,0){25}}
\put(25,90){\line(1,0){25}}
\put(25,115){\line(1,0){25}}
}
\put(4,150){$a_1$}
\put(10,100){$\vdots$}
\put(4,62){$a_s$}
\put(5,38){$a$}
\put(29,5){$\overline{b_1}$}
\put(35,50){$\vdots$}
\put(29,95){$\overline{b_t}$}
\put(32,121){$b$}
\multiput(50,0)(5,0){10}{\line(1,0){2}}
\multiput(25,170)(5,0){15}{\line(1,0){2}}
\put(80,0){\vector(0,1){170}}
\put(80,170){\vector(0,-1){170}}
\put(85,75){$\scriptstyle n-k+1$}
\multiput(0,55)(-5,0){10}{\line(-1,0){2}}
\multiput(0,115)(-5,0){10}{\line(-1,0){2}}
\put(-30,55){\vector(0,1){60}}
\put(-30,115){\vector(0,-1){60}}
\put(-45,85){$\scriptstyle r$}
\end{picture}
}
}
\end{equation} 
for some $k$ and $r$ with $1 \le k \le n$, $1 \le r \le \min \{ s, t\}$, 
$n-k+1=s+t-r$, and 
\begin{gather}
  a_1 =  k, \qquad 
  \overline{b_1} = \overline{k}, \\
  a \succeq \overline{n} \text{ if $a$ exists}, \qquad 
  b \preceq n \text{ if $b$ exists}, 
\\
a_{i+1} \preceq b'_i, \quad (1 \le i \le t'), \qquad 
\overline{b_{i+1}} \succeq \overline{a'_i}, \quad (1 \le i \le s'), 
\label{eq:intersecting-cond}
\end{gather}
where $a'_1 \prec \dots \prec a'_{s'}$ ($s':=t-r$)
and $b'_1 \prec \dots \prec b'_{t'}$ ($t':=s-r$)
are defined as 
\begin{equation}\label{eq:dual-type-one}
\begin{gathered}
\{a_1, \dots, a_s \} \sqcup \{ a'_1, \dots, a'_{s'} \}
 = \{ k, k+1, \dots, n \}, \\
\{ \overline{b_1}, \dots, \overline{b_{t}} \} \sqcup 
\{ \overline{b'_1}, \dots, \overline{b'_{t'}} \} 
 = \{ \overline{k}, \overline{k+1}, \dots, \overline{n} \}. 
\end{gathered}
\end{equation}
See Figure \ref{fig:tableau-path-correspondence} 
for the corresponding part in the paths. 
In particular, if 
$r$ is odd, then we say that $(L,U)$ is {\it odd}.  

\begin{figure}
{\small 
\begin{psfrags}
\psfrag{n}{$n$}
\psfrag{-n}{$-n$}
\psfrag{k}{$n-k+1$}
\psfrag{-k}{$-n+k-1$}
\psfrag{k'}{$n-k'$}
\psfrag{-k'}{$-n+k'$}
\psfrag{0}{$0$}
\psfrag{v}{$\vdots$}
\psfrag{d}{$\cdots$}
\psfrag{a1}{$a_1$}
\psfrag{as}{$a_s$}
\psfrag{a}{$a$}
\psfrag{a'1}{$\overline{a'_1}$}
\psfrag{a's}{$\overline{a'_{s'}}$}
\psfrag{b1}{$\overline{b_1}$}
\psfrag{bt}{$\overline{b_t}$}
\psfrag{b}{$b$}
\psfrag{b'1}{$b'_1$}
\psfrag{b't}{$b'_{t'}$}
\psfrag{height}{height}
\psfrag{pj}{$p_j$}
\psfrag{pj+1}{$p_{j+1}$}
\includegraphics[height=8cm, clip]{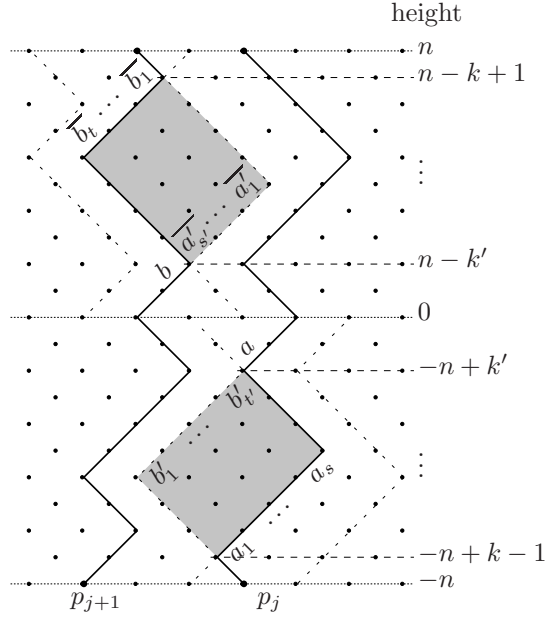}
\end{psfrags}
}
\caption{An example of adjacent paths $(p_j,p_{j+1})$ such that
a part of it corresponds to an LU-configuration
of type 2 as in \eqref{eq:type-two}.}\label{fig:tableau-path-correspondence-two}
\end{figure}

Condition 2. {\it LU-configuration of type 2}. 
$(L,U)$ has the form  
\begin{equation}\label{eq:type-two}
\text{
{\setlength{\unitlength}{0.2mm}
\begin{picture}(50,210)
{\thicklines
\put(0,130){\line(0,1){80}}
\put(25,130){\line(0,1){80}}
\put(25,0){\line(0,1){80}}
\put(50,0){\line(0,1){80}}
\put(0,130){\line(1,0){25}}
\put(0,155){\line(1,0){25}}
\put(0,185){\line(1,0){25}}
\put(0,210){\line(1,0){25}}
\put(25,0){\line(1,0){25}}
\put(25,25){\line(1,0){25}}
\put(25,55){\line(1,0){25}}
\put(25,80){\line(1,0){25}}
}
\put(5,191){$a_1$}
\put(10,160){$\vdots$}
\put(5,135){$a_s$}
\put(8,115){$a$}
%
%
%
\put(31,5){$\overline{b_1}$}
\put(35,30){$\vdots$}
\put(31,61){$\overline{b_t}$}
\put(34,85){$b$}
\multiput(50,0)(5,0){10}{\line(1,0){2}}
\multiput(25,210)(5,0){15}{\line(1,0){2}}
\put(80,0){\vector(0,1){210}}
\put(80,210){\vector(0,-1){210}}
\put(85,90){$\scriptstyle n-k+1$}
\multiput(25,80)(-5,0){15}{\line(-1,0){2}}
\multiput(0,130)(-5,0){10}{\line(-1,0){2}}
\put(-30,80){\vector(0,1){50}}
\put(-30,130){\vector(0,-1){50}}
\put(-70,100){$\scriptstyle n-k'$}
\multiput(25,80)(0,5){10}{\line(0,1){2}}
\end{picture}
}
}
\end{equation}
for some $k$ and $k'$ with $1 \le k  < k' \le n$,  
$n-k+1=n-k'+s+t$, and  
\begin{gather}
a_1 =k, \qquad 
\overline{b_1} = \overline{k}, \qquad 
a'_{s'}=k', \qquad 
\overline{b'_{t'}}=\overline{k'}, \qquad 
a \not\preceq k', \qquad 
b \not\succeq \overline{k'}, 
\\
a_{i+1} \preceq b'_i, \quad (1 \le i < s), \qquad
\overline{b_{i+1}} \succeq \overline{a'_i}, \quad 
(1 \le i < t),
\end{gather}
where  
$a'_1 \prec \dots \prec a'_{s'}$ ($s':=t$) and  
$b'_1 \prec \dots \prec b'_{t'}$ ($t':=s$) are defined by 
\begin{equation}\label{eq:dual-type-two}
\begin{gathered}
\{a_1, \dots, a_s \} \sqcup \{ a'_1, \dots, a'_{s'} \}
 = \{ k, k+1, \dots, k'\}, \\
\{ \overline{b_1}, \dots, \overline{b_{t}} \} \sqcup 
\{ \overline{b'_1}, \dots, \overline{b'_{t'}} \}
 = \{ \overline{k}, \overline{k+1}, \dots,  \overline{k'} \}. 
\end{gathered}
\end{equation}
See Figure \ref{fig:tableau-path-correspondence-two} 
for the corresponding part in the paths. 

We say that an $L$-configuration $L=C(a_1, \dots, a_m)$ 
in the $j$th column of $T$ is {\it boundary} if 
$a_1=T(\mu'_j+1,j)$, i.e., if $a_1$ is at the top of the 
$j$th column, and $m$ is the largest number such that 
$L \cap L'=\emptyset$ for any $LU$-configuration $(L',U')$. 
Similarly, a $U$-configuration $U=C(a_1, \dots, a_m)$ 
in the $j$th column of $T$ is {\it boundary} if 
$a_m=T(\lambda'_j,j)$, i.e., if $a_m$ is at the bottom of the 
$j$th column, and $m$ is the largest number such that 
$U \cap U'=\emptyset$ for any $LU$-configuration $(L',U')$.

Let $(L,U)=(C(a_1, \dots, a_s), 
C(\overline{b_t}, \dots, \overline{b_1}))$ be 
an LU-configuration, and set 
$a'_1\prec \dots \prec a'_{s'}$ and 
$b'_{1}\prec \dots \prec b'_{t'}$ as in 
\eqref{eq:dual-type-one} (resp.~as in \eqref{eq:dual-type-two})
if $(L,U)$ is of type 1 (resp.~of type 2). 
We say that an L-configuration  
$L'$ is {\it right-adjacent} to $(L,U)$ if 
$L'$ is in the right-next column to $L$; furthermore, 
there exists some pair of an entry $e$ of $L'$ and an entry $a_i$ of $L$
such that $e$ is right-next to $a_i$ and $e\prec b'_i$. 
Similarly, 
we say that a U-configuration $U'$
is {\it left-adjacent} to $(L,U)$ if 
$U'$ is in the left-next column to $U$; furthermore, 
there exists some pair of an entry $e$ of $U'$ and an entry $\overline{b_i}$ 
of $U$ such that $e$ is left-next to $\overline{b_i}$ 
and $e\succ \overline{a'_i}$. 
Then, we say that an LU-configuration 
$(L',U')$ is {\it adjacent} to $(L,U)$ if one of the following conditions
is satisfied, and write it by $(L, U)\diamond (L',U')$: 
\begin{enumerate}[(i)]
\item $L'$ is right-adjacent to $(L, U)$.  
\item $L$ is right-adjacent to $(L', U')$.  
\item $U'$ is left-adjacent to $(L, U)$.  
\item $U$ is left-adjacent to $(L', U')$.  
\end{enumerate}

For any tableau $T$, let $\cLU(T)$ be 
the set of all LU-configurations of $T$. 
Then, the adjacent relation $\diamond$ of the LU-configurations 
generates an equivalence relation $\sim$ in $\cLU(T)$. 

\begin{defn}
For any $(L,U) \in \cLU(T)$, 
let $[(L,U)]\subset \cLU(T)$ be the equivalence class of $(L, U)$ 
with respect to $\sim$, 
and let $R=R(L,U):=\bigcup_{(L',U') \in [(L,U)]}(L',U')$ 
be the corresponding configuration in $T$. 
We call $R$  a {\it $\II$-region} of $T$, if 
the following conditions is satisfied: 
\begin{enumerate}
\item
No boundary L-configuration $L$ 
is right-adjacent to $L'$ for any 
LU-configuration $(L',U')$ in $R$. 
\item
No boundary U-configuration 
$U$ is left-adjacent to $U'$ for any 
LU-configuration $(L',U')$ in $R$. 
\end{enumerate}
Moreover, 
we say that $R$ is {\it odd} if 
the number of the odd type 1 LU-configurations 
in $R$ is odd.
\end{defn}
  
Then, an odd $\II$-region of $T=\cTv(\bp)$ 
corresponds to an odd $\II$-region of $\bp$, and therefore, 
Theorem \ref{thm:tableau-description} is rewritten as follows: 

\begin{thm}\label{thm:configuration}
For any skew diagram $\lambda/\mu$ satisfying
the positivity condition \eqref{eq:positivity},
we have 
$$\chi_{\lambda/\mu, a} = \sum_{T \in \Tab(\lambda/\mu)} z_a^{T}, $$
where $\Tab(\lambda/\mu)$ is 
the set of all the tableaux of shape $\lambda/\mu$ which satisfy 
$(\bH)$, $(\bV)$, 
and the following extra rule $(\bE')$:
\vspace*{5pt}

\noindent
\begin{tabular}{ll}
$(\bE')$ & $T$ does not have any odd $\II$-region. 
\end{tabular}
\vspace{5pt}
\end{thm}

\subsection{Explicit list of odd $\II$-regions}

Let us give an explicit list of all the possible odd $\II$-regions 
for $\lambda/\mu$ of at most two columns or of at most 
three rows. 

\begin{exmp}
Let $\lambda/\mu$ be a skew diagram of two columns 
satisfying the positivity condition \eqref{eq:positivity}.
In this case, 
an odd $\II$-region of $T \in \HVTab(\lambda/\mu)$ is nothing but 
an odd type 1 LU-configuration without  
any boundary L-configuration $L$  
which is right-adjacent to it, or 
any boundary U-configuration $U$ 
which is left-adjacent to it. 
Therefore, 
the extra rule $(\bE')$ is given as follows: 

\smallskip
({\bf E-2C}) \quad $T$ does not include any 
odd type 1 LU-configuration as
\begin{center}
{\setlength{\unitlength}{0.2mm}
\begin{picture}(50,220)
{\thicklines
\put(0,80){\line(0,1){140}}
\put(25,0){\line(0,1){220}}
\put(50,0){\line(0,1){140}}
\put(0,80){\line(1,0){25}}
\put(0,105){\line(1,0){25}}
\put(0,195){\line(1,0){25}}
\put(0,220){\line(1,0){25}}
\put(25,0){\line(1,0){25}}
\put(25,25){\line(1,0){25}}
\put(25,115){\line(1,0){25}}
\put(25,140){\line(1,0){25}}
}
\put(4,200){$a_1$}
\put(10,140){$\vdots$}
\put(4,87){$a_s$}
\put(0,60){$\overline{d_{s'}}$}
\put(10,30){$\vdots$}
\put(2,5){$\overline{d_{1}}$}
\put(29,5){$\overline{b_1}$}
\put(35,65){$\vdots$}
\put(29,120){$\overline{b_t}$}
\put(32,146){$c_{t'}$}
\put(35,169){$\vdots$}
\put(33,200){$c_{1}$}
\multiput(50,0)(5,0){10}{\line(1,0){2}}
\multiput(25,220)(5,0){15}{\line(1,0){2}}
\put(80,0){\vector(0,1){220}}
\put(80,220){\vector(0,-1){220}}
\put(85,100){$\scriptstyle n-k+1$}
\multiput(0,80)(-5,0){10}{\line(-1,0){2}}
\multiput(0,140)(-5,0){10}{\line(-1,0){2}}
\put(-30,80){\vector(0,1){60}}
\put(-30,140){\vector(0,-1){60}}
\put(-45,100){$\scriptstyle r$}
\end{picture}
}
\end{center} 
where, if $t'=t-r\ge 1$, then $c_{t'} \preceq n$ and 
$c_i \succeq b'_i$ for any $i=1, \dots, t'$, 
and if $s'=s-r \ge 1$, then $d_{s'}\succeq \overline{n}$ and 
$\overline{d_i} \preceq \overline{a'_i}$ for any $i=1, \dots, s'$. 
\smallskip
\end{exmp}

\begin{exmp}\label{exmp:one-row}
Let $\lambda$ be a Young diagram of one row.  
In this case, the odd $\II$-region of $T \in \HVTab(\lambda)$ 
is the configuration 
{\setlength{\unitlength}{0.2mm}
\begin{picture}(50,25)(0,6)
\multiput(0,0)(25,0){3}{\line(0,1){25}}
\multiput(0,0)(0,25){2}{\line(1,0){50}}
\put(6,6){$n$}
\put(31,6){$\overline{n}$}
\end{picture}
}, and therefore, 
the extra rule $(\bE')$ is given as follows: 

\smallskip
({\bf E-1R}) \quad $T$ does not include 
{\setlength{\unitlength}{0.2mm}
\begin{picture}(50,25)(0,6)
\multiput(0,0)(25,0){3}{\line(0,1){25}}
\multiput(0,0)(0,25){2}{\line(1,0){50}}
\put(6,6){$n$}
\put(31,6){$\overline{n}$}
\end{picture}
}.
\smallskip
\end{exmp}

\begin{rem}\label{rem:E-one}
The extra rule ({\bf E-1R}) is applied 
for any $\lambda/\mu$, since 
{\setlength{\unitlength}{0.2mm}
\begin{picture}(50,25)(0,6)
\multiput(0,0)(25,0){3}{\line(0,1){25}}
\multiput(0,0)(0,25){2}{\line(1,0){50}}
\put(6,6){$n$}
\put(31,6){$\overline{n}$}
\end{picture}
}
is an odd $\II$-region of $T\in \HVTab(\lambda/\mu)$. 
\end{rem}

\begin{table}
\begin{center}
\begin{tabular}{p{3.8cm}|p{8.5cm}}
\noalign{\hrule height0.8pt}
\shortstack{odd type 1\\
LU-configurations}
\shortstack{($a\succeq \overline{n}$, $b \preceq n$)} & 
{\setlength{\unitlength}{0.4mm}
\begin{picture}(32,13)
\multiput(0,0)(16,0){3}{\line(0,1){13}}
\multiput(0,0)(0,13){2}{\line(1,0){32}}
\put(6,4){$\scriptstyle n$}
\put(22,4){$\scriptstyle \overline{n}$}
\end{picture}
\quad 
\begin{picture}(32,26)
\put(0,0){\line(1,0){32}}
\put(0,26){\line(1,0){16}}
\put(0,13){\line(1,0){32}}
\put(0,0){\line(0,1){26}}
\put(16,0){\line(0,1){26}}
\put(32,0){\line(0,1){13}}
\text{
\put(6,3){$\scriptstyle n$}
\put(1,16){$\scriptstyle n-1$}
\put(17,3){$\scriptstyle \overline{n-1}$}
\put(22,16){$\scriptstyle b$}
}
\end{picture}
\quad 
\begin{picture}(32,26)
\put(0,26){\line(1,0){32}}
\put(0,13){\line(1,0){32}}
\put(16,0){\line(1,0){16}}
\put(0,13){\line(0,1){13}}
\put(16,0){\line(0,1){26}}
\put(32,0){\line(0,1){26}}
\text{
\put(6,3){$\scriptstyle a$}
\put(1,16){$\scriptstyle n-1$}
\put(17,3){$\scriptstyle \overline{n-1}$}
\put(22,16){$\scriptstyle \overline{n}$}
}
\end{picture}
\quad 
\begin{picture}(32,45)
\put(0,39){\line(1,0){32}}
\put(0,26){\line(1,0){32}}
\put(16,13){\line(1,0){16}}
\put(16,0){\line(1,0){16}}
\put(0,26){\line(0,1){13}}
\put(16,0){\line(0,1){39}}
\put(32,0){\line(0,1){39}}
\text{
\put(6,16){$\scriptstyle a$}
\put(1,29){$\scriptstyle n-2$}
\put(17,3){$\scriptstyle \overline{n-2}$}
\put(17,16){$\scriptstyle \overline{n-1}$}
\put(22,29){$\scriptstyle \overline{n}$}
}
\end{picture}
\quad 
\begin{picture}(32,45)
\put(0,39){\line(1,0){16}}
\put(0,26){\line(1,0){32}}
\put(0,13){\line(1,0){32}}
\put(16,0){\line(1,0){16}}
\put(0,13){\line(0,1){26}}
\put(16,0){\line(0,1){39}}
\put(32,0){\line(0,1){26}}
\text{
\put(6,3){$\scriptstyle a$}
\put(1,16){$\scriptstyle n-1$}
\put(1,29){$\scriptstyle n-2$}
\put(17,3){$\scriptstyle \overline{n-2}$}
\put(17,16){$\scriptstyle \overline{n-1}$}
\put(22,29){$\scriptstyle b$}
}
\end{picture}
\quad 
\begin{picture}(32,45)
\put(0,39){\line(1,0){16}}
\put(0,26){\line(1,0){32}}
\put(0,13){\line(1,0){32}}
\put(16,0){\line(1,0){16}}
\put(0,13){\line(0,1){26}}
\put(16,0){\line(0,1){39}}
\put(32,0){\line(0,1){26}}
\text{
\put(6,3){$\scriptstyle a$}
\put(1,16){$\scriptstyle n-1$}
\put(1,29){$\scriptstyle n-2$}
\put(17,3){$\scriptstyle \overline{n-2}$}
\put(22,16){$\scriptstyle \overline{n}$}
}
\end{picture}
\quad 
\begin{picture}(32,45)
\put(0,39){\line(1,0){16}}
\put(0,26){\line(1,0){32}}
\put(0,13){\line(1,0){32}}
\put(16,0){\line(1,0){16}}
\put(0,13){\line(0,1){26}}
\put(16,0){\line(0,1){39}}
\put(32,0){\line(0,1){26}}
\text{
\put(6,16){$\scriptstyle n$}
\put(1,29){$\scriptstyle n-2$}
\put(17,3){$\scriptstyle \overline{n-2}$}
\put(17,16){$\scriptstyle \overline{n-1}$}
\put(22,29){$\scriptstyle b$}
}
\end{picture}
\quad 
\begin{picture}(32,45)
\put(0,39){\line(1,0){16}}
\put(0,26){\line(1,0){16}}
\put(0,13){\line(1,0){32}}
\put(0,0){\line(1,0){32}}
\put(0,0){\line(0,1){39}}
\put(16,0){\line(0,1){39}}
\put(32,0){\line(0,1){13}}
\text{
\put(6,3){$\scriptstyle n$}
\put(1,16){$\scriptstyle n-1$}
\put(1,29){$\scriptstyle n-2$}
\put(17,3){$\scriptstyle \overline{n-2}$}
\put(22,16){$\scriptstyle b$}
}
\end{picture}
\quad 
\begin{picture}(32,45)
\put(0,39){\line(1,0){32}}
\put(0,26){\line(1,0){32}}
\put(0,13){\line(1,0){32}}
\put(0,0){\line(1,0){32}}
\put(0,0){\line(0,1){39}}
\put(16,0){\line(0,1){39}}
\put(32,0){\line(0,1){39}}
\text{
\put(6,3){$\scriptstyle n$}
\put(1,16){$\scriptstyle n-1$}
\put(1,29){$\scriptstyle n-2$}
\put(17,3){$\scriptstyle \overline{n-2}$}
\put(17,16){$\scriptstyle \overline{n-1}$}
\put(22,29){$\scriptstyle \overline{n}$}
}
\end{picture}
}
\\
\hline 
\shortstack{other\\
LU-configurations}
\shortstack{
($a\succeq \overline{n}$, $b \preceq n$,\\
$a'\succ n-1$,
$b' \prec \overline{n-1}$)} &
{\setlength{\unitlength}{0.4mm}
\begin{picture}(32,26)
\put(0,26){\line(1,0){32}}
\put(0,13){\line(1,0){32}}
\put(0,0){\line(1,0){32}}
\put(0,0){\line(0,1){26}}
\put(16,0){\line(0,1){26}}
\put(32,0){\line(0,1){26}}
\text{
\put(6,3){$\scriptstyle n$}
\put(1,16){$\scriptstyle n-1$}
\put(17,3){$\scriptstyle \overline{n-1}$}
\put(22,16){$\scriptstyle \overline{n}$}
}
\end{picture}
\quad 
\begin{picture}(32,26)
\put(0,26){\line(1,0){16}}
\put(0,13){\line(1,0){32}}
\put(16,0){\line(1,0){16}}
\put(0,13){\line(0,1){13}}
\put(16,0){\line(0,1){26}}
\put(32,0){\line(0,1){13}}
\text{
\put(6,3){$\scriptstyle a$}
\put(1,16){$\scriptstyle n-1$}
\put(17,3){$\scriptstyle \overline{n-1}$}
\put(22,16){$\scriptstyle b$}
} 
\end{picture}
\quad 
\begin{picture}(32,45)
\put(0,39){\line(1,0){32}}
\put(0,26){\line(1,0){32}}
\put(0,13){\line(1,0){32}}
\put(16,0){\line(1,0){16}}
\put(0,13){\line(0,1){26}}
\put(16,0){\line(0,1){39}}
\put(32,0){\line(0,1){39}}
\text{
\put(6,3){$\scriptstyle a$}
\put(1,16){$\scriptstyle n-1$}
\put(1,29){$\scriptstyle n-2$}
\put(17,3){$\scriptstyle \overline{n-2}$}
\put(17,16){$\scriptstyle \overline{n-1}$}
\put(22,29){$\scriptstyle \overline{n}$}
}
\end{picture}
\quad   
\begin{picture}(32,45)
\put(0,39){\line(1,0){32}}
\put(0,26){\line(1,0){32}}
\put(0,13){\line(1,0){32}}
\put(16,0){\line(1,0){16}}
\put(0,13){\line(0,1){26}}
\put(16,0){\line(0,1){39}}
\put(32,0){\line(0,1){39}}
\text{
\put(6,16){$\scriptstyle n$}
\put(1,29){$\scriptstyle n-2$}
\put(17,3){$\scriptstyle \overline{n-2}$}
\put(17,16){$\scriptstyle \overline{n-1}$}
\put(22,29){$\scriptstyle \overline{n}$}
}
\end{picture}
\quad    
\begin{picture}(32,45)
\put(0,39){\line(1,0){16}}
\put(0,26){\line(1,0){32}}
\put(0,13){\line(1,0){32}}
\put(0,0){\line(1,0){32}}
\put(0,0){\line(0,1){39}}
\put(16,0){\line(0,1){39}}
\put(32,0){\line(0,1){26}}
\text{
\put(6,3){$\scriptstyle n$}
\put(1,16){$\scriptstyle n-1$}
\put(1,29){$\scriptstyle n-2$}
\put(17,3){$\scriptstyle \overline{n-2}$}
\put(17,16){$\scriptstyle \overline{n-1}$}
\put(22,29){$\scriptstyle b$}
}
\end{picture}
\quad    
\begin{picture}(32,45)
\put(0,39){\line(1,0){16}}
\put(0,26){\line(1,0){32}}
\put(0,13){\line(1,0){32}}
\put(0,0){\line(1,0){32}}
\put(0,0){\line(0,1){39}}
\put(16,0){\line(0,1){39}}
\put(32,0){\line(0,1){26}}
\text{
\put(6,3){$\scriptstyle n$}
\put(1,16){$\scriptstyle n-1$}
\put(1,29){$\scriptstyle n-2$}
\put(17,3){$\scriptstyle \overline{n-2}$}
\put(22,16){$\scriptstyle \overline{n}$}
}
\end{picture}
\quad  
\begin{picture}(32,45)
\put(0,39){\line(1,0){16}}
\put(0,26){\line(1,0){32}}
\put(16,13){\line(1,0){16}}
\put(16,0){\line(1,0){16}}
\put(0,26){\line(0,1){13}}
\put(16,0){\line(0,1){39}}
\put(32,0){\line(0,1){26}}
\text{
\put(6,16){$\scriptstyle a$}
\put(1,29){$\scriptstyle n-2$}
\put(17,3){$\scriptstyle \overline{n-2}$}
\put(17,16){$\scriptstyle \overline{n-1}$}
\put(22,29){$\scriptstyle b$}
}
\end{picture}
\quad    
\begin{picture}(32,45)
\put(0,39){\line(1,0){16}}
\put(0,26){\line(1,0){16}}
\put(0,13){\line(1,0){32}}
\put(16,0){\line(1,0){16}}
\put(0,13){\line(0,1){26}}
\put(16,0){\line(0,1){39}}
\put(32,0){\line(0,1){13}}
\text{
\put(6,3){$\scriptstyle a$}
\put(1,16){$\scriptstyle n-1$}
\put(1,29){$\scriptstyle n-2$}
\put(17,3){$\scriptstyle \overline{n-2}$}
\put(22,16){$\scriptstyle b$}
}
\end{picture}
\quad  
\begin{picture}(32,45)
\put(0,39){\line(1,0){16}}
\put(0,26){\line(1,0){16}}
\put(16,13){\line(1,0){16}}
\put(16,0){\line(1,0){16}}
\put(0,26){\line(0,1){13}}
\multiput(16,0)(0,4){10}{\line(0,1){2}}
\put(16,0){\line(0,1){13}}
\put(16,26){\line(0,1){13}}
\put(32,0){\line(0,1){13}}
\text{
\put(6,16){$\scriptstyle a'$}
\put(1,29){$\scriptstyle n-2$}
\put(17,3){$\scriptstyle \overline{n-2}$}
\put(22,16){$\scriptstyle b'$}
}
\end{picture}
}
\\
\noalign{\hrule height0.8pt}
\end{tabular}
\end{center}
\caption{The LU-configurations in a tableau in $\Tab(\lambda/\mu)$ 
for a skew diagram $\lambda/\mu$ of at most three rows. 
We omit $a$ and $b$ if the inequalities are satisfied 
by the vertical rule $(\bV)$. 
}\label{tab:LU-configurations}
\end{table}

\begin{exmp}
Let $\lambda/\mu$ be a skew diagram of two rows, i.e., $\lambda'_1 = 2$. 
Since an odd $\II$-region of $T \in \HVTab(\lambda/\mu)$ is a combination
of the LU-configurations of at most two rows 
in Table \ref{tab:LU-configurations}
with one odd type 1 LU-configuration. Therefore, 
the extra rule $(\bE')$ 
is given by 
({\bf E-1R}) and the following condition: 

\smallskip
({\bf E-2R}) \quad 
$T$ does not include 
\begin{center}
\begin{psfrags}
\psfrag{n}{$\scriptstyle n$}
\psfrag{n-1}{$\scriptstyle n-1$}
\psfrag{bn}{$\scriptstyle \overline{n}$}
\psfrag{bn-1}{$\scriptstyle \overline{n-1}$}
\psfrag{d}{$\cdots$}
\psfrag{a}{$\scriptstyle \overline{n}$}
\psfrag{b}{$\scriptstyle n$}
\psfrag{j}{$\scriptstyle j$}
\includegraphics[height=1.9cm, clip]{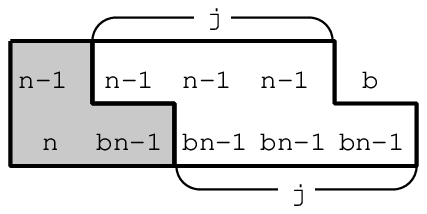}
\ , \qquad 
\includegraphics[height=1.9cm, clip]{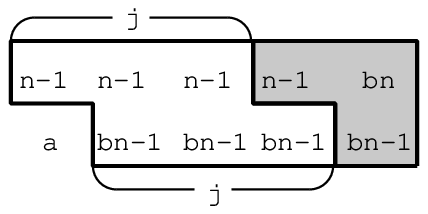}
\ , 
\end{psfrags}
\end{center}
where $j\ge 0$. 

This rule appears for the proposed crystals of 
the representations of two-row rectangle \cite{SS}.
\end{exmp}

\begin{exmp}
Let $\lambda/\mu$ be a skew diagram of three rows, i.e., $\lambda'_1 = 3$.   
Since an odd $\II$-region of $T \in \HVTab(\lambda/\mu)$ is a combination 
of the LU-configurations of at most three rows 
in Table \ref{tab:LU-configurations}
with one or three odd type 1 LU-configurations, 
the extra rule $(\bE')$ 
for a tableau $T \in \HVTab(\lambda/\mu)$ is given by 
({\bf E-1R}), ({\bf E-2R}), and the following condition: 

\smallskip
({\bf E-3R}) \quad 
$T$ does not include 
\begin{equation}\label{eq:three-one}
\text{
\begin{psfrags}
\psfrag{n}{$\scriptstyle n$}
\psfrag{n-1}{$\scriptstyle n-1$}
\psfrag{n-2}{$\scriptstyle n-2$}
\psfrag{bn}{$\scriptstyle \overline{n}$}
\psfrag{bn-1}{$\scriptstyle \overline{n-1}$}
\psfrag{bn-2}{$\scriptstyle \overline{n-2}$}
\psfrag{d}{$\cdots$}
\psfrag{a}{$\scriptstyle a$}
\psfrag{b}{$\scriptstyle b$}
\psfrag{j1}{$\scriptstyle j_1$}
\psfrag{j2}{$\scriptstyle j_2$}
\psfrag{j3}{$\scriptstyle j_3$}
\psfrag{j4}{$\scriptstyle j_4$}
\psfrag{j5}{$\scriptstyle j_5$}
\includegraphics[height=2.5cm, clip]{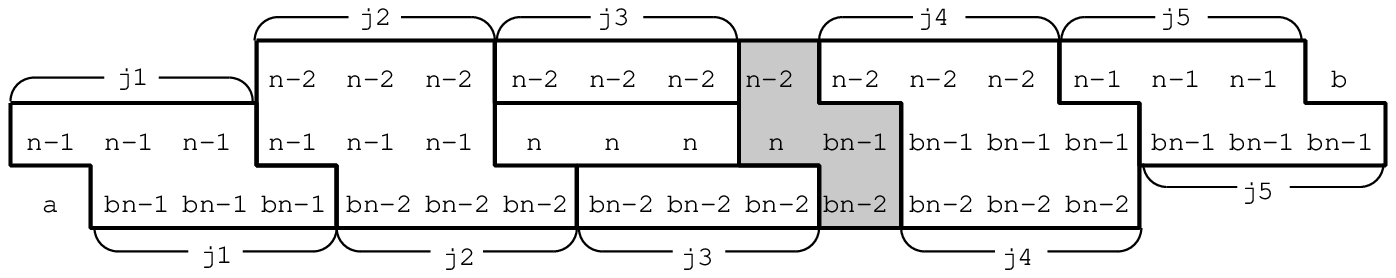}
\end{psfrags}
}
\end{equation}
where $j_i\ge 0$, $j_1=0$ if $j_2=0$, and 
\begin{equation*}
(a,b)=
\begin{cases}
(\overline{n}, n), (n, \overline{n}), & \text{if $j_2 \ne 0$}, \\
(\overline{n},n),  (\overline{n-1},n) & \text{if $j_2=0$},\\
\end{cases}
\end{equation*}
\begin{equation}\label{eq:three-two}
\text{
\begin{psfrags}
\psfrag{n}{$\scriptstyle n$}
\psfrag{n-1}{$\scriptstyle n-1$}
\psfrag{n-2}{$\scriptstyle n-2$}
\psfrag{bn}{$\scriptstyle \overline{n}$}
\psfrag{bn-1}{$\scriptstyle \overline{n-1}$}
\psfrag{bn-2}{$\scriptstyle \overline{n-2}$}
\psfrag{d}{$\cdots$}
\psfrag{a}{$\scriptstyle a$}
\psfrag{b}{$\scriptstyle b$}
\psfrag{j1}{$\scriptstyle j_1$}
\psfrag{j2}{$\scriptstyle j_2$}
\psfrag{j3}{$\scriptstyle j_3$}
\psfrag{j4}{$\scriptstyle j_4$}
\psfrag{j5}{$\scriptstyle j_5$}
\includegraphics[height=2.5cm, clip]{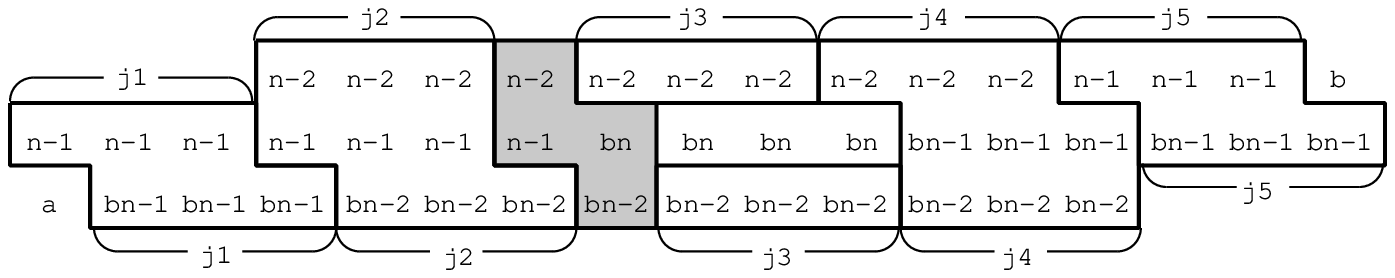}
\end{psfrags}
}
\end{equation}
where $j_i\ge 0$, $j_5=0$ if $j_4=0$, and 
\begin{equation*}
(a,b)=
\begin{cases}
(\overline{n}, n), (n, \overline{n}), & \text{if $j_4\ne 0$}, \\
(\overline{n},n-1), (\overline{n},n), & \text{if $j_4=0$}, \\
\end{cases}
\end{equation*}
and 
\begin{equation}\label{eq:three-three}
\text{
\begin{psfrags}
\psfrag{n}{$\scriptstyle n$}
\psfrag{n-1}{$\scriptstyle n-1$}
\psfrag{n-2}{$\scriptstyle n-2$}
\psfrag{bn}{$\scriptstyle \overline{n}$}
\psfrag{bn-1}{$\scriptstyle \overline{n-1}$}
\psfrag{bn-2}{$\scriptstyle \overline{n-2}$}
\psfrag{d}{$\cdots$}
\psfrag{a}{$\scriptstyle a$}
\psfrag{b}{$\scriptstyle b$}
\psfrag{j1}{$\scriptstyle j_1$}
\psfrag{j2}{$\scriptstyle j_2$}
\psfrag{j4}{$\scriptstyle j_4$}
\psfrag{j5}{$\scriptstyle j_5$}
\includegraphics[height=2.5cm, clip]{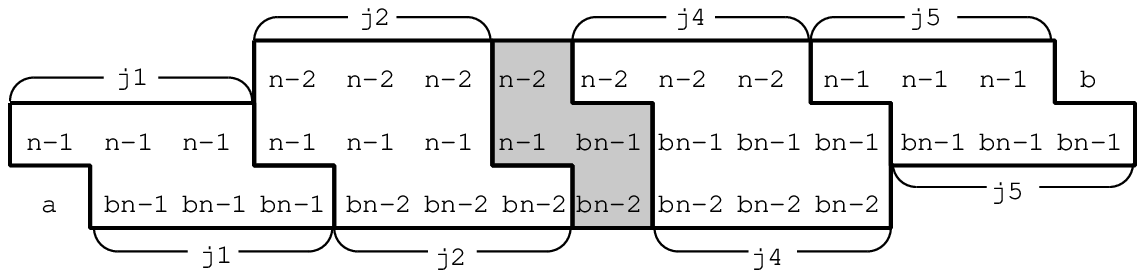}
\end{psfrags}
}
\end{equation}
where $j_i \ge 0$ and $(a,b)=(n,\overline{n})$ or $(\overline{n},n)$. 
\smallskip

The odd $\II$-regions \eqref{eq:three-one}, \eqref{eq:three-two}, and 
\eqref{eq:three-three} include three (resp.~one)
odd type 1 LU-configurations 
if $(a,b)=(n, \overline{n})$, $j_1+j_2\ge 1$, 
and $j_4+j_5 \ge 1$ are satisfied (resp.~otherwise). 
\end{exmp}

\begin{figure}
\begin{psfrags}
\psfrag{j1-1}{\small $\scriptscriptstyle j_1-1$}
\psfrag{j5-1}{\small $\scriptscriptstyle j_5-1$}
\psfrag{j2}{\small $\scriptscriptstyle j_2$}
\psfrag{j4}{\small $\scriptscriptstyle j_4$}
\psfrag{dots}{$\dots$}
\includegraphics[width=12.5cm, clip]{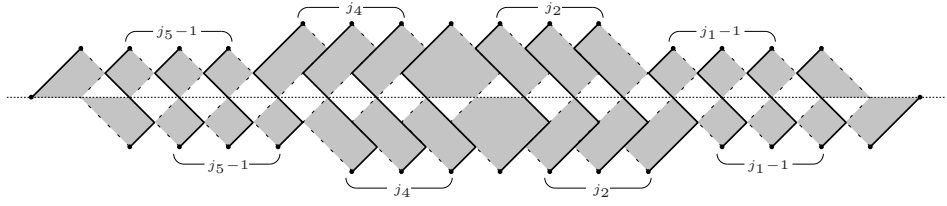}
\end{psfrags}
\caption{The odd $\II$-region corresponding to 
the rule \eqref{eq:three-three} for 
$(a,b)=(n,\overline{n})$, $j_1, j_5 \ne 0$. }\label{fig:odd-exmp}
\end{figure}

As mentioned in the introduction, 
though these rules look rather complicated, 
they are more easily recognizable in the path picture. 
For example, the rule \eqref{eq:three-three} 
for $(a,b)=(n, \overline{n})$, $j_1, j_5 \ne 0$ corresponds 
to the odd $\II$-region in Figure \ref{fig:odd-exmp}. 

\section{Construction of $\phi$ and proof of Proposition 
\ref{prop:folding-map}}\label{sec:folding-map}
In this section, we construct the folding map 
$\phi: P_2(\lambda/\mu) \to P(\lambda/\mu)$ 
in Proposition \ref{prop:folding-map}, 
which is a key to derive the tableaux description 
(see Theorem \ref{thm:tableau-description}).

\subsection{$k$-expansion and $k$-folding}
To define $\phi$, we introduce the {\it $k$-expansion} 
and the {\it $k$-folding}, which are generalizations of the  
the expansion and the folding in Section \ref{sec:expansion-folding}. 
The original corresponds to $k=1$. 
We also generalize related notions.  

\begin{defn}\label{def:k-unit}
Let 
$(\alpha; \beta) \in \cH(\lambda/\mu)$. 
For any unit $U \subset S_{\pm}$, let $\pm r= \Ht(U)$ and 
let $a$ and $a'=a+1$ be the horizontal positions of the 
left and the right vertices of $U$. 
Then, for any $k=1 ,\dots , l-1$,
\begin{enumerate}
\item 
$U$ is called a $\I_k$-{\it unit} of $(\alpha; \beta)$ 
if there exists some $i$ ($0\le i \le l$) such that 
\begin{equation}\label{eq:Ik-unit}
\begin{aligned}
& \alpha^*_i(r) \le a < a' \le \beta_{i+k}(r), \quad \text{if $U \subset S_+$}, \\
& \alpha_i(-r) \le a < a' \le \beta_{i+k}^*(-r), \quad \text{if $U \subset S_-$}. 
\end{aligned}
\end{equation}
\item 
$U$ is called a $\II_k$-{\it unit} of $(\alpha; \beta)$ 
if there exists some $i$ ($0 \le i \le l$) such that 
\begin{equation}\label{eq:IIk-unit}
\begin{aligned}
& \beta_{i+k}(r) \le a < a' \le \alpha^*_i(r), \quad \text{if $U \subset S_+$}, \\
& \beta^*_{i+k}(-r) \le a < a' \le \alpha_i(-r), \quad \text{if $U \subset S_-$}.  
\end{aligned}
\end{equation}
\end{enumerate}
Here, we set $\beta_i(r)=\beta^*_i(-r)=-\infty$ 
for $i \ge l+1$ 
and $\alpha_i(-r)=\alpha^*_i(r)=+\infty$  for 
$i \le 0$. 
Furthermore, a $\II_k$-unit $U$ of $(\alpha; \beta)$ is called 
a {\it boundary} $\II_k$-unit if \eqref{eq:IIk-unit} holds 
for $i = 0$, $i \ge l-k+1$, or $r=n$. 
\end{defn}

As in the $k=1$ case, actually \eqref{eq:Ik-unit} 
does not hold for $i = 0$ or $i \ge l-k+1$. Also 
it does not hold for $r=n$ if $\lambda/\mu$ satisfies the 
positivity condition. 

As Lemma \ref{lem:I-II-unit}, it is easy to see that 
\begin{lem}\label{lem:unit}
\begin{enumerate}
\item $U$ is a $\I_k$-unit (resp.~a $\II_k$-unit) 
if and only if $U^*$ is a $\I_k$-unit (resp.~a $\II_k$-unit). 
\item 
No unit is simultaneously 
a $\I_k$- and $\II_{k'}$-unit for any $k$ and $k'$ 
such that $k+1 \ge k'$. 
\item \label{item:not-connected}
If $U$ is a $\I_k$-unit and $U'$ is a $\II_{k'}$-unit 
for any $k \ge k'$, then $U$ and $U'$ are not adjacent. 
\item \label{item:complementary} 
The set of all the $\I_k$-units 
and that of all the $\II_{k+1}$-units are complementary to 
each other in $S_+\sqcup S_-$.  
\end{enumerate}
\end{lem}

We define $\cU_{\I_k}$ and $\tilde{\cU}_{\I_k}$
(resp.~$\cU_{\II_k}$ and $\tilde{\cU}_{\II_k}$)
for $\I_k$-units (resp.~$\II_k$-units) 
of $(\alpha; \beta)$ and a {\it connected component} 
of $\tilde{\cU}_{\I_k}$ and $\tilde{\cU}_{\II_k}$
similarly as in the $k=1$ case. 

\begin{defn} 
Let $\lambda/\mu$ be a skew diagram 
satisfying the positivity condition 
\eqref{eq:positivity}, 
and let $(\alpha; \beta)\in \cH(\lambda/\mu)$. 
\begin{enumerate}
\item 
A connected component $V$ of $\tilde{\cU}_{\I_k}$ 
is called a $\I_k$-{\it region} 
of $(\alpha; \beta)$ 
if it contains at least one $\I_k$-unit of height $0$. 
\item 
A connected component $V$ of $\tilde{\cU}_{\II_k}$ 
is called a $\II_k$-{\it region} 
of $(\alpha; \beta)$ if it satisfies the following conditions: 
\begin{enumerate}[(i)]
\item $V$ contains at least one $\I_k$-unit of height $0$. 
\item $V$ does not contain any boundary $\II_k$-unit. 
\end{enumerate}
\end{enumerate}
\end{defn}

As Proposition \ref{prop:duality}, we have 
\begin{prop}
If $V$ is a $\I_k$- or $\II_k$-region of $(\alpha; \beta)$, then 
$V^*=V$. 
\end{prop}

If $U$ is a $\I_k$-unit of $(\alpha; \beta)$, then 
it is also a $\I_{k'}$-unit $(\alpha; \beta)$ 
for any $1 \le k' \le k$, 
while if $U$ is a $\II_k$-unit of $(\alpha; \beta)$, 
then it is also a $\II_{k'}$-unit of $(\alpha; \beta)$ 
for any $k \le k' \le l-1$. Then, it follows that

\begin{lem}\label{lem:inclusion}
\begin{enumerate}
\item If $V$ is a $\I_k$-region and $V'$ is a $\I_{k'}$-region 
for $k' \ge k$, then $V' \subset V$ or $V \cap V '= \emptyset$.  
\item \label{item:inclusion-II}
If $V$ is a $\II_k$-region and $V'$ is a $\II_{k'}$-region
for $k' \le k$, then $V' \subset V$ or $V \cap V '= \emptyset$.  
\end{enumerate}
\end{lem}

For any 
$(\alpha; \beta)
\in \cH(\lambda/\mu)$, let $V$ be any $\I_k$- or 
$\II_k$-region of $(\alpha; \beta)$. 
Let $\alpha'_i$ be the lower path obtained from 
$\alpha_i$ by replacing the part 
$\alpha_i \cap V$ with $\beta^*_{i+k}\cap V$, and let 
$\beta'_i$ be the upper path obtained from 
$\beta_i$ by replacing the part
$\beta_i \cap V$ with $\alpha^*_{i-k}\cap V$. 
Set $\varepsilon_V^k(\alpha; \beta):= 
(\alpha'_1, \dots, \alpha'_l; \beta'_1, \dots, \beta'_l)$. 
Then, Proposition \ref{prop:I-II-region} 
is generalized as follows: 

\begin{prop}\label{prop:region}
Let $\lambda/\mu$ be a skew diagram 
satisfying the positivity condition \eqref{eq:positivity}. 
Then, for any $(\alpha; \beta) \in \cH(\lambda/\mu)$, 
we have 
\begin{enumerate}
\item For any $\I_k$- or $\II_k$-region $V$
of $(\alpha; \beta)$, $\varepsilon_V^k(\alpha; \beta) \in \cH(\lambda/\mu)$. 
\item \label{item:expansion-result}
If $V$ is a $\I_k$-region and 
$V'\subset V$ is a $\I_r$-region of $(\alpha; \beta)$ 
for $r \ge k$, then $V'$  is a $\II_{2k-r}$-region 
of $\varepsilon_V^k(\alpha; \beta)$. 
\item \label{item:folding-result-zero}
If $V$ is a $\II_k$-region and 
$V'\subset V$ is a $\II_r$-region of $(\alpha;\beta)$ 
for $r \le k$, then $V'$ is a $\I_{2k-r}$-region of 
$\varepsilon_V^k(\alpha; \beta)$.
\end{enumerate}
\end{prop}

Set $(\alpha'; \beta')=\varepsilon_V^k(\alpha; \beta)$ 
for any $\I_k$- or $\II_k$-region $V$. 
The following equalities are useful: 
\begin{equation}\label{eq:k-folding}
\begin{aligned}
\alpha'_{i}(0) & =
\begin{cases}
\beta^*_{i+k}(0), & \text{if $\beta_{i+r}$ intersects with $V$ at hight 0}, \\
\alpha_i(0), & \text{otherwise}, 
\end{cases}
\\
\beta'_{i}(0) & =
\begin{cases}
\alpha^*_{i-k}(0),  & \text{if $\alpha_{i-k}$ intersects with $V$ at height 0}, \\
\beta_i(0), & \text{otherwise}.
\end{cases} 
\end{aligned}
\end{equation}

We call the correspondence 
$(\alpha; \beta) \mapsto \varepsilon_V^k(\alpha; \beta)$ the 
{\it $k$-expansion} (resp.~the {\it $k$-folding}) 
with respect to $V$, if $V$ is a $\I_k$-region 
(resp. a $\II_k$-region).  
As the $k=1$ case, 
we have $\varepsilon^k_V \circ \varepsilon^k_V = \id$ 
for any $V$. 

Let $(\alpha; \beta) \in \cH(\lambda/\mu)$. 
If $h:=\alpha_i(0) - \beta_{i+k}(0)$ is a 
non-positive number (resp.~a positive number), 
then we call a pair $(\alpha_i, \beta_{i+k})$ 
a {\it $k$-overlap}  
(resp.~a {\it $k$-hole}).  
Furthermore, if $h$ is an even number (resp.~an odd number), 
then $(\alpha_i, \beta_{i+k})$ is called {\it even} 
(resp.~{\it odd}). 

For any $\I_k$-region $V$ (resp.~$\II_k$-region $V$) of 
$(\alpha;\beta) \in \cH(\lambda/\mu)$, 
let $n(V)$ be the number of 
the even $k$-overlaps 
(resp.~the even $k$-holes) which 
intersect with $V$ at height $0$. 
\begin{defn}
A $\I_k$- or $\II_k$-region $V$ is called {\it even} 
(resp.~{\it odd}) if $n(V)$ is even (resp.~odd). 
\end{defn}

\subsection{Outline of construction}
In this subsection, 
we give the outline of the construction of  
the folding map $\phi$ whose existence is admitted in 
Section \ref{sec:tableau-description}. 
Proofs of Propositions \ref{prop:one}--\ref{prop:four}
below will be given in the following sections. 

\begin{figure}
\begin{psfrags}
\psfrag{0}{$\scriptstyle 0$}
\psfrag{1}{$\scriptstyle 1$}
\psfrag{2}{$\scriptstyle 2$}
\psfrag{3}{$\scriptstyle 3$}
\psfrag{4}{$\scriptstyle 4$}
\psfrag{5}{$\scriptstyle 5$}
\psfrag{6}{$\scriptstyle 6$}
\psfrag{7}{$\scriptstyle 7$}
\psfrag{9}{$\scriptstyle 9$}
\psfrag{12}{$\scriptstyle 12$}
\psfrag{14}{$\scriptstyle 14$}
\psfrag{17}{$\scriptstyle 17$}
\psfrag{19}{$\scriptstyle 19$}
\psfrag{20}{$\scriptstyle 20$}
\psfrag{21}{$\scriptstyle 21$}
\psfrag{23}{$\scriptstyle 23$}
\psfrag{24}{$\scriptstyle 24$}
\psfrag{25}{$\scriptstyle 25$}
\psfrag{26}{$\scriptstyle 26$}
\psfrag{27}{$\scriptstyle 27$}
\psfrag{28}{$\scriptstyle 28$}
\psfrag{29}{$\scriptstyle 29$}
\psfrag{30}{$\scriptstyle 30$}
\psfrag{32}{$\scriptstyle 32$}
\psfrag{a1}{\small $\alpha_1$}
\psfrag{a2}{\small $\alpha_2$}
\psfrag{a3}{\small $\alpha_3$}
\psfrag{a4}{\small $\alpha_4$}
\psfrag{a5}{\small $\alpha_5$}
\psfrag{a6}{\small $\alpha_6$}
\psfrag{a7}{\small $\alpha_7$}
\psfrag{a8}{\small $\alpha_8$}
\psfrag{sa}{\small $s(\alpha_i)$}
\psfrag{ma}{\small $m_1(\alpha_i)$}
\psfrag{b1}{\small $\beta_1$}
\psfrag{b2}{\small $\beta_2$}
\psfrag{b3}{\small $\beta_3$}
\psfrag{b4}{\small $\beta_4$}
\psfrag{b5}{\small $\beta_5$}
\psfrag{b6}{\small $\beta_6$}
\psfrag{b7}{\small $\beta_7$}
\psfrag{b8}{\small $\beta_8$}
\psfrag{sb}{\small $s(\beta_i)$}
\psfrag{mb}{\small $m_1(\beta_i)$}
\psfrag{p1}{$p_1$}
\psfrag{p2}{$p_2$}
\psfrag{p3}{$p_3$}
\psfrag{p4}{$p_4$}
\psfrag{p5}{$p_5$}
\psfrag{p6}{$p_6$}
\psfrag{p7}{$p_7$}
\psfrag{p8}{$p_8$}
\psfrag{dots}{$\cdots$}
\psfrag{pi}{$\pi$}
\psfrag{pi-1}{$\pi^{_1}$}
\psfrag{height}{\small height}
\vspace*{12pt}
\includegraphics[width=12.5cm, clip]{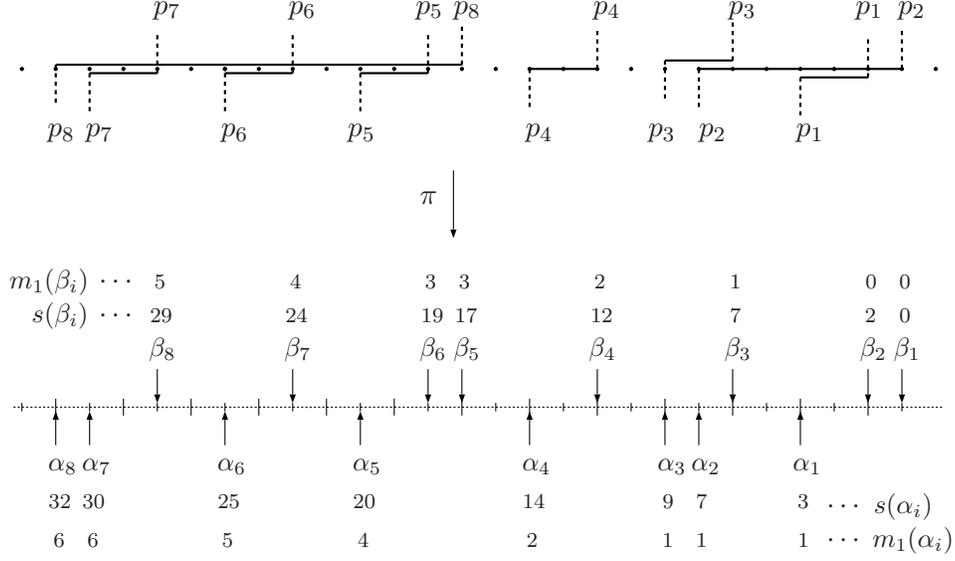}
\end{psfrags}
\caption{An example of  
$s(\alpha_i)$, $s(\beta_i)$,  $m(\alpha_i)$, and $m(\beta_i)$ 
in the condition $\six_t$ of $Q_t(\lambda/\mu)$ 
for $t=1$. 
}\label{fig:s-m}
\end{figure}

Let $t_0$ be the minimal number that satisfies $2^{t_0} > l$. 
For any $t=1, 2, \dots, t_0$, 
we define 
\begin{align*}
  & Q_t(\lambda/\mu):=
  \{ (\alpha; \beta) \in \cH(\lambda/\mu) \mid 
  \text{$(\alpha; \beta)$ satisfies Conditions 
    $\one_t$--$\six_t$} \}, \\
  & \hat{Q}_t(\lambda/\mu):=
\{ (\alpha; \beta) \in \cH(\lambda/\mu) \mid 
\text{$(\alpha; \beta)$ satisfies Conditions 
  $\one_t$--$\seven_t$} \}, 
\end{align*}
where Conditions $\one_t$--$\seven_t$
are given as follows: 
\begin{itemize}
\item[$\one_t$] $\alpha_i(0) \le \beta_i(0)$ for any $i=1, \dots, l$. 
\item[$\two_t$]
  $(\alpha; \beta)$ does not contain any odd $\II_1$-region. 
\item[$\three_t$]
  $(\alpha; \beta)$ does not contain any odd $\I_{2^t-1}$-region. 
\item[$\four_t$]
  $(\alpha; \beta)$ does not contain any $2^t$-overlap. 
\item[$\five_t$]
  If $t \ge 2$, then $(\alpha; \beta)$ contains at least one $2^{t-1}$-overlap. 
\item[$\six_t$]
  $s(\alpha_i) \equiv m_t(\alpha_i)$ and 
  $s(\beta_i) \equiv m_t(\beta_i)$, where 
  \begin{align*}
    s(\alpha_i):= & \alpha_i(0) -\beta_1(0) + i-1, \\
    s(\beta_i):= & \beta_i(0) -\beta_1(0) + i-1, \\
    m_t(\alpha_i):= & \# \{ j \mid 
    j \le i, \ 
    \text{$(\alpha_j, \beta_{j+2^t-1})$ is an even $(2^t-1)$-overlap} \}\\
    & + \# \{ j \mid 
    j < i, \ 
    \text{$(\alpha_{j}, \beta_{j+1})$ is an even $1$-hole} \}, \\
    m_t(\beta_i):= & \# \{ j \mid 
    j < i, \ 
    \text{$(\alpha_{j-2^t+1}, \beta_{j})$ is an even $(2^t-1)$-overlap} \} \\
    & + \# \{ j \mid 
    j \le i, \ 
    \text{$(\alpha_{j-1}, \beta_{j})$ is an even $1$-hole} \}. 
  \end{align*}
Here and the rest part of this section, 
$\equiv$ denotes the congruence modulo $2$. See Figure \ref{fig:s-m}. 
\item[$\seven_t$]
  $(\alpha; \beta)$ has at least one even $(2^t-1)$-overlap 
  (Then,  it has at least two, 
  because of Condition $\three_t$). 
\end{itemize}

\begin{rem}
If $t=t_0$, then Conditions 
$\three_t$ and $\four_t$ are void. Also, 
$\hat{Q}_{t_0}(\lambda/\mu)=\emptyset$, 
because $\seven_t$ is not satisfied for any 
$(\alpha; \beta) \in \cH(\lambda/\mu)$. 
\end{rem}

Now, the folding map 
$\phi: P_2(\lambda/\mu) \to P(\lambda/\mu)$ 
is constructed as follows: 
For any $\bp \in  P_2(\lambda/\mu)$ or any $\bp \in P(\lambda/\mu)$, 
one can associate  
$\pi(\bp):=(\alpha; \beta) \in \cH(\lambda/\mu)$ by 
removing all the E-steps from $\bp$. 
Then, we have 

\begin{prop}\label{prop:one}
The map $\pi$ 
gives a bijection  
$P_2(\lambda/\mu) \to Q_1(\lambda/\mu)$. 
\end{prop}

The following claim is the main part of the construction of $\phi$: 
\begin{prop}\label{prop:two}
There exists a bijection 
$\varphi_t: \hat{Q}_t(\lambda/\mu) \to Q_{t+1}(\lambda/\mu)$ 
for any $t=1, \dots, t_0 - 1$. 
\end{prop}
The map $\varphi_t$ is defined by using the $2^t$-folding. 
See \eqref{eq:t-folding}. 
Applying the bijections $\varphi_1$, $\varphi_2$, \dots repeatedly, 
we obtain a bijection  
$\varphi:Q_1(\lambda/\mu) \to R(\lambda/\mu)$, where 
\begin{equation}\label{eq:R}
R(\lambda/\mu):=
\bigsqcup_{t=1}^{t_0}\left(Q_t(\lambda/\mu)-\hat{Q}_t(\lambda/\mu)\right). 
\end{equation}
Furthermore, we have  
\begin{prop}\label{prop:three}
The map $\pi$ gives a bijection $P(\lambda/\mu)\to  R(\lambda/\mu)$. 
\end{prop}

Thanks to Propositions 
\ref{prop:one}--\ref{prop:three}, we now have a bijection $\phi$ 
by the commutative diagram, 
\begin{equation}\label{eq:folding-map}
\begin{CD}
P_2(\lambda/\mu) @>{\phi}>> P(\lambda/\mu)\\
@V{\pi}VV  @VV{\pi}V \\
Q_1(\lambda/\mu) @>>{\varphi}> R(\lambda/\mu). 
\end{CD}
\end{equation}

Then, Proposition \ref{prop:folding-map} follows from 
\begin{prop}\label{prop:four}
The map $\phi$ is weight-preserving. 
\end{prop}

\subsection{Proof of Proposition \ref{prop:one}}
We use the following lemma 
which immediately follows 
by the definition of $s(\alpha_i)$ and 
$s(\beta_i)$:  
\begin{lem}\label{lem:even}
For any $i=1, \dots, l$, 
\begin{enumerate}
\item \label{item:even-alpha-beta}
$\alpha_i(0)-\beta_j(0)$ is even 
if and only if $s(\alpha_i)-i\equiv s(\beta_j)-j$. 
\item \label{item:even-alpha}
$\alpha_i(0)-\alpha_j(0)$ is even 
if and only if $s(\alpha_i)-i\equiv s(\alpha_j)-j$. 
\item \label{item:even-beta}
$\beta_i(0)-\beta_j(0)$ is even 
if and only if $s(\beta_i)-i\equiv s(\beta_j)-j$. 
\end{enumerate}
\end{lem}

First, we prove that 
\begin{equation}\label{eq:inclusion}
\pi(P_2(\lambda/\mu)) \subset Q_1(\lambda/\mu). 
\end{equation}

Let $\bp \in P_2(\lambda/\mu)$
and set $(\alpha; \beta)=\pi(\bp)$. 
It holds for any $i=1, \dots, l$ that 
\begin{equation}\label{eq:half-paths-number}
\# \{ j \mid \alpha_i(0) \le \beta_j(0) \} \ge i.  
\end{equation}
Therefore, 
$\one_t$ (for $t=1$) is satisfied. 
By Condition \eqref{item:Ptwo-no-odd}
of $P_2(\lambda/\mu)$, 
$\two_t$ and $\three_t$ are satisfied, and  
by Condition \eqref{item:Ptwo-no-ord},  
$\four_t$ is satisfied. 

We prove Condition $\six_t$ for $\beta_i$ 
by induction with respect to $i$. 
The proof for $\alpha_i$ is similar. 
For $i=1$, 
we have $s(\beta_1)=m_1(\beta_1)=0$. 
Assume $s(\beta_{i-1}) \equiv m_1(\beta_{i-1})$. 
The following four cases should be considered: 
\begin{enumerate}[(a)]
\item \label{item:a}
$(\alpha_{i-1}, \beta_i)$ is an even 1-hole, and 
$(\alpha_{i-2}, \beta_{i-1})$ is an even 1-overlap.  
\item \label{item:b}
$(\alpha_{i-1}, \beta_i)$ is an even 1-hole, and 
$(\alpha_{i-2}, \beta_{i-1})$ is not an even 1-overlap.  
\item \label{item:c}
$(\alpha_{i-1}, \beta_i)$ is not an even 1-hole, and 
$(\alpha_{i-2}, \beta_{i-1})$ is an even 1-overlap.  
\item \label{item:d}
$(\alpha_{i-1}, \beta_i)$ is not an even 1-hole, and 
$(\alpha_{i-2}, \beta_{i-1})$ is not an even 1-overlap.  
\end{enumerate}
By the definition of $m_1(\beta_i)$, we have 
$m_1(\beta_i) = m_1(\beta_{i-1})+2$ for \eqref{item:a},
$m_1(\beta_{i-1})+1$ 
for \eqref{item:b} and \eqref{item:c}, and 
$m_1(\beta_{i-1})$ for \eqref{item:d}.  
On the other hand, 
by Condition \eqref{item:Ptwo-no-ord} of 
$P_2(\lambda/\mu)$, we have 
\begin{lem}\label{lem:alpha_i-beta_i}
$\alpha_{i-1}(0) - \beta_{i-1}(0)$ is odd
if and only if 
one and the only one of $(\alpha_{i-2}, \beta_{i-1})$ and  
$(\alpha_{i-1}, \beta_{i})$ is an even overlap. 
\end{lem}

By Lemma \ref{lem:alpha_i-beta_i}, 
$$\beta_i(0) -\beta_{i-1}(0)=
\left(
  \beta_i(0)-\alpha_{i-1}(0)
\right) + 
\left(
  \alpha_{i-1}(0) - \beta_{i-1}(0)
\right)
$$ 
is odd for \eqref{item:a} and \eqref{item:d}, 
and even  
for \eqref{item:b} and \eqref{item:c}.  
By the assumption of induction and 
Lemma \ref{lem:even} \eqref{item:even-beta}, 
we obtain $s(\beta_i) \equiv m_1(\beta_i)$ in each case. 

Next, we define the inverse map 
$\pi^{-1}:Q_1(\lambda/\mu) \to P_2(\lambda/\mu)$. 
For any $(\alpha;\beta) \in Q_1(\lambda/\mu)$, 
set $\bp:=\pi^{-1}(\alpha;\beta)$ as follows: 
Let $p_i:=[\alpha_i, \beta_j]$ denote the 
path defined by $\alpha_i$, $\beta_j$ and the 
consecutive E-steps from $\alpha_i(0)$ to $\beta_j(0)$. 
\begin{enumerate}[Step 1.]
\item \label{item:path-one}
First, for any $i$ such that  
$(\alpha_i, \beta_{i+1})$ is an even 1-overlap, 
set $p_i=[\alpha_i, \beta_{i+1}]$ 
(see $p_1$, $p_5$,  $p_6$, and  $p_7$ in Figure \ref{fig:s-m} for examples). 
Let $\Lambda_1$ (resp.~$\Lambda_2$)
be the set of all $i$ such that 
$\alpha_i$ (resp.~$\beta_i$) do not form even 1-overlaps. 
\item \label{item:path-two}
Next, for any $i \in \Lambda_1 \cap \Lambda_2$ 
such that $h:=\alpha_i(0)-\beta_i(0)$ 
is even (and non-positive by $\one_t$), set $p_i=[\alpha_i, \beta_{i}]$ 
(see  $p_3$ and  $p_4$ in Figure \ref{fig:s-m} for examples). 
Let $\Lambda'_1\subset \Lambda_1$ (resp.~$\Lambda'_2\subset \Lambda_2$)
be the set of all $i$ such that $h$ are not even. 
\item \label{item:path-three} 
Finally, for any $i\in \Lambda'_1 \cap \Lambda'_2$, 
there exists some $k>0$ such that 
$h':=\alpha_i(0)-\beta_{i-k}(0)<0$, 
by \eqref{eq:half-paths-number} and the previous steps. 
Let $k$ be the minimum of such numbers. 
Since $(\alpha_{i-k'}, \beta_{i-k'+1})$ for all $1 \le k' \le k$ 
is an even 1-overlap, we have 
$m_1(\alpha_i)=m_t(\beta_{i-k})+k$. 
Then, by Condition $\six_t$, 
$h'$ is even by Lemma \ref{lem:even} \eqref{item:even-alpha-beta}.  
Set $p_i=[\alpha_i, \beta_{i-k}]$ 
(see $p_2$ and $p_8$ in 
Figure \ref{fig:s-m} for examples). 
\end{enumerate}

Next, we prove that 
$\bp$ satisfies Condition 
\eqref{item:Ptwo-no-ord} of $P_2(\lambda/\mu)$. 
Namely, for any intersecting pair $(p_i, p_j)$ ($i<j$), 
we prove that $\alpha_i(0)-\alpha_j(0)$ is odd. 
The following two cases should be considered 
(see Figure \ref{fig:s-m}; 
the other cases do not occur by Condition $\four_t$):  
\begin{enumerate}[(A)]
\item {\it The case where $p_i$ is defined in 
Step \ref{item:path-one} and 
$p_j$ is defined in Step \ref{item:path-three}}. 
In this case, all $p_k$ ($i<k<j$) are defined in Step \ref{item:path-one}, 
by the definition of $p_j$. Therefore, $(\alpha_k, \beta_{k+1})$ 
for any $i<k<j$ is an even 1-overlap, and we have 
$m_1(\alpha_j)=m_1(\alpha_i)+j-i-1$. Thus, 
$\alpha_i(0)-\alpha_j(0)$ is odd by Condition $\six_t$ and 
Lemma \ref{lem:even} \eqref{item:even-alpha}. 
\item {\it The case where $(p_i, p_j)$ is a pair of paths defined in 
Steps \ref{item:path-two} and \ref{item:path-three}}. 
In this case, $p_j=[\alpha_j, \beta_{i+1}]$ and 
$(\alpha_i, \beta_{i+1})$ is a 1-overlap 
(otherwise, it contradicts to $\four_t$). 
Moreover, $(\alpha_i, \beta_{i+1})$ is odd 
(otherwise, $p_i$ should be $[\alpha_i, \beta_{i+1}]$). 
Since $\alpha_j(0)-\beta_{i+1}(0)$ is even by $p_j$, we obtain that 
$\alpha_i(0)-\alpha_j(0)$ is odd. 
\end{enumerate}

Finally, the facts that $\bp$ satisfies 
Condition \eqref{item:Ptwo-no-odd} of $P_2(\lambda/\mu)$ and 
$\pi^{-1}$ is the inverse of $\pi$ are 
obvious by construction. 

\subsection{Proof of Proposition \ref{prop:three}}
Set 
$$
m_{\infty}(\alpha_i)=
m_{\infty}(\beta_i) := 
\# \{ j \mid j \le i  \text{ and $(\alpha_{j-1}, \beta_j)$ 
is an even $1$-hole} \}.
$$
Then 
\begin{lem}\label{lem:condition}
Let $t=1, \dots, t_0$. 
For any $(\alpha; \beta) \in \cH(\lambda/\mu)$,  
$(\alpha; \beta) \in 
Q_t(\lambda/\mu) \backslash \hat{Q}_t(\lambda/\mu)$ 
if and only if $(\alpha; \beta)$ satisfies conditions 
$\one_t$, $\two_t$, $\four_t$, 
$\five_t$ and the following condition: 

\noindent
\begin{tabular}{ll}
$\six'_t$ & $s(\alpha_i) \equiv s(\beta_i) \equiv 
m_{\infty}(\alpha_i) \equiv m_{\infty}(\beta_i)$. 
\end{tabular}
\end{lem}

\begin{proof}
($\Rightarrow$) If $(\alpha; \beta) \in
Q_t(\lambda/\mu) \backslash \hat{Q}_t(\lambda/\mu)$, 
then 
$$m_t(\alpha_i)=m_t(\beta_i)=
m_{\infty}(\alpha_i)=m_{\infty}(\beta_i),$$ 
and we obtain $\six'_t$. 

($\Leftarrow$) 
We show that  
$(\alpha; \beta)$ does not have any even $(2^t-1)$-overlap. 
If a pair $(\alpha_i, \beta_{i+2^t-1})$ is 
an even $(2^t-1)$-overlap, then  
$s(\alpha_i) \not\equiv s(\beta_{i+2^t-1})$ 
by Lemma \ref{lem:even} \eqref{item:even-alpha-beta}. 
On the other hand, there does not exist any $1$-hole 
between $\alpha_i$ and $\beta_{i+2^t-1}$, 
which implies 
$m_{\infty}(\alpha_i)=m_{\infty}(\beta_{i+2^t-1})$, 
and it contradicts to $\six'_t$. 
Therefore, $\seven_t$ is satisfied, and moreover, we also have 
$m_{\infty}(\alpha_i)=m_t(\alpha_i)$ and 
$m_{\infty}(\beta_i)=m_t(\beta_i)$. 
Thus, $\six_t$ is satisfied. 
Finally, $\three_t$ is satisfied  
because every $\I_{2^t-1}$-region $V$ of $(\alpha; \beta)$ 
satisfies $n(V)=0$ 
by the definition of $n(V)$.  
\end{proof}

By Lemma \ref{lem:condition}, 
the set 
$R(\lambda/\mu)$ in \eqref{eq:R} is described as follows: 
\begin{equation*}
R(\lambda/\mu)=
\{ (\alpha; \beta) \in \cH(\lambda/\mu) \mid  
\text{
$(\alpha; \beta)$ satisfies 
$\one_t$, $\two_t$ and $\six'_t$}
\}. 
\end{equation*}

Next, we prove  
$\pi(P(\lambda/\mu)) \subset R(\lambda/\mu)$. 
Conditions $\one_t$ and $\two_t$ are satisfied, 
similarly as in the proof of \eqref{eq:inclusion}. 
We prove $\six'_t$.  
Since $p_i=[\alpha_i, \beta_i]$, 
$\alpha_i(0)-\beta_i(0)$ is even. Thus, 
$s(\alpha_i)\equiv s(\beta_i)$ by 
Lemma \ref{lem:even} \eqref{item:even-alpha-beta}. 
We prove $s(\alpha_i) \equiv m_{\infty}(\alpha_i)$ by induction 
with respect to $i$. 
The proof for $s(\alpha_i)\equiv m_{\infty}(\beta_i)$ is similar.  
For $i=1$, we have 
$s(\beta_1)=m_{\infty}(\alpha_1)=0$, and we obtain 
$s(\alpha_1) \equiv s(\beta_1) \equiv m_{\infty}(\alpha_1)$. 
Assume 
$s(\alpha_{i-1}) \equiv m_{\infty}(\alpha_{i-1})$.  
Then the following two cases should be considered: 
\begin{enumerate}[(a)]
\item \label{item:two-one} 
$(\alpha_{i-1}, \beta_i)$ is an even 1-hole. 
\item \label{item:two-two} 
$(\alpha_{i-1}, \beta_i)$ is not an even 1-hole. 
\end{enumerate}
By the definition, 
$m_{\infty}(\alpha_i)=m_{\infty}(\alpha_{i-1})+1$ 
for \eqref{item:two-one} and 
$m_{\infty}(\alpha_{i-1})$ 
for \eqref{item:two-two}. 
On the other hand, 
we have 
$s(\beta_i) \not\equiv s(\alpha_{i-1})$ 
for \eqref{item:two-one}, and 
$s(\beta_i) \equiv s(\alpha_{i-1})$ 
for \eqref{item:two-two}.  
Therefore, we obtain 
$s(\alpha_i) \equiv s(\beta_i) \equiv m_{\infty}(\alpha_i)$ 
in each case.  

Finally, we define the inverse map 
$\pi^{-1}:R(\lambda/\mu) \to P(\lambda/\mu)$ by 
$\bp:=\pi^{-1}(\alpha; \beta)$, 
$p_i:=[\alpha_i, \beta_i]$ ($i=1, \dots, l$)
for any $(\alpha; \beta) \in R(\lambda/\mu)$.  
The path $p_i$ is well defined, because of $\six'_t$ and 
Lemma \ref{lem:even} \eqref{item:even-alpha-beta}.  
If $(p_i,p_{i+1})$ is ordinarily intersecting, 
then we have $s(\alpha_i) \not\equiv s(\alpha_{i+1})$ 
by Lemma \ref{lem:even} \eqref{item:even-alpha}
and $m_{\infty}(\alpha_i)=m_{\infty}(\alpha_{i+1})$, 
which contradicts to $\six'_t$. 
Therefore, 
$\bp$ satisfies Condition \eqref{item:p-no-ord} of $P(\lambda/\mu)$. 
The facts that $\bp$ satisfies Condition \eqref{item:p-no-odd} of $P(\lambda/\mu)$  
and $\pi^{-1}$ is the inverse of $\pi$ are  
obvious. 

\subsection{Definition of $\varphi_t$}
Here, we give the definition of the map $\varphi_t$ 
in Proposition \ref{prop:two}. 
Fix $k=2, \dots, l$. 
Suppose that $(\alpha; \beta) \in \cH(\lambda/\mu)$ 
has an even number $2m$ of even $(k-1)$-overlaps 
$(\alpha_{i_1}, \beta_{i_1+k-1})$, \dots, 
$(\alpha_{i_{2m}}, \beta_{i_{2m}+k-1})$
($i_1< \dots < i_{2m}$). 
We say $(\alpha_{i_j}, \beta_{i_j+k-1})$ 
($1 \le j \le 2m$) is of {\it R-type} if $j$ is odd, and of  
{\it L-type} if $j$ is even. 
Let $(\alpha_{i}, \beta_{i+k-1})$ and 
$(\alpha_{i'}, \beta_{i'+k-1})$ be a nearest pair of 
even $(k-1)$-overlaps with $i=i_j$, $i'=i_{j+1}$. 
Then, we say the height $0$ units $U \subset S_+$ between 
$\beta_{i'+k-1}$ and $\alpha^*_{i}$ and 
their duals $U^*\subset S_-$ 
are of {\it LR-type} if 
$j$ is odd, and of {\it RL-type} otherwise. 
Remark that any height $0$ $\II_k$-unit 
is either of LR-type or RL-type. 

The next lemma is the key for  
the definition and the bijectivity of $\varphi_t$. 
\begin{lem}\label{lem:complementary}
$(\alpha; \beta) \in \cH(\lambda/\mu)$ 
does not have any odd $\I_{k-1}$-region 
if and only if 
the following conditions are satisfied: 
\begin{enumerate}[(i)]
\item \label{item:complementary-one}
$(\alpha; \beta)$ 
has an even number of even $(k-1)$-overlaps. 
\item \label{item:2t-one}
No connected component of $\tilde{\cU}_{\II_k}$ 
of $(\alpha; \beta)$
contains height $0$ $\II_k$-units of both LR- and RL-type, 
simultaneously. 
\item \label{item:c-two}
Any connected component of $\tilde{\cU}_{\II_k}$ 
of $(\alpha; \beta)$ 
which contains a 
height $0$ $\II_k$-unit of LR-type 
is a $\II_k$-region. 
\end{enumerate}
\end{lem}

The proofs of the lemma requires 
some graph-theoretical consideration, 
and it is given in Appendix \ref{sec:graph}. 

\begin{figure}
\begin{psfrags}
\psfrag{L}{\tiny $L$}
\psfrag{R}{\tiny $R$}
\psfrag{V1}{\small $V_1$}
\psfrag{V2}{\small $V_2$}
\psfrag{V3}{\small $V_3$}
\psfrag{V4}{\small $V_4$}
\psfrag{V5}{\small $V_5$}
\psfrag{a1}{\small $\alpha_1$}
\psfrag{a2}{\small $\alpha_2$}
\psfrag{a3}{\small $\alpha_3$}
\psfrag{a4}{\small $\alpha_4$}
\psfrag{a5}{\small $\alpha_5$}
\psfrag{a6}{\small $\alpha_6$}
\psfrag{a7}{\small $\alpha_7$}
\psfrag{a8}{\small $\alpha_8$}
\psfrag{b1}{\small $\beta_1$}
\psfrag{b2}{\small $\beta_2$}
\psfrag{b3}{\small $\beta_3$}
\psfrag{b4}{\small $\beta_4$}
\psfrag{b5}{\small $\beta_5$}
\psfrag{b6}{\small $\beta_6$}
\psfrag{b7}{\small $\beta_7$}
\psfrag{b8}{\small $\beta_8$}
\psfrag{pi}{\small $\varphi_t=\varepsilon^{2^t}_{V_2} \circ \varepsilon^{2^t}_{V_5}$}
\includegraphics[width=8.5cm, clip]{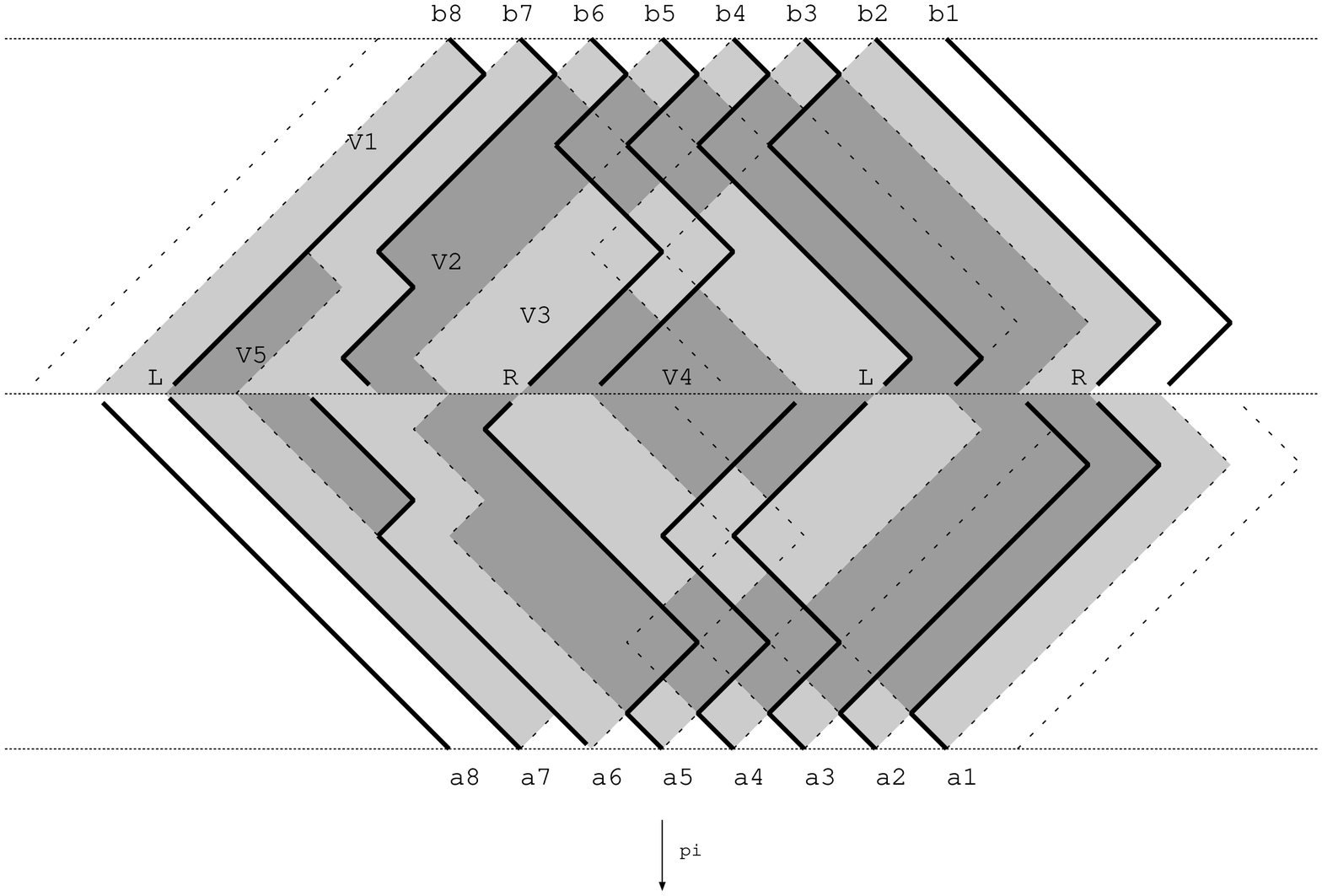}
\end{psfrags}
\vspace{5pt}\\
\begin{psfrags}
\psfrag{V1}{\small $V_1$}
\psfrag{V2}{\small $V_2$}
\psfrag{V3}{\small $V_3$}
\psfrag{V4}{\small $V_4$}
\psfrag{V5}{\small $V_5$}
\psfrag{a1}{\small $\alpha'_1$}
\psfrag{a2}{\small $\alpha'_2$}
\psfrag{a3}{\small $\alpha'_3$}
\psfrag{a4}{\small $\alpha'_4$}
\psfrag{a5}{\small $\alpha'_5$}
\psfrag{a6}{\small $\alpha'_6$}
\psfrag{a7}{\small $\alpha'_7$}
\psfrag{a8}{\small $\alpha'_8$}
\psfrag{b1}{\small $\beta'_1$}
\psfrag{b2}{\small $\beta'_2$}
\psfrag{b3}{\small $\beta'_3$}
\psfrag{b4}{\small $\beta'_4$}
\psfrag{b5}{\small $\beta'_5$}
\psfrag{b6}{\small $\beta'_6$}
\psfrag{b7}{\small $\beta'_7$}
\psfrag{b8}{\small $\beta'_8$}
\includegraphics[width=8.5cm, clip]{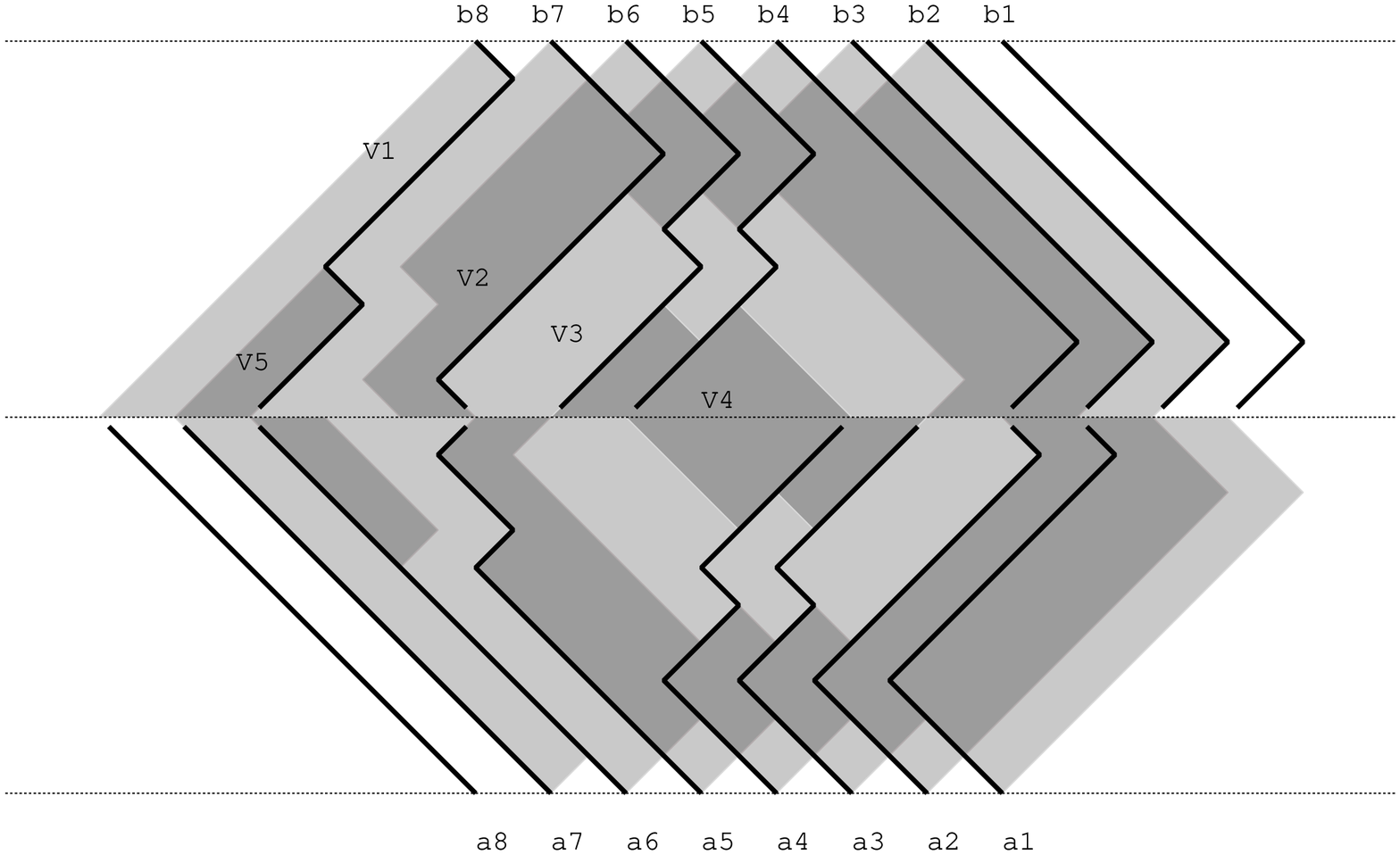}
\end{psfrags}
\caption{An example of $(\alpha; \beta) \in \hat{Q}_t(\lambda/\mu)$, 
its $\I_{2^t-1}$-regions $V_1$, $V_3$ 
and $\II_{2^t}$-regions $V_2$, $V_4$, $V_5$, 
and the map $\varphi_t$ for $t=1$. 
The even $(2^t-1)$-overlaps $(\alpha_1, \beta_2)$ and $(\alpha_5, \beta_6)$ 
are of R-type, $(\alpha_3, \beta_4)$ and $(\alpha_7, \beta_8)$
are of L-type, and   
$V_2$ and $V_5$ are $\II_{2^t}$-regions of LR-type. 
}\label{fig:2t-folding}
\end{figure}

Now, let $(\alpha; \beta)\in \hat{Q}_t(\lambda/\mu)$ and 
let $V$ be a $\II_{2^t}$-region of $(\alpha; \beta)$. 
Since $(\alpha; \beta)$ does not have any 
odd $\I_{2^t-1}$-region by $\three_t$, 
one of the following is satisfied  
by Lemma \ref{lem:complementary}
with $k=2^t$: 
\begin{enumerate}[(a)]
\item \label{item:LR-region} 
All the height $0$ units in $V$ are of LR-type. 
\item 
All the height $0$ units in $V$ are of RL-type. 
\end{enumerate}
If \eqref{item:LR-region} is satisfied, 
we say $V$ is of {\it LR-type}. 
See Figure \ref{fig:2t-folding} for an example. 
We define the map $\varphi_t$ in Proposition \ref{prop:two} by 
{\it (the composition of) the $2^t$-foldings 
with respect to all the $\II_{2^t}$-regions 
$V_1$, \dots, $V_p$ of LR-type of $(\alpha; \beta)$}, i.e., 
\begin{equation}\label{eq:t-folding}
\varphi_t: (\alpha; \beta) \mapsto 
\varepsilon_{V_1}^{2^t} \circ 
\dots \circ \varepsilon_{V_p}^{2^t}(\alpha; \beta). 
\end{equation}

\subsection{Proof of Proposition \ref{prop:two}}
To begin with, we give two lemmas. 

\begin{lem}\label{lem:odd-even}
Let $(\alpha;\beta)\in\hat{Q}_t(\lambda/\mu)$ and 
$W=V_1 \cup \dots \cup V_p$ for $V_i$'s in \eqref{eq:t-folding}. 
Then,  
\begin{enumerate}
\item \label{item:alpha}
$\alpha_i$ intersects with $W$ at height $0$ if and only if 
$$\# \{ j \mid 
j \le i 
\text{ and $(\alpha_j, \beta_{j+2^t-1})$ is an even $(2^t-1)$-overlap} \} $$
is odd. 
\item \label{item:beta}
$\beta_i$ intersects with $W$ at height $0$ if and only if
$$\# \{ j \mid
j < i
\text{ and $(\alpha_{j-2^t+1}, \beta_{j})$ is an even $(2^t-1)$-overlap} \} $$
is odd.
\end{enumerate}
\end{lem}

\begin{proof}
We prove it for  \eqref{item:alpha}. 
Let $U\subset W$ be a height $0$ unit
which intersects with $\alpha_i$. 
Since $U$ is of LR-type, and 
the number of the even $(2^t-1)$-overlaps 
$(\alpha_{j}, \beta_{j+2^t-1})$ to the right of $U$ 
is odd. 

Conversely, suppose the above number is odd. 
Then, 
the unit $U$ whose right vertex intersects with 
$\alpha_i$ is of LR-type. 
Moreover, $U$ is a $\II_{2^t}$-unit because 
$\beta_{i+2^t}(0) \le \alpha_i(0)-2$, 
by Conditions $\four_t$ and $\six_t$. 
Then, $U$ is included in $W$ by Lemma \ref{lem:complementary}; 
therefore, $\alpha_i$ intersects  
with $W$. 
\end{proof}

\begin{lem}\label{lem:folding-result}
Let 
$V$ be a $\II_k$-region of $(\alpha; \beta)\in \cH(\lambda/\mu)$ and 
$(\alpha'; \beta')=\varepsilon_V^k(\alpha; \beta)$. Then,  
\begin{enumerate}
\item \label{item:folding-result-one}
For any  $r=1, \dots, k$, 
both $\alpha_i$ and $\beta_{i+r}$ 
intersect with $V$ at height $0$ and 
$\alpha_i(0)\ge \beta_{i+r}(0) + 2$ if and only if 
both $\alpha'_{i+r-k}$ and $\beta'_{i+k}$ 
intersect with $V$ at height $0$ and 
$\alpha'_{i+r-k}(0) \le \beta'_{i+k}(0)$. 
\item \label{item:folding-result-two}
Suppose that $(\alpha; \beta)$ does not have any odd $\I_{2k-1}$-region. 
Then, 
$V'$ is an odd $\II_1$-region of $(\alpha; \beta)$ 
if and only if 
$V'$ is an odd $\II_1$- or odd $\I_{2k-1}$-region of $(\alpha'; \beta')$.  
\end{enumerate}
\end{lem}

\begin{proof}
\eqref{item:folding-result-one} It is obtained by 
\eqref{eq:k-folding}. 
\eqref{item:folding-result-two} 
If $V'$ is an odd $\II_1$-region of $(\alpha; \beta)$, 
then by Lemma \ref{lem:inclusion} \eqref{item:inclusion-II}, 
$V'\subset V$ or $V'\cap V=\emptyset$. 
In the former case, 
$V'$ is a $\I_{2k-1}$-region of $(\alpha'; \beta')$ 
by Proposition \ref{prop:region} \eqref{item:folding-result-zero}. 
In the latter case, $V'$ is also a $\II_1$-region of 
$(\alpha'; \beta')$. The converse is similar. 
\end{proof}

Now, let us prove $\varphi_t: \hat{Q}_t(\lambda/\mu) \to Q_{t+1}(\lambda/\mu)$ 
is a bijection. 
First, we prove 
$\varphi_t(\hat{Q}_t(\lambda/\mu)) \subset Q_{t+1}(\lambda/\mu)$. 
Set $(\alpha'; \beta')=\varphi_t(\alpha; \beta)$. 
Let $V_1$, \dots, $V_p$ be the set of all the 
$\II_{2^t}$-regions of LR-type of $(\alpha; \beta)$, and 
set $W=V_1 \cup \dots \cup V_p$. 
We remark that $V_i$ is a $\I_{2^t}$-region of $(\alpha';\beta')$ 
by Proposition \ref{prop:region} \eqref{item:folding-result-zero}. 

\noindent \, $\bullet$ $\one_{t+1}$. 
By the definition of the $k$-folding and 
Condition $\one_t$ of $(\alpha; \beta)$, 
we have $\alpha'_i(0) \le \alpha_i(0) \le \beta_i(0) \le \beta'_i(0)$. 

\noindent \, $\bullet$ {\it $\two_{t+1}$ and $\three_{t+1}$.} 
This is 
obtained by Conditions $\two_{t}$ and $\four_{t}$ of $(\alpha; \beta)$
and Lemma \ref{lem:folding-result} \eqref{item:folding-result-two}. 

\noindent \, $\bullet$ $\four_{t+1}$.
Suppose that  $(\alpha'_i, \beta'_{i+2^{t+1}})$ is a $2^{t+1}$-overlap 
of $(\alpha'; \beta')$. 
Since $(\alpha_i, \beta_{i+2^{t+1}})$ is not 
a $2^{t+1}$-overlap by $\four_t$, 
$\alpha'_i$ and $\beta'_{i+2^{t+1}}$ 
should intersect with $W$ at height $0$ (otherwise,  
$(\alpha_i, \beta_{i+2^{t+1}})=(\alpha'_i, \beta'_{i+2^{t+1}})$). 
As \eqref{eq:k-folding}, 
we have $\alpha'_i(0)=\beta_{i+2^t}^*(0)$ and 
$\beta'_{i+2^{t+1}}(0)=\alpha^*_{i+2^t}(0)$, 
and therefore, $\beta_{i+2^t}(0) \le \alpha_{i+2^t}(0)-2$, 
which contradicts to Condition $\one_t$ of $(\alpha; \beta)$. 

\noindent \, $\bullet$ $\five_{t+1}$. 
Since $(\alpha; \beta)$ satisfies $\seven_{t}$, there exists 
an even $(2^t-1)$-overlap 
$(\alpha_i, \beta_{i+2^t-1})$.  
By Condition $\four_{t}$, 
$(\alpha_{i-1}, \beta_{i+2^t-1})$ is not a $2^t$-overlap. 
If $\alpha_{i-1}(0) = \beta_{i+2^t-1}(0)+1$, 
then $m_t(\alpha_{i-1}) = m_t(\beta_{i+2^t-1})$, 
therefore, $\alpha_{i-1}(0) \equiv \beta_{i+2^t-1}(0)$ 
by Condition $\six_{t}$, which is a contradiction. 
Thus, $\alpha_{i-1}(0) \ge  \beta_{i+2^t-1}(0)+2$, 
and $\alpha_{i-1}$ and $\beta_{i+2^t-1}$ intersect with 
$W$ at height $0$. 
Therefore, by Lemma \ref{lem:folding-result} \eqref{item:folding-result-one}, 
we have $\alpha'_{i-1}(0) \le \beta'_{i+2^t-1}(0)$, 
namely, $(\alpha'_{i-1}, \beta'_{i+2^t-1})$ is a 
$2^t$-overlap. 

\noindent \, $\bullet$ $\six_{t+1}$. 
We prove $s(\alpha'_i) = m_t(\alpha'_i)$ 
(the proof for $s(\beta'_i) = m_t(\beta'_i)$
is similar). 

{\it Case 1}. If $\alpha_i$ does not intersect with 
$W$ at height $0$, which implies $\alpha'_i(0) = \alpha_i(0)$, 
then $s(\alpha'_i) = s(\alpha_i)$. 
On the other hand, we have 
\begin{equation}\label{eq:mt}
m_{t+1}(\alpha'_i)-m_t(\alpha_i) = 
- \# \{ j \mid
    j \le i, \
    \text{$(\alpha_j, \beta_{j+2^t-1})$ is an even $(2^t-1)$-overlap} \}, 
\end{equation}
by Lemma \ref{lem:folding-result} \eqref{item:folding-result-one} 
and Condition $\four_t$ of $(\alpha; \beta)$. Moreover, 
the right-hand side of \eqref{eq:mt} is even 
by Lemma \ref{lem:odd-even} \eqref{item:alpha}. 
By Condition $\six_{t}$ of $(\alpha; \beta)$, 
we obtain $m_{t+1}(\alpha'_i) \equiv m_t(\alpha_i)
\equiv s(\alpha_i)=s(\alpha'_i)$. 

{\it Case 2}. If $\alpha_i$ intersects with $W$ 
at height $0$, which implies that $\alpha'_i(0)=\beta^*_{i+2^t}(0)$, 
then $s(\alpha'_i) = s(\beta_{i+2^t})+1$. 
On the other hand, 
we have 
\begin{equation}\label{eq:mt-2}
m_{t+1}(\alpha'_i)-m_t(\beta_{i+2^t}) =
- \# \{ j \mid
    j < i, \
    \text{$(\alpha_{j-2^t+1}, \beta_{j})$ is an even $(2^t-1)$-overlap} \}, 
\end{equation}
by Lemma \ref{lem:folding-result} \eqref{item:folding-result-one}
and Condition $\four_t$ of $(\alpha; \beta)$. Moreover, 
the right-hand side of \eqref{eq:mt} is odd 
by Lemma \ref{lem:odd-even} \eqref{item:alpha}. 
By Condition $\six_{t}$ of $(\alpha; \beta)$,
we obtain $m_{t+1}(\alpha'_i) \not\equiv m_t(\beta_{i+2^t})
\equiv s(\beta_{i+2^t})\not\equiv s(\alpha'_i)$.

Next, let us define the inverse map 
$$\varphi_t^{-1}: Q_{t+1}(\lambda/\mu) \to \hat{Q}_t(\lambda/\mu)$$
as follows: 
For any $(\alpha'; \beta') \in Q_{t+1}(\lambda/\mu)$, set 
$$\varphi^{-1}_t: (\alpha'; \beta') \mapsto (\alpha; \beta):=
\varepsilon_{V_1}^{2^t} \circ \dots \circ \varepsilon_{V_p}^{2^t}(\alpha'; \beta'), $$ 
where $V_1, \dots, V_p$ are the set of all the 
$\I_{2^t}$-regions of $(\alpha'; \beta')$. 
By $\five_{t+1}$, we have $p \ge 1$. 
Set $W=V_1\cup \dots \cup V_p$. 
We prove that $(\alpha; \beta) \in \hat{Q}_t(\lambda/\mu)$. 

\noindent \, $\bullet$ $\one_{t}$. 
By \eqref{eq:k-folding} wherein 
$(\alpha_i)$, $(\alpha'_i)$ and $(\beta_i)$, $(\beta'_i)$ are interchanged, 
$\alpha_i(0)={\beta'}^*_{i+2^t}(0)$ or $\alpha'_i(0)$, and 
$\beta_i(0)={\alpha'}^*_{i-2^t}(0)$ or $\beta'_i(0)$. 
Only the case 
$(\alpha_i(0), \beta_i(0)) = ({\beta'}^*_{i+2^t}(0), {\alpha'}^*_{i-2^t}(0))$ 
is nontrivial. 
By $\four_{t+1}$, we have $\alpha'_{i-2^t}(0) > \beta'_{i+2^t}(0)$. 
Then, we have $\alpha'_{i-2^t}(0) \le  \beta'_{i+2^t}(0)$ 
by the same argument in the proof of 
$\five_{t+1}$ of $\varphi_t(\alpha; \beta)$.

\noindent \, $\bullet$ $\two_{t}$. 
This is because of Lemma 
\ref{lem:folding-result} \eqref{item:folding-result-two} and 
the Conditions $\two_{t+1}$ and $\three_{t+1}$ of $(\alpha';\beta')$. 

\noindent \, $\bullet$ $\four_{t}$. 
Suppose that $(\alpha; \beta)$ has a $2^t$-overlap 
$(\alpha_i, \beta_{i+2^t})$. 
By Proposition \ref{prop:region} \eqref{item:expansion-result}, 
$V_1, \dots, V_p$ are $\II_{2^t}$-regions of $(\alpha; \beta)$, 
and therefore, 
$\alpha_i$ and $\beta_{i+2^t}$ do not 
intersect with $W$ at height $0$ by Lemma \ref{lem:unit} 
\eqref{item:not-connected}. 
Then,  
$(\alpha'_i, \beta'_{i+2^t})$ is also a $2^t$-overlap 
because 
$\alpha'_i(0)=\alpha_i(0)$ and  $\beta'_{i+2^t}(0)=\beta_{i+2^t}(0)$, which 
implies that 
there exists a height $0$ $\I_{2^t}$-unit $U\not\subset W$ of $(\alpha'; \beta)$ 
between $\alpha'_i$ and $\beta'_{i+2^t}$,  
and then  contradicts to the definition of $W$. 

\noindent \, $\bullet$ {\it $\three_{t}$ and $\seven_t$. } 
This is the most non-trivial part of the proof of Proposition \ref{prop:two}. 
We prove that  
Conditions 
\eqref{item:complementary-one}--\eqref{item:c-two}  
in Lemma \ref{lem:complementary} 
(for $k=2^t$) are satisfied. 

Let us study when a pair $(\alpha_i, \beta_{i+2^t-1})$ is 
an even $(2^t-1)$-overlap. 
There are four cases \ref{item:A}--\ref{item:D} 
to be considered. 
We prove that $(\alpha_i, \beta_{i+2^t-1})$ is even only in 
Cases \ref{item:A} and \ref{item:B}. 

\begin{enumerate}[A.]
\item \label{item:A}
{\it The case where $\alpha_i$ does not intersect with $W$ 
at height $0$ and 
$\beta_{i+2^t-1}$ intersects with $W$ at height $0$}. 
In this case, 
we have $\beta^*_{i+2^t-1}(0)=\alpha'_{i-1}(0)$ 
and $\alpha_i(0)=\alpha'_i(0)$, and then
$(\alpha_i, \beta_{i+2^t-1})$ is a $(2^t-1)$-overlap. 
We also have $m_{t+1}(\alpha'_i) =m_{t+1}(\alpha'_{i-1})$. 
This is because $(\alpha'_{i-1}, \beta'_i)$ is not an even $1$-hole,   
and $(\alpha'_{i}, \beta'_{i+2^{t+1}-1})$ is not an (even) 
$(2^{t+1}-1)$-overlap (otherwise, 
$(\alpha'_{i}, \beta'_{i+2^{t}})$ is also a $2^t$-overlap, 
and $\alpha_i$ intersects with $W$ at height $0$, which 
contradicts to the assumption). 
Therefore, $s(\alpha'_i)\equiv s(\alpha'_{i-1})$, i.e., 
$\alpha'_i(0) \not\equiv \alpha'_{i-1}(0)$. 
So, $(\alpha_{i}, \beta_{i+2^{t}-1})$ is an even $(2^t-1)$-overlap.
\item \label{item:B}
{\it The case where $\alpha_i$ intersects with $W$ at height $0$ 
and $\beta_{i+2^t-1}$ does not intersect with $W$ at height $0$}. 
As in Case \ref{item:A}, one can show that 
$(\alpha_{i}, \beta_{i+2^{t}-1})$ 
is an even $(2^{t}-1)$-overlap. 
\item \label{item:C}
{\it The case where both $\alpha_i$ and
$\beta_{i+2^t-1}$ intersect with $W$ at height $0$}. 
In this case, we have 
$\beta^*_{i+2^t-1}(0)=\alpha'_{i-1}(0)$, $\alpha^*_i(0)=\beta'_{i+2^t}(0)$, 
and $(\alpha'_{i-1}, \beta'_{i+2^t-1})$, 
$(\alpha'_{i}, \beta'_{i+2^t})$ are $2^t$-overlaps. 
If $(\alpha_i, \beta_{i+2^t-1})$ is a $(2^t-1)$-overlap, 
then there does not exist any path at height 0 
between $\beta'_{i+2^t}(0)$ and $\alpha'_{i-1}(0)$.  
We have $m_{t+1}(\alpha'_{i-1})=m_{t+1}(\beta'_{i+2^t})$, 
therefore, $s(\alpha'_{i-1})\equiv s(\beta'_{i+2^t})$, 
namely, $(\alpha_i, \beta_{i+2^t-1})$ is odd. 
\item \label{item:D}
{\it The case where both $\alpha_i$ and 
$\beta_{i+2^t-1}$ do not intersect with $W$ at height $0$}.
In this case, we have $\alpha_i(0)=\alpha'_i(0)$ and 
$\beta'_{i+2^t-1}(0)=\beta_{i+2^t-1}(0)$. 
If $(\alpha_i, \beta_{i+2^t-1})$ is 
a $(2^t-1)$-overlap, then 
there does not exist any path 
at height 0 between $\alpha'_i$ and $\beta'_{i+2^t-1}$ 
because of the assumption. 
Then we have $m_{t+1}(\beta'_{i+2^t-1})\equiv m_{t+1}(\alpha'_i)$, 
therefore, $(\alpha_i, \beta_{i+2^t-1})$ is odd. 
\end{enumerate}

To summarize, the even $(2^t-1)$-overlaps 
of L-type (resp.~R-type) of $(\alpha; \beta)$  
are the ones in Case \ref{item:A} 
(resp.~Case \ref{item:B}); 
furthermore, 
the number of the even $(2^t-1)$-overlaps 
of $(\alpha; \beta)$ is even, and 
a  height 0 $\II_{2^t}$-unit $U$  
of $(\alpha; \beta)$ is of LR-type 
if and only if $U \subset W$. 
Also, $W$ is a union of $\II_{2^t}$-regions of $(\alpha; \beta)$ 
by Proposition \ref{prop:region} \eqref{item:expansion-result}. 
Thus, Conditions 
\eqref{item:complementary-one}--\eqref{item:c-two} in 
Lemma \ref{lem:complementary}
(for $k=2^t$) are satisfied.

\noindent \, $\bullet$ $\five_{t}$. 
If $(\alpha_i, \beta_{i+2^t-1})$ is an overlap, then 
$(\alpha_i, \beta_{i+2^{t-1}})$ is also an overlap, 
and therefore, $\five_{t}$ holds. 

\noindent \, $\bullet$ $\six_{t}$. 
We prove $s(\alpha_i)\equiv m_{t}(\alpha_i)$. 
\begin{enumerate}[1.]
\item {\it The case where $\alpha_i$ does not intersect with 
$W$ at height $0$}. 
This is the cases \ref{item:A} and \ref{item:D} in the 
proof of $\three_t$. By $\alpha_i(0) = \alpha'_i(0)$, 
we have $s(\alpha_i)=s(\alpha'_i)$. 
Since every even $(2^{t+1}-1)$-overlap of $(\alpha'; \beta')$ 
intersects with $W$ at height $0$, we have 
$(\alpha_i, \beta_{i+1})$ is an even $1$-hole  
by Lemma \ref{lem:folding-result} \eqref{item:folding-result-one}.  
Therefore, 
\begin{equation}\label{eq:AD}
m_t(\alpha_i)-m_{t+1}(\alpha'_i)= 
\#\{j \mid j \le i, \text{$(\alpha_j, \beta_{j+2^t-1})$ is an even 
$(2^t-1)$-overlap} \}. 
\end{equation}
{}From 
the proof of $\three_{t}$, the right-hand side of 
\eqref{eq:AD} is even. Thus, $m_t(\alpha_i)\equiv s(\alpha_i)$. 
\item {\it The case where $\alpha_i$ intersects with $W$ at height $0$}.
This is the cases \ref{item:B} and \ref{item:C} in the proof of 
$\three_t$. 
By $\alpha^*_i(0)=\beta'_{i+2^t}(0)$ and $\beta'_1(0)=\beta_1(0)$, 
we have $s(\alpha_i) \not\equiv s(\beta'_{i+2^t})$. 
As in the former case, we have 
\begin{equation}\label{eq:CD}
m_t(\alpha_i)-m_{t+1}(\beta'_{i+2^t})=
\#\{j \mid j \le i, \text{$(\alpha_j, \beta_{j+2^t-1})$ is an even
$(2^t-1)$-overlap} \}, 
\end{equation}
and the right-hand side of
\eqref{eq:CD} is odd. Thus, $m_t(\alpha_i)\not\equiv m_{t+1}(\beta'_{i+2^t})$, 
which implies $m_t(\alpha_i)\equiv s(\alpha_i)$.
\end{enumerate}

Finally, the fact that $\pi^{-1}$ is the inverse of $\pi$ 
is obvious by construction.       

\subsection{Proof of Proposition \ref{prop:four}}
Finally, we show that $\phi$ is weight-preserving 
and then complete the proof of Proposition \ref{prop:folding-map}. 

For $\bp\in P_2(\lambda/\mu)$ and
$\bp'=\phi(\bp)\in P(\lambda/\mu)$,
we shall show that 
$z_a^{\bp}=z_a^{\bp'}$.

Let $(\alpha;\beta)\in Q_1({\lambda/\mu})$ and
$(\alpha';\beta')\in R({\lambda/\mu})$ be
the corresponding ones to $\bp$ and $\bp'$
under the identification $\pi$.
We decompose the monomial $z_a^{\bp}$ in \eqref{eq:weight} into
two parts as $z_a^{\bp}=HE$,
where $H$ is the factor coming from the
lower and upper paths $(\alpha;\beta)$ for $\bp$,
and $E$ is the one from the height 0 parts (E-steps) of
$\bp$.
We do the same for $z_a^{\bp'}$ as $z_a^{\bp'}=H'E'$.

By a similar argument for the proof of Proposition \ref{prop:second-involution},
one obtains that
\begin{align}\label{eq:H-weight}
H'=H\delta^{-1},
\end{align}
where
\begin{equation}\label{eq:delta-weight}
\begin{aligned}
\delta&=
\prod_{i=1}^l
\left(
\prod_{k=\alpha'_i(0)}^{\alpha_i(0)-1}z_{n,a-2k}
\prod_{k=\beta_i(0)}^{\beta'_i(0)-1}z_{\overline{n},a-2k}
   \right) 
\\
&=
\prod_{i=1}^l\prod_{k=\beta_i(0)}^{\beta'_i(0)-1}
z_{\overline{n},a-2k} z_{n,a-2k-2}.
\end{aligned}
\end{equation}
Therefore, it is enough to show that
\begin{align}\label{eq:E-weight}
E'=E\delta.
\end{align}
To see it, we notice two facts
following from \eqref{eq:delta-weight}, and
the definitions of $E$, $E'$, and $\phi$:

Fact 1. Both $E'$ and $E\delta$ are products of the
factors $z_{\overline{n},a-2k}z_{n,a-2k+2}$ with $k\in  
\mathbb{Z}$.

Fact 2.  For each $k\in \mathbb{Z}$, the total degree of 
$z_{\overline{n},a-2k}$ and $z_{n,a-2k}$ for $E'$ and the one
for $E\delta$ are equal.

Now, the equality \eqref{eq:E-weight} easily follows from these facts.

\appendix
\section{Proof of Lemma \ref{lem:complementary}}\label{sec:graph}
In this section, we prove Lemma \ref{lem:complementary}. 
To do so, we use a lemma (Lemma \ref{lem:graph-two}) 
for a certain graph. 
 
Let $H\subset \bR^2$ be the upper half plane.  
A {\it graph} $\Gamma$ 
consists of a set $V(\Gamma)$ of vertices 
on the boundary of $H$, denoted by $\times$, 
and a set $A(\Gamma)$ of arcs in $H$ which satisfies the following 
conditions: 
\begin{enumerate}[(i)]
\item \label{item:def-one}
An arc in $A(\Gamma)$ connects two vertices in $V(\Gamma)$.
\item \label{item:def-two}
In each side of a vertex $v$, 
the number of the vertices which are connected with $v$ is at most one. 
\item \label{item:def-three}
The arcs are nonintersecting with each other at $H\backslash V(\Gamma)$. 
\end{enumerate}
The vertices are labeled as $L$, $R$, $L$, $R$, \dots 
{}from left to right. The following 
is an example of a graph: 
\begin{equation}\label{eq:graph}
\begin{psfrags}
\psfrag{L}{$\scriptstyle L$}
\psfrag{R}{$\scriptstyle R$}
\includegraphics[width=11cm, clip]{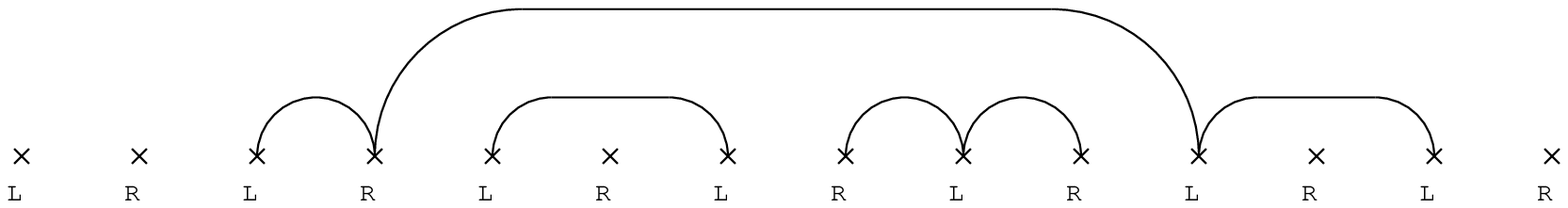}
\end{psfrags}
\end{equation}
Let $v,v' \in V(\Gamma)$. 
If $v$ and $v'$ are connected by an arc in $A(\Gamma)$, 
then we write $v \frown v'$. 
Then $\frown$ generates an equivalence relation $\sim$ in $V(\Gamma)$. 
We call each equivalence class $S$ of $V(\Gamma)$ a 
{\it segment} of $\Gamma$. 
We say that a segment $S$ is {\it even} (resp.~{\it odd})
if the number of vertices $|S|$ is even (resp.~odd). 

Let $S$ and $S'$ be segments of $\Gamma$, and let $v_1$ and $v_2$ be 
the leftmost and the rightmost vertices in $S$. 
If all the vertices $v \in S$ are between $v_1$ and $v_2$, 
then we say that $S'$ is {\it inside} $S$. 
Let us write $S'\lhd S$ if $S'$ is inside $S$. 

\begin{lem}\label{lem:graph}
A graph $\Gamma$ does not have any odd segment if and only if 
each segment $S$ of $\Gamma$ satisfies one of the following conditions: 
\begin{enumerate}
\item \label{item:graph-one}
The leftmost vertex is $L$ and the rightmost vertex is $R$. 
\item \label{item:graph-two}
The leftmost vertex is $R$ and the rightmost vertex is $L$. 
\end{enumerate}
\end{lem}

\begin{proof}
($\Rightarrow$) 
Let $S$ be an even segment of $\Gamma$ 
and let $v_1$ and $v_2$ be its leftmost and the rightmost 
vertices. 
We can show, by induction with respect to $\lhd$, 
that the number $h(S)$ of the vertices between 
$v_1$ and $v_2$ is even, which implies that 
$S$ satisfies \eqref{item:graph-one} or \eqref{item:graph-two}. 

($\Leftarrow$) If $\Gamma$ has an odd segment, 
then there exists an odd segment $S$ which does not have 
any odd segment inside. 
Since the number $h(S)$ is odd, 
$S$ satisfies neither 
\eqref{item:graph-one} nor \eqref{item:graph-two}. 
\end{proof}

For any graph $\Gamma$, we define  
the {\it dual graph} $\Gamma^*=(V(\Gamma^*), A(\Gamma^*))$ 
as follows: 
Let $V(\Gamma^*)$ be the set of all the vertices $\diamond$ which are 
placed between each nearest pair of vertices in $V(\Gamma)$, 
and the vertex `$\infty$'
which is on the right to the rightmost vertex of $\Gamma$. 
Let $A(\Gamma^*)$ be the set of all the dotted arcs 
which connects each nearest pair of vertices in $V(\Gamma^*)$ 
without intersecting with the arcs in $A(\Gamma^*)$. 
For example, the dual graph $\Gamma^*$ of $\Gamma$ in \eqref{eq:graph} 
is given as follows: 
\begin{equation*}
\begin{psfrags}
\psfrag{L}{$\scriptstyle L$}
\psfrag{R}{$\scriptstyle R$}
\psfrag{infty}{$\scriptstyle \infty$}
\includegraphics[width=11cm, clip]{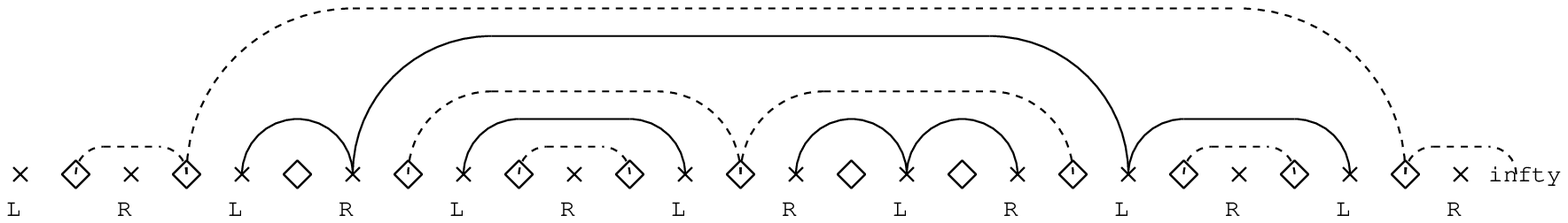}
\end{psfrags}
\end{equation*}
A {\it segment} $S^*$ of $\Gamma^*$
is similarly defined as that of $\Gamma$. 
We say the vertices of $\Gamma^*$ 
between $L$ and $R$ (resp.~$R$ and $L$) 
is of {\it LR-type} 
(resp.~of {\it RL-type}). 

\begin{defn}Let $S^*$ be a segment of $\Gamma^*$. Then, 
\begin{enumerate}
\item 
We say that $S^*$ is {\it mixed} if it contains  vertices 
of both LR-type and RL-type simultaneously, 
and {\it unmixed} otherwise. 
We say an unmixed segment $S^*$ is 
{\it LR-type} (resp.~{\it RL-type}) if 
it consists of vertices of LR-type (resp.~RL-type) 
(also, $\infty$ may be included). 
\item 
We say that $S^*$ is {\it bounded} if it does not contain $\infty$. 
\end{enumerate}
\end{defn}

\begin{lem}\label{lem:graph-two}
A graph $\Gamma$ does not have any odd segment 
if and only if the following conditions are satisfied: 
\begin{enumerate}[(i)]
\item \label{item:grpha-two1} $|V(\Gamma)|$ is even. 
\item \label{item:grpha-two2} $\Gamma^*$ is unmixed. 
\item \label{item:grpha-two3} Any segment of $LR$-type of $\Gamma^*$ is bounded.
\end{enumerate}
\end{lem}

\begin{proof}
($\Rightarrow$) 
\eqref{item:grpha-two2}
If there exists a mixed segment of $\Gamma^*$, 
then one of the following occurs, 
which contradicts to Lemma \ref{lem:graph}: 
\begin{center}
\begin{psfrags}
\psfrag{L}{$\scriptstyle L$}
\psfrag{R}{$\scriptstyle R$}
\psfrag{dots}{\dots}
\includegraphics[width=5.5cm, clip]{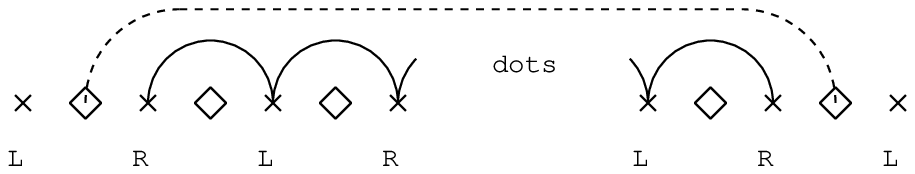}
\end{psfrags}
\qquad  
\qquad  
\begin{psfrags}
\psfrag{L}{$\scriptstyle R$}
\psfrag{R}{$\scriptstyle L$}
\psfrag{dots}{\dots}
\includegraphics[width=5cm, clip]{graph-3.eps}
\end{psfrags}
\end{center}
\eqref{item:grpha-two3} If there exits an unmixed segment of LR-type 
which is not bounded, 
then there exists an LR-type vertex 
$w\in V(\Gamma^*)$ which belongs to the same segment with $\infty$. 
\begin{equation*}
\begin{psfrags}
\psfrag{L}{$\scriptstyle L$}
\psfrag{R}{$\scriptstyle R$}
\psfrag{dots}{\dots}
\psfrag{infty}{$\scriptstyle \infty$}
\psfrag{v}{$\scriptstyle w$}
\includegraphics[width=7cm, clip]{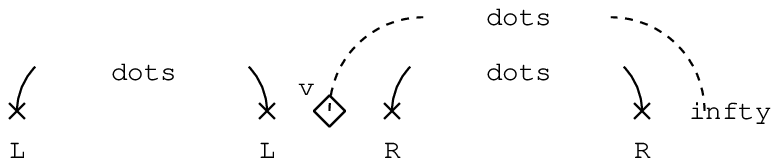}
\end{psfrags}
\end{equation*}
The number of all the vertices to the left of $w$ 
is odd, 
which means that there exists an odd segment 
of $\Gamma$ to the left of $w$. 

($\Leftarrow$)
If there exists an odd segment in $\Gamma$, 
then there exists a segment $S$ of $\Gamma$ 
that satisfies one of the 
following conditions, by Lemma \ref{lem:graph}: 
\begin{enumerate}
\item \label{item:graph-1}
Both the leftmost and the rightmost vertices of $S$ are $L$. 
\item \label{item:graph-2}
Both the leftmost and the rightmost vertices of $S$ are $R$. 
\end{enumerate}
For example, suppose that $S$ satisfies \eqref{item:graph-1}, 
and let $v \in V(\Gamma^*)$ be the one right-next to $S$. 
Then, one of the following occurs (by \eqref{item:grpha-two1}): 
\begin{equation*}
\begin{psfrags}
\psfrag{L}{$\scriptstyle R$}
\psfrag{R}{$\scriptstyle L$}
\psfrag{dots}{\dots}
\psfrag{v}{$\scriptstyle v'$}
\psfrag{v'}{$\scriptstyle v$}
\includegraphics[width=5.3cm, clip]{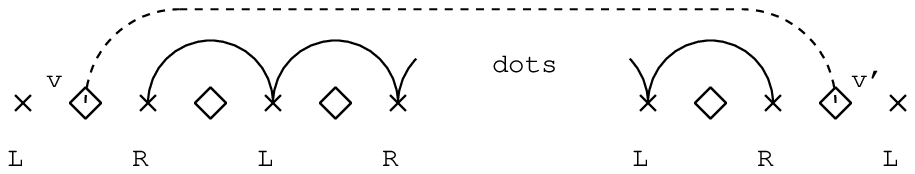}
\end{psfrags}
\qquad  \qquad 
\begin{psfrags}
\psfrag{L}{$\scriptstyle R$}
\psfrag{R}{$\scriptstyle L$}
\psfrag{dots}{\dots}
\psfrag{infty}{$\scriptstyle \infty$}
\psfrag{v'}{$\scriptstyle v$}
\includegraphics[width=6cm, clip]{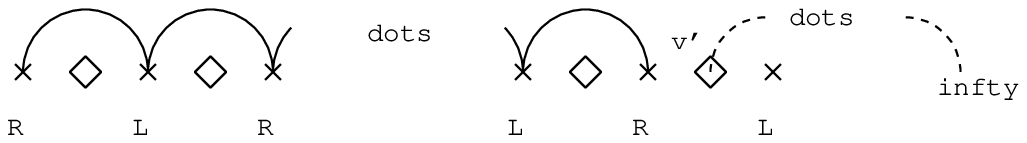}
\end{psfrags}
\end{equation*}
Namely, if there exists an RL-type vertex $v' \in V(\Gamma^*)$ 
left-next to $S$, then $v$ and $v'$ belongs to the same 
segment of $\Gamma^*$, otherwise, $v$ belongs to 
the same segment with the vertex $\infty$. 
The \eqref{item:graph-2} case is similar. 
\end{proof}

Now, let us prove 
Lemma \ref{lem:complementary}. 
With any $(\alpha; \beta) \in \cH(\lambda/\mu)$, 
we associate a graph $\Gamma$ as follows: 
Each vertex of $\Gamma$ naturally corresponds to 
each even $(k-1)$-overlaps of $(\alpha; \beta)$. 
An arc of $\Gamma$ connects a nearest pair of 
even $(k-1)$-overlaps (under the above correspondence) 
belonging to the same $\I_{k-1}$-region. 
Then, an odd segment of $\Gamma$ corresponds to 
an odd $\I_{k-1}$-region of $(\alpha; \beta)$. 
Furthermore, Conditions \eqref{item:complementary-one}, 
\eqref{item:2t-one}, and \eqref{item:c-two}
of Lemma \ref{lem:graph-two} are equivalent to the 
ones of Lemma \ref{lem:complementary} due to the complementarity 
of the $I_{k-1}$- and the $\II_k$-units 
(Lemma \ref{lem:unit} \eqref{item:complementary}). 
This completes the proof of Lemma \ref{lem:complementary}.

\end{document}